# RENORMALIZATION OF THE TWO-DIMENSIONAL LOTKA–VOLTERRA MODEL

By J. Theodore Cox[1] and Edwin A. Perkins[2]

*Syracuse University and University of British Columbia*

We show that renormalized two-dimensional Lotka–Volterra models near criticality converge to a super-Brownian motion. This is used to establish long-term survival of a rare type for a range of parameter values near the voter model.

**1. Introduction.** We consider here the two-dimensional version of a spatially explicit, stochastic Lotka–Volterra model for competition introduced by Neuhauser and Pacala in [7]. The idea there was to formulate and study a model for use in plant ecology which was based on individual, stochastic short range interactions between plants. Most classical competition models are "mean field" differential equations models, and do not take into account the spatial locations of individual plants or individual dynamical effects. Neuhauser and Pacala proved that their model differs in interesting ways from the classical differential equations models. We refer the reader to [7] for a discussion of the biological significance of their findings.

In our previous papers [2] and [3] we studied this model in dimensions $d \geq 3$. In [2] we proved that suitably rescaled nearly critical sequences of these processes converge to super-Brownian motion. In [3] we used this convergence and a renormalization argument to prove that survival and/or coexistence hold for certain parameter regions. In fact, these results were proved for a more general class of models we called *voter model perturbations*.

Our goal here is to extend this work to the more biologically relevant case $d = 2$, the critical dimension. The fact that the two-dimensional random walk is recurrent requires that we use a different mass normalization than in the $d \geq 3$ case, and this complicates the analysis considerably. We believe

Received January 2007; revised September 2007.
[1]Supported in part by NSF Grant DMS-05-05439. Part of the research was done while the author was visiting the University of British Columbia.
[2]Supported in part by an NSERC Discovery Grant.

that appropriate versions of our main results, Theorems 1.2 and 1.3 below, hold for the general voter model perturbations of [2], but in order to keep the presentation as simple as possible, we will consider only the Lotka–Volterra models which we now define.

The state space for our process is $\{0,1\}^{\mathbb{Z}^d}$, where $\mathbb{Z}^d$ is the $d$-dimensional integer lattice. For $\xi \in \{0,1\}^{\mathbb{Z}^d}$, we interpret $\xi(x) = i$ to mean there is a plant of type $i$ at site $x \in \mathbb{Z}^d$, $i = 0, 1$, and we will sometimes identify $\xi$ with the set $\{x \in \mathbb{Z}^d : \xi(x) = 1\}$. The parameters for our process are two nonnegative numbers $\alpha_0, \alpha_1$, and a probability mass function $p : \mathbb{Z}^d \to [0,1]$ which satisfies $p(0) = 0$, $p$ is symmetric with covariance matrix $\sigma^2 I$, and the kernel $p(x,y) = p(y - x)$ is irreducible.

Define the *local densities* $f_i = f_i(\xi) = (f_i(x,\xi), x \in \mathbb{Z}^d)$,

$$(1.1) \quad f_i(x,\xi) = \sum_y p(y - x) 1\{\xi(y) = i\}, \qquad i = 0, 1, x \in \mathbb{Z}^d, \xi \in \{0,1\}^{\mathbb{Z}^d},$$

and the Lotka–Volterra *rate function* $c(x,\xi)$ by

$$(1.2) \qquad c(x,\xi) = \begin{cases} f_1(f_0 + \alpha_0 f_1), & \text{if } \xi(x) = 0, \\ f_0(f_1 + \alpha_1 f_0), & \text{if } \xi(x) = 1. \end{cases}$$

The Lotka–Volterra process $\xi_t$ is the unique $\{0,1\}^{\mathbb{Z}^d}$-valued Feller process with rate function $c(x,\xi)$, meaning that the generator of $\xi_t$ is the closure of the operator $\Omega$:

$$\Omega \phi(\xi) = \sum_x c(x,\xi)(\phi(\xi^x) - \phi(\xi))$$

on the set of functions $\phi : \{0,1\}^{\mathbb{Z}^d} \to \mathbb{R}$ depending on only finitely many coordinates (see, e.g., Remark 2.5 of [3]). Here $\xi^x(y) = \xi(y)$ for $y \neq x$ and $\xi^x(x) = 1 - \xi(x)$.

One can interpret the rate function in the following way. A plant of type $i$ at site $x$ in configuration $\xi$ dies at rate $f_i + \alpha_i f_{1-i}$ and is immediately replaced by a plant of type $\xi(y)$, where $y$ is chosen with probability $p(y - x)$. The death rate incorporates both interspecific and intraspecfic effects. The parameter $\alpha_i$ measures the competitive effect of the neighboring type $1 - i$ plants on type $i$, while we set the self-competition parameter equal to one.

Since $f_0 + f_1 = 1$, $c(x,\xi)$ can also be written in the form

$$(1.3) \qquad c(x,\xi) = \begin{cases} f_1 + (\alpha_0 - 1)f_1^2, & \text{if } \xi(x) = 0, \\ f_0 + (\alpha_1 - 1)f_0^2, & \text{if } \xi(x) = 1. \end{cases}$$

Setting $\alpha_0 = \alpha_1 = 1$ results in the well-known *voter model* (see Chapter 4 of [6]). In [1] an invariance principle was proved for the voter model. Namely, appropriately rescaled voter models converge to super-Brownian motion. The above form for $c(x,\xi)$ suggests the possibility of a similar result holding for the Lotka–Volterra model for parameters $\alpha_i$ sufficiently close to one. This is the case, as was proved in [2] for high dimensions, $d \geq 3$. We briefly recall the main result of that paper.

We define a sequence of rescaled Lotka–Volterra models as follows. Consider a sequence $\{\alpha_i^N : N \in \mathbb{N}\}$ of $\alpha_i$ values for $i = 1, 2$ and let $\xi_t = \xi_t^{(N)}$ denote the Lotka–Volterra model with $\alpha_i = \alpha_i^N$. For $N = 1, 2, \ldots$, let $\mathbf{S_N} = \mathbb{Z}^d/\sqrt{N}$ and define the kernels $p_N : \mathbf{S}_N \to [0, 1]$ by

$$(1.4) \qquad p_N(x) = p(x\sqrt{N}), \qquad x \in \mathbf{S}_N.$$

For $\xi \in \{0,1\}^{\mathbf{S}_N}$, the rescaled densities $f_i^N = f_i^N(\xi) = (f_i^N(x,\xi), x \in \mathbf{S}_N)$, are given by

$$(1.5) \qquad f_i^N(x,\xi) = \sum_{y \in \mathbf{S}_N} p_N(y-x) 1\{\xi(y) = i\}, \qquad i = 0, 1.$$

Then $\xi_t^N(x) = \xi_{Nt}^{(N)}(x\sqrt{N}), x \in \mathbf{S}_N$, is the unique Feller process taking values in $\{0,1\}^{\mathbf{S}_N}$ with rate function

$$(1.6) \qquad c_N(x,\xi) = \begin{cases} N(f_1^N + (\alpha_0^N - 1)(f_1^N)^2), & \text{if } \xi(x) = 0, \\ N(f_0^N + (\alpha_1^N - 1)(f_0^N)^2), & \text{if } \xi(x) = 1. \end{cases}$$

Given a sequence $N' = N'(N)$, we define the measure-valued processes $X_t^N$ by

$$(1.7) \qquad X_t^N = \frac{1}{N'} \sum_{x \in \mathbf{S}_N} \xi_t^N(x) \delta_x,$$

where $\delta_x$ is the unit point mass at $x$. That is, we place an atom of size $1/N'$ at each site $x$ with $\xi_t^N(x) = 1$. (We will see below that the appropriate choice for $N'$ is dimension dependent.) If $\sum_x \xi_0^N(x) < \infty$, then for each $t \geq 0$, $X_t^N \in \mathcal{M}_f = \mathcal{M}_f(\mathbb{R}^d)$, the space of finite Borel measures on $\mathbb{R}^d$, which we endow with the topology of weak convergence. Let $D([0,\infty), \mathcal{M}_f)$ be the Skorokhod space of cadlag $\mathcal{M}_f$-valued paths, and let $\Omega_{X,C}$ be the space of continuous $\mathcal{M}_f$-valued paths with the topology of uniform convergence on compacts. In either case, $X_t$ will denote the coordinate function, $X_t(\omega) = \omega(t)$. Integration of a function $\phi$ with respect to a measure $\mu$ will be denoted by $\mu(\phi)$. Also, we will use $\mathbf{1}$ to denote the function identically one on $\mathbb{R}^d$.

We make the following assumptions about the initial states $\xi_0^N$:

$$(1.8) \qquad \begin{array}{ll} \text{(a)} & \sum_{x \in \mathbf{S}_N} \xi_0^N(x) < \infty, \\ \text{(b)} & X_0^N \to X_0 \qquad \text{in } \mathcal{M}_f(\mathbb{R}^d) \text{ as } N \to \infty. \end{array}$$

Our basic assumption concerning the rates $\alpha_i^N$ is, for $i = 0, 1$,

$$(1.9) \qquad \theta_i^N \equiv N'(\alpha_i^N - 1) \to \theta_i \in \mathbb{R} \qquad \text{as } N \to \infty.$$

In some places we only require that

(1.10) $$\bar{\theta} = 1 \vee \sup\{|\theta_i^N| : i = 0, 1, N = 1, 2, \ldots\} < \infty.$$

Let $P_N$ be the law of $X_\cdot^N$ on $D([0, \infty), \mathcal{M}_f)$ and for $X_0 \in \mathcal{M}_f$, let $P_{X_0}^{\gamma, \theta, \sigma^2}$ be the law of super-Brownian motion with branching rate $\gamma$, drift $\theta$ and diffusion coefficient $\sigma^2$ on $\Omega_{X,C}$ (and on $D([0, \infty), \mathcal{M}_f)$). (See Section 3 below for a characterization of $P_{X_0}^{\gamma, \theta, \sigma^2}$.)

In [2] it was shown that, for dimension $d \geq 3$ and with $N' \equiv N$, $P_N \Rightarrow P_{X_0}^{\gamma, \theta, \sigma^2}$ as $N \to \infty$ for certain parameters $\gamma$ and $\theta$, where $\Rightarrow$ denotes weak convergence on $D([0, \infty), \mathcal{M}_f)$. To define these parameters, we introduce random walk systems $\{B_t^x, t \geq 0, x \in \mathbb{Z}^d\}$ and $\{\hat{B}_t^x, t \geq 0, x \in \mathbb{Z}^d\}$. The walks $B_t^x$ and $\hat{B}_t^x$ are rate one walks with step distribution $p(\cdot)$ and $B_0^x = \hat{B}_0^x = x$. The system $\{B_t^x, t \geq 0, x \in \mathbb{Z}^d\}$ is a system of independent random walks. The system $\{\hat{B}_t^x, t \geq 0, x \in \mathbb{Z}^d\}$ is a system of coalescing random walks, meaning that the walks move independently of one another until they meet, at which time they coalesce and move together. We define the collision times $\tau(x, y) = \inf\{t \geq 0 : B_t^x = B_t^y\}$ and $\hat{\tau}(x, y) = \inf\{t \geq 0 : \hat{B}_t^x = \hat{B}_t^y\}$, and the constants

$$\gamma_e = \sum_e p(e) P(\hat{\tau}(0, e) = \infty),$$

$$\gamma_0 = \sum_{e, e'} p(e) p(e') P(\hat{\tau}(0, e) = \hat{\tau}(0, e') = \infty),$$

$$\gamma_1 = \sum_{e, e'} p(e) p(e') P(\hat{\tau}(0, e) = \hat{\tau}(0, e') = \infty, \hat{\tau}(e, e') < \infty).$$

In [2] we proved the following.

THEOREM 1.1. *Assume $N' \equiv N$ and $d \geq 3$, and (1.8) and (1.9) hold. Then*

$$P_N \Rightarrow P_{X_0}^{2\gamma_e, \theta, \sigma^2} \qquad as \ N \to \infty,$$

*where $\theta = \theta_0 \gamma_0 - \theta_1 \gamma_1$.*

The strategy used for the proof of Theorem 1.1 in [2] was the following. First, derive a semimartingale representation for $X_\cdot^N$. Second, obtain $L^2$ bounds on $X_t^N(\mathbf{1})$, which, along with the semimartingale representation, lead to a proof of tightness of the laws $P_N$. Finally, show that any limiting martingale problem for $X_\cdot^N$ takes the form of the martingale problem characterizing super-Brownian motion. Unfortunately, there are significant difficulties implementing this approach when $d = 2$.

For the two-dimensional voter model, the appropriate mass renormalization factor is $N' = N/(\log N)$ (see Theorem 1.2 in [1]). The normalization $N' \equiv N$ leads only to deterministic heat flow in the limit, as in [9]. Adopting this $N'$, the first problem we encounter is that the estimates in [1] no longer imply even $L^1$ boundedness of $X_t^N(\mathbf{1})$, let alone $L^2$ boundedness. This means that even tightness is a more complicated issue than it was before. The method used for $d \geq 3$ in [2] depended on the fact that over short time scales $\varepsilon$ the Lotka–Volterra dynamics are quite close to voter dynamics, as the rates of the "perturbation terms" are of smaller order. By conditioning back time $\varepsilon$, we allow the quadratic terms in the drift to relax to voter model equilibria values, and this produces the constants in the limiting super-Brownian motion. We needed to choose $\varepsilon = \varepsilon(N) \to 0$ so that the Lotka–Volterra process is still well approximated by the voter dynamics, but slow enough so that the system has a chance to relax. Here, however, our errors in the voter model approximation to the Lotka–Volterra model are multiplied by a factor of $\log N$. This effectively puts an upper bound on $\varepsilon$ so that the voter approximation over time $\varepsilon$ is good enough. At this point we must verify that this still gives the system enough time to relax to its equilibrium values. The factor of $\log N$ also makes our total mass bounds problematic. Again, the key is the above voter comparison, as once we pass to the voter model over short time intervals, the voter model *clustering* in $d = 2$ effectively absorbs this factor, providing $\varepsilon$ is enough time for the system to cluster.

Another new issue for $d = 2$ is that even to define the parameters of our limiting super-Brownian motion some new two-dimensional random walk estimates are required (see, e.g., Lemma 2.5 and Proposition 2.1 below).

To state our results, we introduce the two-dimensional potential kernel $a(x)$,

$$(1.11) \qquad a(x) = \int_0^\infty [P(B_t^0 = 0) - P(B_t^x = 0)]\, dt.$$

Note that $a(x) \geq 0$ since symmetry of $p(\cdot)$ implies $P(B_t^0 = 0) \geq P(B_t^x = 0)$. We may now define

$$\gamma^* = 2\pi\sigma^2 \int_0^\infty \sum_{x,y,e,e' \in \mathbb{Z}^2} p(e)p(e')$$

$$\times P(\tau(0,e) \wedge \tau(0,e') > \tau(e,e') \in du,$$

$$B_u^0 = x, B_u^e = y)$$

$$\times a(y - x).$$

(1.12)

The fact that $\gamma^*$ is finite is contained in Proposition 2.1, proved in Section 2.

Our two-dimensional Lotka–Volterra invariance principle is the following:

THEOREM 1.2. *Assume $d = 2$, $N' = N/(\log N)$, and (1.8) and (1.9) hold. Then*

$$P_N \Rightarrow P_{X_0}^{4\pi\sigma^2, \theta, \sigma^2} \qquad as\ N \to \infty,$$

*where $\theta = \gamma^*(\theta_0 - \theta_1)$.*

Theorem 1.1 was used in [3] to prove, for dimensions $d \geq 3$, that survival holds for a region of parameter values $\alpha = (\alpha_0, \alpha_1)$ near $(1, 1)$. If $P^\alpha$ denotes the dependence of the Lotka–Volterra model on $\alpha$, survival for parameter values $\alpha$ means that

$$P^\alpha(|\xi_t| > 0 \text{ for all } t > 0 \,|\, |\xi_0| = 1) > 0,$$

where $|\xi| = \sum_x \xi(x)$. A similar result holds here. Let $S$ be the set of all $(\alpha_0, \alpha_1)$ for which survival occurs. For $0 < \eta < 1$, define $S^\eta$ to be the set of all $(\alpha_0, \alpha_1) \neq (1, 1)$ such that

$$\alpha_1 - 1 < \begin{cases} (1 - \eta)(\alpha_0 - 1), & \text{if } \alpha_0 \geq 1, \\ (1 + \eta)(\alpha_0 - 1), & \text{if } \alpha_0 < 1. \end{cases}$$

THEOREM 1.3. *For $0 < \eta < 1$, there exists $r(\eta) > 0$ such that survival holds for all $(\alpha_0, \alpha_1) \in S^\eta$ such that $1 - r(\eta) < \alpha_0$ and $\alpha_1 < 1 + r(\eta)$.*

If $\tilde{S}^\eta = \{(\alpha_0, \alpha_1) \in S^\eta : 1 - r(\eta) < \alpha_0 \text{ and } \alpha_1 < 1 + r(\eta)\}$, Theorem 1.3 shows survival holds on the region $\tilde{S} = \bigcup_{0 < \eta < 1} \tilde{S}^\eta$ illustrated in Figure 1.

Theorem 1.3 follows from Theorem 1.2 by a cavalier interchange of limits. Long-term survival for the limiting super-Brownian motion in Theorem 1.2 occurs iff $\theta > 0$, that is, $\theta_0 > \theta_1$. Interchanging limits as $N \to \infty$ and $t \to \infty$ leads to survival of the original particle system for $\alpha_0 > \alpha_1$ and $\alpha$ near $(1,1)$. The monotonicity of the Lotka–Volterra models (increasing in $\alpha_0$ and decreasing in $\alpha_1$) established in Section 1 of [3] [see (1.3) of that work] allows one to infer survival for larger values of $\alpha_0$ and smaller values of $\alpha_1$ as stated in Theorem 1.3.

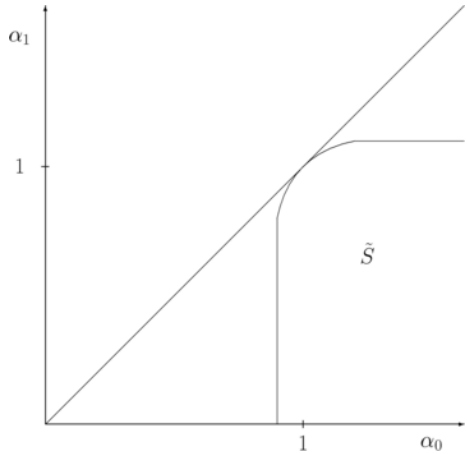

FIG. 1.  *Survival region of Theorem 1.3.*

The above interchange of limits argument is carried out for $d \geq 3$ in Theorem 1 of [3], where now $\theta > 0$ in Theorem 1.1 leads to survival for $|\alpha_0 - 1| < r(\eta)$ satisfying

$$(1.13) \qquad \alpha_1 - 1 < \begin{cases} (m_0 + \eta)(\alpha_0 - 1), & \text{if } \alpha_0 < 1, \\ (m_0 - \eta)(\alpha_0 - 1), & \text{if } \alpha_0 \geq 1, \end{cases}$$

where $m_0 = \gamma_1/\gamma_0 < 1$. In the proof of Theorem 1.3 the interchange of limits is justified by a comparison with super-critical oriented percolation as in [3]. To obtain suitable independence in our percolation events, we must study the Lotka–Volterra model with 0 boundary conditions outside a large box and show the effect of these boundary conditions is small on an appropriate space–time region (Lemma 9.1). This argument is now more involved than the corresponding one in [3] due to the different mass normalization and the fact that the conditioning technique used in the convergence theorem must be adapted to handle this new type of bound.

Note that in Theorem 1.3 survival fails at the point $\alpha = (1, 1)$ itself. We conjecture that Theorem 1.3 is sharp in the following sense:

CONJECTURE. *There is a continuous curve $\alpha_1 = g(\alpha_0)$, tangent to $\alpha_1 = \alpha_0$ at $(1, 1)$, so that survival fails above the curve for $\alpha_0$ close to 1.*

The corresponding conjecture has been proved for $d \geq 3$ to the left of $(1, 1)$ with the slope $m_0$ in place of 1, but in fact is false to the right of $(1, 1)$ (with slope $m_0$) in this higher dimensional setting. The exact state of affairs will appear in a forthcoming article with Rick Durrett.

We say coexistence holds for the Lotka–Volterra model if there is a stationary distribution under which there are both 0's and 1's (necessarily infinitely many of each) a.s. For $d \geq 3$, it follows from (1.13) and symmetry that survival of 0's will occur if

$$(1.14) \qquad \alpha_1 - 1 > (m_0 + \eta)^{-1}(\alpha_0 - 1), \qquad \alpha_0 < 1 \text{ and close to } 1.$$

We have restricted $\alpha_0 < 1$ so that the regions in (1.13) and (1.14) intersect in a local nonempty wedge to the left of $(1, 1)$ containing a local piece of the diagonal $\alpha_1 = \alpha_0$ (see Figure 3 in [3]). The survival of both types for parameter values in this local wedge easily leads to coexistence in this local wedge for $d \geq 3$ (see Theorem 4 of [3]). The wedge is nonempty because $m_0 < 1$. For $d = 2$ in Theorem 1.3, we are in the critical case $m_0 = 1$, and the above proof of coexistence just fails. Neuhauser and Pacala have conjectured coexistence holds along the diagonal $\alpha_0 = \alpha_1$ for $\alpha_0 < 1$ (see Conjecture 1 of [7]) and proved it for small $\alpha_0$ (see Theorem 1 of [7]). Theorem 4 in [3] confirms this for $d \geq 3$ for $\alpha_0$ close to 1, but Theorem 1.3 and the above Conjecture suggest that for $d = 2$ their conjecture is quite delicate near $\alpha_0 = 1$. One approach to establishing coexistence near $\alpha_0 = 1$ would be to derive second order (concave up) asymptotics for the survival region near $(1, 1)$ and so deduce survival in an open thorn with tip at $(1, 1)$ with slope 1.

Section 2 gives an alternative description of $\gamma^*$ and proves it is finite (Proposition 2.1). The martingale problems for both the approximating Lotka–Volterra processes and limiting super-Brownian motion are provided in Section 3. In Section 4 we state the key bounds, including $L^1$ and $L^2$ bounds on the total mass (Proposition 4.3), a mean measure bound (Proposition 4.4), and asymptotics for the increasing processes arising in the martingale problem for $X^N_\cdot$ obtained in Section 3 (Proposition 4.7). We then show how Theorem 1.2 follows from these estimates. Section 5 contains mass and mean measure bounds for the voter and biased voter models. Both upper and lower biased voter models are used because they are simpler to handle than the Lotka–Volterra model, and bound it above and below, respectively. The voter model estimates are obtained by direct duality calculations, and the biased voter bounds are then obtained by conditioning back over short time intervals, arguing that the voter dynamics are close over such intervals. In Section 6 the above bounds and techniques are used to prove Propositions 4.3–4.5. The proof of Proposition 4.7, which is more delicate as the demands on our $\varepsilon$ relaxation time are more severe, is split over Sections 7 and 8. Finally, the proof of Theorem 1.3 is given in Section 9.

NOTE. Henceforth, we will assume that $d = 2$ and $N' = N/(\log N)$.

**2. Characterization of $\gamma^*$.** Recall the independent system of random walks $\{B^x_t, x \in \mathbb{Z}^2\}$ and the coalescing random walk system $\{\hat{B}^x_t, x \in \mathbb{Z}^2\}$ from Section 1, and also the collision times $\tau(x, y)$ and $\hat{\tau}(x, y)$. For $e, e' \in \mathbb{Z}^2$, define the event $\Gamma_T(e, e') = \{\hat{\tau}(e, e') < T, \hat{\tau}(0, e) \wedge \hat{\tau}(0, e') > T\}$, and let

$$q_T = \sum_{e,e' \in \mathbb{Z}^2} p(e)p(e')P(\Gamma_T(e, e')). \tag{2.1}$$

We will need the following characterization of $\gamma^*$ defined in (1.12).

PROPOSITION 2.1.

$$\gamma^* = \lim_{T \to \infty} (\log T)q_T < \infty. \tag{2.2}$$

Before beginning the proof of Proposition 2.1, we assemble some facts about two-dimensional random walk. Let $\tau_x = \inf\{t \geq 0 : B_t^0 = x\}$, and write $P^x$ to indicate the law of the walk $B_\cdot^x$. Let $\tilde{P}(\cdot) = \sum_{e \in \mathbb{Z}^2} p(e) P^e(\cdot)$, and define

$$H(t) = \tilde{P}(\tau_0 > t). \tag{2.3}$$

Let $|x|$ be the Euclidean norm of $x \in \mathbb{R}^2$.

PROPOSITION 2.2.

$$\lim_{t \to \infty} H(t) \log t = 2\pi\sigma^2. \tag{2.4}$$

$$\frac{P^x(\tau_0 > t)}{H(t)} \leq 2a(x) \qquad \text{for all } x \in \mathbb{Z}^2, t > 0. \tag{2.5}$$

$$\lim_{t \to \infty} \frac{P^x(\tau_0 > t)}{H(t)} = a(x) \qquad \text{for all } x \in \mathbb{Z}^2. \tag{2.6}$$

$$a(x)/|x|, \qquad x \neq 0 \text{ is bounded.} \tag{2.7}$$

REMARK 2.3. The trivial bound (2.7) is derived below. In fact, it is not hard to show $a(x)/\log|x| \to 1/\pi\sigma^2$ as $|x| \to \infty$. If $\sum_e |e|^r p(e) < \infty$ for some $r > 2$, this limit is a simple consequence of the far more precise P12.3 in [10].

PROOF OF PROPOSITION 2.2. The limit (2.4) is well known, for example, see Lemma A.3 in [1] for a proof. Now let $Y_n, n = 0, 1, 2, \ldots$, be a random walk with step distribution $p(\cdot)$, and let $\sigma_x = \inf\{n \geq 1 : Y_n = x\}$. We will abuse notation slightly and also let $P^x$ denote the law of the walk starting at $Y_0 = x$. We note that $a(x)$ defined in (1.11) is also given by $\sum_{n=0}^{\infty} [P^0(Y_n = 0) - P^0(Y_n = x)]$.

By P11.5 in [10], $P^0(\sigma_x < \sigma_0) = 1/2a(x)$. Since the sequence of states visited by the walk $B_t^0$ is equal in law to the sequence visited by the walk $Y_n$ (with $Y_0 = 0$), it follows that $\tilde{P}(\tau_x < \tau_0) = 1/2a(x)$. By the strong Markov property,

$$H(t) \geq \sum_{e \in \mathbb{Z}^2} p(e) P^e(\tau_x < \tau_0 \text{ and } \tau_0 > t)$$

$$\geq \sum_{e \in \mathbb{Z}^2} p(e) P^e(\tau_x < \tau_0) P^x(\tau_0 > t),$$

and (2.5) follows.

For (2.6), we recall the (discrete time) result T16.1 of [10],

$$\lim_{n \to \infty} \frac{P^x(\sigma_0 > n)}{P^0(\sigma_0 > n)} = a(x). \tag{2.8}$$

A standard result (see the local limit theorem P7.9 of [10]) is

(2.9) $$(\log n)P^0(\sigma_0 > n) \to 2\pi\sigma^2 \quad \text{as } n \to \infty.$$

In view of (2.4), $H(t)/P^0(\sigma_0 > t) \to 1$ as $t \to \infty$, and therefore, in order to prove (2.6), it suffices to show that

(2.10) $$\lim_{t\to\infty} \frac{P^x(\tau_0 > t)}{P^x(\sigma_0 > t)} = 1.$$

To do this, let $S(t), t \geq 0$, be a rate one Poisson process, independent of the walk $Y_n$. Then $Y_{S(\cdot)}$ is a realization of $B_\cdot^0$. By a standard large deviations estimate, there is a constant $C > 0$ such that $e^{Ct}P(S(t) \notin [t/2, 2t]) \to 0$ as $t \to \infty$. Consequently,

$$P^x(\tau_0 > t)$$
$$= o(e^{-Ct}) + \sum_{k \in [t/2, 2t]} \frac{e^{-t}t^k}{k!} P^x(\sigma_0 > k)$$
$$\leq o(e^{-Ct}) + P^x(\sigma_0 > t/2)$$

as $t \to \infty$, and similarly,

$$P^x(\tau_0 > t) \geq (1 - o(e^{-Ct}))P^x(\sigma_0 > 2t).$$

These inequalities, together with (2.8) and (2.9), imply (2.10).

Finally, let $\psi(\theta)$ be the characteristic function of $Y_1$, $\psi(\theta) = \sum_x e^{ix\cdot\theta}p(x)$, $\theta \in \mathbb{R}^2$. Then (see equation (12.3) in [10])

$$a(x) = (2\pi)^{-2} \int_{[-\pi,\pi]^2} \frac{1 - e^{ix\cdot\theta}}{1 - \psi(\theta)} d\theta.$$

Since $|1 - e^{ix\cdot\theta}| \leq |x||\theta|$ and $|\theta|/|1 - \psi(\theta)|$ is integrable over $[-\pi,\pi]^2$ (see the paragraph after (12.3) in [10]), (2.7) holds. □

The next result is a general inequality which we will apply to bound the probability that walks starting from $0, e, e'$ avoid each other for at least time $T$.

LEMMA 2.4. *Let $X, Y$ be nonnegative random variables. For any constant $c > 0$, $E(XY \wedge c) \leq \sqrt{c}E(X + Y)$.*

PROOF. For any $a > 0$,
$$E(XY \wedge c) = E(XY \wedge c; Y \leq a) + E(XY \wedge c; Y > a)$$
$$\leq aE(X) + cP(Y > a)$$
$$\leq aE(X) + cE(Y)/a.$$

Now let $a = \sqrt{c}$. □

LEMMA 2.5. *For distinct sites $0, e_1, e_2$, let $\tau = \tau(0, e_1) \wedge \tau(0, e_2) \wedge \tau(e_1, e_2)$. Then for $T > 0$,*

$$P(\tau > T) \leq (2H(2T))^{3/2} (a(e_1) + a(e_2)) \sqrt{a(e_1 - e_2)}.$$

PROOF. Define

$$\chi_T = P(\tau > T \mid B_s^0, 0 \leq s \leq T)$$

and

$$\chi_T^i = P(\tau(0, e_i) > T \mid B_s^0, 0 \leq s \leq T), \qquad i = 1, 2.$$

Note that (2.5) implies that, for $i = 1, 2$, $E(\chi_T^i) \leq 2a(e_i)H(2T)$ [since $B_t^0 - B_t^{e_i}$ is a rate-two random walk with step distribution $p(\cdot)$]. By independence,

$$P(\tau(0, e_1) > T, \tau(0, e_2) > T \mid B_s^0, s \leq T) = \chi_T^1 \chi_T^2,$$

and also

$$P(\tau(e_1, e_2) > T \mid B_s^0, s \leq T) = P(\tau(e_1, e_2) > T) \leq 2a(e_1 - e_2)H(2T).$$

Consequently, since each $\chi_T^i \leq 1$,

$$\chi_T \leq \chi_T^1 \chi_T^2 \wedge 2a(e_1 - e_2)H(2T).$$

We may now apply Lemma 2.4 to complete the proof. □

PROOF OF PROPOSITION 2.1. Since $\Gamma_T(e, e') \subset \{\hat{\tau}(0, e) > T\}$, summation over $e'$ gives $q_T \leq \sum_{e \in \mathbb{Z}^2} p(e) P(\hat{\tau}(0, e) > T) = \sum_{e \in \mathbb{Z}^2} p(e) P(\tau(0, e) > T) = H(2T)$. It follows from (2.4) that $q_T(\log T)$ is bounded as $T \to \infty$.

By the Markov property,

$$q_T = \int_0^T \sum_{x,y,e,e' \in \mathbb{Z}^2} p(e)p(e') P(\hat{\tau}(0, e) \wedge \hat{\tau}(0, e') > \hat{\tau}(e, e') \in du,$$

(2.11)
$$\hat{B}_u^0 = x, \hat{B}_u^e = y)$$

$$\times P(\hat{\tau}(x, y) > T - u).$$

Below we will use the fact that this expression is unchanged if we replace the coalescing walks and corresponding hitting times $\hat{B}$ and $\hat{\tau}$ with the independent walks $B$ and $\tau$. Since

$$(\log T) P(\tau(x, y) > T - u) \to 2\pi\sigma^2 a(x - y) \qquad \text{as } T \to \infty \text{ for fixed } u,$$

by (2.4) and (2.6), it follows from Fatou that $\liminf_{T \to \infty} q_T(\log T) \geq \gamma^*$, and therefore, $\gamma^*$ is finite.

The next step is to apply Lemma 2.5 and see that we may restrict the integral in (2.11) to $u < T/2$. By Lemma 2.5,

$$\log T \sum_{e,e' \in \mathbb{Z}^2} p(e)p(e')P(\tau(0,e) \wedge \tau(0,e') \wedge \tau(e,e') > T/2)$$

$$\leq (\log T)(2H(T))^{3/2} \sum_{e,e' \in \mathbb{Z}^2} p(e)p(e')(a(e) + a(e'))\sqrt{a(e - e')}.$$

By (2.7), the sum above is finite, and so by (2.4), the right-hand side above tends to 0 as $T \to \infty$.

Now for $u \leq T/2$, (2.4) and (2.6) imply

$$(\log T)P(\tau(x,y) > T - u) \to 2\pi\sigma^2 a(x - y) \quad \text{as } T \to \infty,$$

and (2.5) implies

$$(\log T)P(\tau(x,y) > T - u) \leq 2a(x - y)(\log T)H(T/2).$$

The right-hand side above is no larger than a fixed multiple of $a(x - y)$ by (2.4). Since we have already established the integral defining $\gamma^*$ is finite, we may apply the dominated convergence theorem to conclude that

$$\log T \int_0^{T/2} \sum_{x,y,e,e' \in \mathbb{Z}^2} p(e)p(e')P(\tau(0,e) \wedge \tau(0,e') > \tau(e,e') \in du,$$

(2.12)
$$B_u^0 = x, B_u^e = y)$$
$$\times P(\tau(x,y) > T - u) \to \gamma^*,$$

which completes the proof of Proposition 2.1. □

**3. Semimartingale decompositions.** Let $\xi_t^N$ be the rescaled Lotka–Volterra model as defined in Section 1. Following [2], we introduce the following notation. If

$$\psi \in C_b(\mathbf{S}_N), \qquad \phi = \phi_s(x), \qquad \dot{\phi}_s(x) \equiv \frac{\partial}{\partial s}\phi(s,x) \in C_b([0,T] \times \mathbf{S}_N),$$

and $s \leq T$, define

(3.1) $\quad \mathcal{A}_N(\psi)(x) = \sum_{y \in \mathbf{S}_N} Np_N(y - x)(\psi(y) - \psi(x)),$

(3.2) $\quad D_t^{N,1}(\phi) = \int_0^t X_s^N(\mathcal{A}_N\phi_s + \dot{\phi}_s)\,ds,$

(3.3) $\quad D_t^{N,2}(\phi) = \frac{N(\alpha_0 - 1)}{N'} \int_0^t \sum_{x \in \mathbf{S}_N} \phi_s(x)\mathbf{1}\{\xi_s^N(x) = 0\}(f_1^N(x,\xi_s^N))^2\,ds,$

(3.4) $$D_t^{N,3}(\phi) = \frac{N(\alpha_1 - 1)}{N'} \int_0^t \sum_{x \in \mathbf{S}_N} \phi_s(x) 1\{\xi_s^N(x) = 1\}(f_0^N(x, \xi_s^N))^2 \, ds,$$

(3.5) $$\langle M^N(\phi) \rangle_{1,t} = \frac{N}{(N')^2} \int_0^t \sum_{x \in \mathbf{S}_N} \phi_s^2(x) \\ \times \sum_{y \in \mathbf{S}_N} p_N(y - x)(\xi_s^N(y) - \xi_s^N(x))^2 \, ds,$$

(3.6) $$\langle M^N(\phi) \rangle_{2,t} = \frac{1}{(N')^2} \int_0^t \sum_{x \in \mathbf{S}_N} \phi_s^2(x) \\ \times [(\alpha_0^N - 1)1\{\xi_s^N(x) = 0\}(f_1^N(x, \xi_s^N))^2 \\ + (\alpha_1^N - 1)1\{\xi_s^N(x) = 1\}(f_0^N(x, \xi_s^N))^2] \, ds.$$

Note that $\langle M^N(\phi) \rangle_{2,t}$ may be negative. If $X$. is a process, let $(\mathcal{F}_t^X, t \geq 0)$ be the right-continuous filtration generated by $X$..

The following result is Proposition 2.3 of [2], which was stated in the case $N' \equiv N$. Only trivial modifications of the proof of that result are necessary to prove the following for general, and hence, our choice of, $N'$. Recall that $X_\cdot^N$ is as in (1.7).

PROPOSITION 3.1.  *For $\phi, \dot\phi \in C_b([0,T] \times \mathbf{S}_N)$ and $t \in [0,T]$,*

(3.7) $$X_t^N(\phi_t) = X_0^N(\phi_0) + D_t^N(\phi) + M_t^N(\phi),$$

*where*

(3.8) $$D_t^N(\phi) = D_t^{N,1}(\phi) + D_t^{N,2}(\phi) - D_t^{N,3}(\phi)$$

*and $M_t^N(\phi)$ is an $(\mathcal{F}_t^{X^N})$ square-integrable martingale with predictable square function*

(3.9) $$\langle M^N(\phi) \rangle_t = \langle M^N(\phi) \rangle_{1,t} + \langle M^N(\phi) \rangle_{2,t}.$$

Having described our approximating martingale problems, it is a good time to recall the target martingale problem. Let $C_b^\infty(\mathbb{R}^d)$ denote the space of bounded infinitely differentiable functions on $\mathbb{R}^d$ with uniformly bounded partial derivatives. An adapted a.s.-continuous $M_f(\mathbb{R}^d)$-valued process $X_t, t \geq 0$, on a complete filtered probability space $(\Omega, \mathcal{F}, \mathcal{F}_t, P)$ is said to be a *super-Brownian motion with branching rate $\gamma \geq 0$, drift $\theta \in \mathbb{R}$ and diffusion coefficient $\sigma^2 > 0$ starting at $X_0 \in \mathcal{M}_f(\mathbb{R}^d)$* iff it solves the following martingale problem:

(MP) For all $\phi \in C_b^\infty(\mathbb{R}^d)$,

$$(3.10) \quad M_t(\phi) = X_t(\phi) - X_0(\phi) - \int_0^t X_s\left(\frac{\sigma^2 \Delta \phi}{2}\right) ds - \theta \int_0^t X_s(\phi) \, ds$$

is a continuous $(\mathcal{F}_t^X)$-martingale, with $M_0(\phi) = 0$ and predictable square function

$$(3.11) \quad \langle M(\phi)\rangle_t = \int_0^t X_s(\gamma \phi^2) \, ds.$$

The existence and uniqueness in law of a solution to this martingale problem is well known (see, e.g., Theorem II.5.1 and Remark II.5.13 of [8]). Recall from Section 1 that $P_{X_0}^{\gamma,\theta,\sigma^2}$ denotes the law of the solution on $\Omega_{X,C}$ (and also on $D([0,\infty), \mathcal{M}_f))$.

**4. Convergence to super-Brownian motion.** The goal of this section is to outline the proofs of the following two key results. Recall that if $S$ is a metric space and $\{Q_N\}$ is a sequence of probabilities on $D([0,\infty), S)$, then $\{Q_N\}$ is $C$-tight iff it is tight and every limit point is supported by $C([0,\infty), S)$.

PROPOSITION 4.1. *The family of laws $\{P_N, N \in \mathbb{N}\}$ is $C$-tight.*

PROPOSITION 4.2. *If $\hat{P}$ is any weak limit point of the sequence $P_N$, then $\hat{P} = P^{4\pi\sigma^2, \theta, \sigma^2}$.*

Clearly, Theorem 1.2 follows from these propositions. We say "outline" because we will state the key technical ingredients as Propositions 4.3–4.7 below, and assuming their validity, derive Propositions 4.1 and 4.2. Propositions 4.3–4.7 will then be established in Sections 6–8.

Proposition 4.3 gives "total mass" bounds and will be used to prove that, for $\phi \in C_b^{1,3}(\mathbb{R}_+ \times \mathbb{R}^3)$, each of the families $X_\cdot^N(\phi)$, $D_\cdot^{N,i}(\phi)$, $M_\cdot^N(\phi)$ and $\langle M^N(\phi)\rangle_\cdot$, $N = 1, 2, \ldots$, is $C$-tight (see Proposition 4.11 below). Propositions 4.4 and 4.5 provide "spatial mass" bounds which are used to prove the required compact containment condition (Proposition 4.12 below). Proposition 4.7 will be used to identify the weak limit points of the $P_N$.

Let

$$(4.1) \quad g(s) = C_{4.1} s^{-1/3} e^{C_{4.1} s},$$

where $C_{4.1}$ is a positive constant which will be chosen in Section 6.

PROPOSITION 4.3. (a) *For $T > 0$, there is a constant $C_{4.2}(T)$ such that, for all $N \in \mathbb{N}$,*

$$(4.2) \quad \sup_{t \leq T} E(X_t^N(\mathbf{1})) \leq C_{4.2}(T) X_0^N(\mathbf{1}),$$

$$(4.3) \quad E\left(\sup_{t \leq T} X_t^N(\mathbf{1})^2\right) \leq C_{4.2}(T)(X_0^N(\mathbf{1})^2 + X_0^N(\mathbf{1})).$$

(b) *For all $s > 0$ and $N \in \mathbb{N}$,*

$$(\log N) E(X_s^N(f_0^N(\cdot, \xi_s^N))) \leq g(s) X_0^N(\mathbf{1}), \tag{4.4}$$

$$(\log N) E(X_s^N(\mathbf{1}) X_s^N(f_0^N(\cdot, \xi_s^N))) \leq g(s)(X_0^N(\mathbf{1})^2 + X_0^N(\mathbf{1})). \tag{4.5}$$

To state the next bounds, we need some additional notation. For $D \subset \mathbb{R}^2$ and $\phi : D \to \mathbb{R}$, define

$$|\phi|_{\mathrm{Lip}} = \sup\left\{ \frac{|\phi(x) - \phi(y)|}{|x - y|} : x \neq y \in D \right\} \quad \text{and} \tag{4.6}$$

$$\|\phi\|_{\mathrm{Lip}} = |\phi|_{\mathrm{Lip}} + \|\phi\|_\infty.$$

Let $\mathcal{A}_N^*(\psi) = (N + \bar{\theta} \log N) \frac{A_N}{N}(\psi)$, with semigroup $P_t^{N,*}$.

PROPOSITION 4.4. *For $p \geq 3$, there is a $C_{4.7}(p)$ such that, for any $t \geq 0$ and $\phi : \mathbb{R}^2 \to \mathbb{R}^+$,*

$$E(X_t^N(\phi)) \leq e^{(\log N)^{1-p}} e^{C_{4.7} t} X_0^N(P_t^{N,*} \phi) \tag{4.7}$$
$$+ C_{4.7} e^{C_{4.7} t} \|\phi\|_{\mathrm{Lip}} (\log N)^{(1-p)/2} X_0^N(\mathbf{1}).$$

PROPOSITION 4.5. *For $p \geq 3$, there is a constant $C_{4.8}(p)$ such that, for all $\phi : \mathbb{R}^2 \to \mathbb{R}_+$, if $\epsilon = (\log N)^{-p}$, then*

$$E(X_\epsilon^N(\log N \phi f_0^N(\cdot, \xi_\varepsilon^N))) \leq C_{4.8} X_0^N(\mathbf{1}) \|\phi\|_{\mathrm{Lip}} (\log N)^{(1-p)/2} \tag{4.8}$$
$$+ C_{4.8} X_0^N(\phi).$$

REMARK 4.6. It will be immediate from the proofs that the constants in Propositions 4.3, 4.4 and 4.5 will depend only on the kernel $p(\cdot)$ and $\bar{\theta}$ (as well as the choice of parameter $p$ in the latter two). This will be convenient in Section 9.

PROPOSITION 4.7. *Let $\sup_{K,T}$ indicate a supremum over all $X_0^N \in \mathcal{M}_f(\mathbf{S}_N)$, $\phi : \mathbb{R}^2 \to \mathbb{R}$ and $t \geq 0$ satisfying $X_0^N(\mathbf{1}) \leq K$, $\|\phi\|_{\mathrm{Lip}} \leq K$ and $t \leq T$. Then for every $K, T > 0$ and $0 < p < 2$:*

(a) $\lim_{N \to \infty} \sup_{K,T} E(|\int_0^t (X_s^N(\log N \phi^2 f_0^N(\cdot, \xi_s^N)) - 2\pi \sigma^2 X_s^N(\phi^2)) \, ds|^p) = 0$;

(b) *for $i = 2, 3$,* $\lim_{N \to \infty} \sup_{K,T} E(|D_t^{N,i}(\phi) - \int_0^t \theta_{i-2} \gamma^* X_s^N(\phi) \, ds|^p) = 0$.

Assuming the validity of these results, we can now prove tightness for the sequence of laws $\{P_N : N = 1, 2, \ldots\}$.

The first step is to obtain more precise information on the terms in the decomposition of $X_t^N(\phi)$ given in Proposition 3.1. Lemma 4.8 below bounds the terms in the increasing process $\langle M^N(\phi) \rangle_t$ and some of the terms in the drift $D_t^N(\phi)$.

LEMMA 4.8.　*There is a constant $C_{4.9}$ such that if $\phi:[0,T] \times \mathbf{S}_N \to \mathbb{R}$ is a bounded measurable function, then for $0 \le t \le T$:*

(a) $\langle M^N(\phi) \rangle_{2,t} = \int_0^t m_{2,s}^N(\phi)\,ds$ *where*

$$(4.9) \qquad |m_{2,s}^N(\phi)| \le C_{4.9} \frac{\|\phi_s\|_\infty^2}{(N')^2} X_s^N(\mathbf{1}), \qquad 0 \le s \le t.$$

(b)

$$(4.10) \quad \langle M^N(\phi) \rangle_{1,t} = 2\int_0^t X_s^N(\log N \phi_s^2 f_0^N(\xi_s^N))\,ds + \int_0^t m_{1,s}^N(\phi_s)\,ds,$$

*where*

$$(4.11) \qquad |m_{1,s}^N(\phi)| \le \frac{C_{4.9} \log N}{\sqrt{N}} \|\phi_s\|_{\mathrm{Lip}}^2 X_s^N(\mathbf{1}), \qquad 0 \le s \le t.$$

(c) *For $i=2,3$, $D_t^{N,i}(\phi) = \int_0^t d_s^{N,i}(\phi)\,ds$, where*

$$|d_s^{N,i}(\phi)| \le C_{4.9} \|\phi_s\|_\infty X_s^N(\log N f_0^N(\xi_s^N)), \qquad 0 \le s \le t.$$

PROOF.　(a) From (3.6) and (1.10) we have

$$|m_{2,s}^N(\phi)| \le \frac{\|\phi_s\|_\infty^2}{(N')^2} \sum_x \frac{2\bar\theta}{N'} \xi_s^N(x) = 2\bar\theta \frac{\|\phi_s\|_\infty^2}{(N')^2} X_s^N(\mathbf{1}).$$

(b) Use $(\xi_s^N(y) - \xi_s^N(x))^2 = \xi_s^N(x)(1 - \xi_s^N(y)) + \xi_s^N(y)(1 - \xi_s^N(x))$ in (3.5) to conclude that

$$\langle M^N(\phi) \rangle_{1,t} = \int_0^t \frac{2\log N}{N'} \sum_x \phi_s^2(x) \xi_s^N(x) \sum_y p_N(y-x)(1-\xi_s^N(y))\,ds$$
$$+ \int_0^t \frac{\log N}{N'} \sum_x \sum_y (\phi_s(x)^2 - \phi_s(y)^2)$$
$$\times p_N(x-y) \xi_s^N(y)(1-\xi_s^N(x))\,ds.$$

This is the required expression, where

$$|m_{1,s}^N(\phi)| = \left| \frac{\log N}{N'} \sum_x \sum_y (\phi_s(x) + \phi_s(y))(\phi_s(x) - \phi_s(y)) \right.$$

$$\times p_N(x-y)\xi_s^N(y)(1-\xi_s^N(x))\bigg|$$

$$\leq (\log N)2\|\phi_s\|_\infty |\phi_s|_{\mathrm{Lip}}\frac{1}{N'}\sum_y \xi_s^N(y)\sum_x |y-x|p_N(y-x)$$

$$\leq C_{4.9}\|\phi_s\|_{\mathrm{Lip}}^2 X_s^N(\mathbf{1})\frac{\log N}{\sqrt{N}},$$

since $\sum_x |x|p_N(x) \leq \sqrt{2\sigma^2/N}$.

(c) Use (3.3) or (3.4), and (1.10) to see that, for $i=2$ or 3,

$$|d_s^{N,i}(\phi)| \leq \bar{\theta}\log N\|\phi_s\|_\infty \frac{1}{N'}\sum_x\sum_y \xi_s^N(x)(1-\xi_s^N(y))p_N(y-x)$$

$$= \bar{\theta}\log N\|\phi_s\|_\infty X_s^N(f_0^N(\xi_s^N)). \qquad \square$$

REMARK 4.9. Inequalities (4.9) and (4.11) imply, in conjunction with (4.2), that for $T>0$, there is a constant $C_{4.12}(T)$ such that if $\phi_s \equiv \psi$, then for $0 \leq s \leq T$,

(4.12) $\quad E(|m_{1,s}^N| + |m_{2,s}^N|) \leq C_{4.12}(T)\|\psi\|_{\mathrm{Lip}}^2(\log N/N^{1/2})X_0^N(\mathbf{1}).$

LEMMA 4.10. *For $T > 0$, there is a constant $C_{4.13}(T)$ so that, for all $0 \leq s \leq t \leq T$,*

(4.13)
$$E\bigg(\bigg[\int_s^t X_r^N(\log N f_0^N(\xi_r^N))\,dr\bigg]^2\bigg)$$
$$\leq C_{4.13}(T)(t-s)^{4/3}(X_0^N(\mathbf{1})^2 + X_0^N(\mathbf{1})).$$

PROOF. Using the Markov property, and then (4.4) and (4.5), the left-hand side above is

$$2E\bigg(\int_s^t\int_r^t X_r^N(\log N f_0^N(\xi_r^N))X_{r'}^N(\log N f_0^N(\xi_{r'}^N))\,dr'\,dr\bigg)$$

$$\leq 2\int_s^t\int_r^t E(X_r^N(\log N f_0^N(\xi_r^N))g(r'-r)X_r^N(\mathbf{1}))\,dr'\,dr$$

$$\leq 2\int_s^t\int_r^t g(r)g(r'-r)\,dr'\,dr(X_0^N(\mathbf{1})^2 + X_0^N(\mathbf{1})).$$

After plugging in the explicit form of $g$ given in (4.1), a little integration yields (4.13). $\square$

PROPOSITION 4.11. *For each $\phi \in C_b^{1,3}(\mathbb{R}_+ \times \mathbb{R}^3)$, each of the families $\{X_\cdot^N(\phi.), N \in \mathbb{N}\}$, $\{D_\cdot^{N,i}(\phi), N \in \mathbb{N}\}$, $i = 1, 2, 3$, $\{\langle M^N(\phi)\rangle_\cdot, N \in \mathbb{N}\}$, and $\{M_\cdot^N(\phi), N \in \mathbb{N}\}$ is C-tight in $D([0,\infty),\mathbb{R})$.*

PROOF. Fix $T > 0$, and $\phi$ as above, and recall the decomposition of $X_t^N(\phi_t)$ in Proposition 3.1. Lemmas 4.8(c) and 4.10, and assumption (1.8) imply there is constant $C_T$ such that, for $i = 2, 3$, $0 \leq s \leq t \leq T$,

$$E((D_t^{N,i}(\phi) - D_s^{N,i}(\phi))^2) \leq C_T \|\phi\|_\infty^2 (t-s)^{4/3}.$$

$C$-tightness of $\{D_\cdot^{N,i}(\phi) : N \in \mathbb{N}\}$ ($i = 2, 3$) is now standard (see, e.g., Theorem 3.8.8 and Proposition 3.10.3 of [4]).

Tightness of $\{\langle M^N(\phi)\rangle_\cdot : N \in \mathbb{N}\}$ follows by similar reasoning using Lemmas 4.8(a), (b) and 4.10, as well as Proposition 4.3(a). $C$-tightness of the remaining terms now follows just as in the argument for $d \geq 3$ in Proposition 3.7 of [2]. In particular, the hypothesis that $\phi \in C_b^{1,3}$ is used to prove $C$-tightness of $\{D_\cdot^{N,1}(\phi)\}$ as in Proposition 3.7 of [2]. $\square$

We come now to the compact containment condition. Let

$$B(x, n) = \{y \in \mathbb{R}^2 : |y - x| \leq n\}.$$

PROPOSITION 4.12. *For all $\epsilon > 0$, there is an $n \in \mathbb{N}$, so that*

$$\sup_N P\left(\sup_{t \leq \epsilon^{-1}} X_t^N(B(0,n)^c) > \epsilon\right) < \epsilon.$$

PROOF. Choose a sequence $h = \{h_n\}$, $h_n : \mathbb{R}^2 \to [0, 1]$, with uniformly (in $n$) bounded continuous partial derivatives of order 3 or less, such that

$$1(|x| > n + 1) \leq h_n(x) \leq 1(|x| > n).$$

Note that $C_h \equiv \sup_n(\|h_n\|_{\mathrm{Lip}} + \|\Delta h_n\|_\infty) < \infty$. By the semimartingale decomposition in Proposition 3.1,

$$(4.14) \quad \sup_{t \leq T} X_t^N(h_n) \leq X_0^N(h_n) + \sum_{i=1}^3 \sup_{t \leq T} |D_t^{N,i}(h_n)| + \sup_{t \leq T} |M_t^N(h_n)|.$$

Our task is to show that the expected value of the right-hand side above tends to 0 as $N, n \to \infty$.

Let $\eta_N = \sup_n \|\mathcal{A}_N(h_n) - \sigma^2 \Delta h_n/2\|_\infty$. As in Lemma 2.6 of [1], our assumptions on $\{h_n\}$ imply that $\lim_{N \to \infty} \eta_N = 0$. It is easy to see that

$$(4.15) \quad |D_t^{N,1}(h_n)| \leq \int_0^t X_s^N(\eta_N \mathbf{1} + C_h h_{n-1}) \, ds$$

and

$$(4.16) \quad |D_t^{N,3}(h_n)| \leq \bar{\theta} \int_0^t X_s^N(h_n \log N f_0^N(\xi_s^N)) \, ds.$$

The term $D_t^{N,2}(h_n)$ is more complicated:

$$|D_t^{N,2}(h_n)| \leq \bar{\theta} \log N \int_0^t \frac{1}{N'} \sum_x h_n(x)(1-\xi_s^N(x)) \sum_y p_N(y-x)\xi_s^N(y)\,ds$$

$$\leq \bar{\theta} \log N \int_0^t \frac{1}{N'} \sum_{x,y} |h_n(x) - h_n(y)| p_N(y-x)\xi_s^N(y)\,ds$$

$$+ \bar{\theta} \log N \int_0^t \frac{1}{N'} \sum_y h_n(y)\xi_s^N(y) f_0^N(y,\xi_s^N)\,ds$$

$$\leq \bar{\theta}|h_n|_{\text{Lip}} \log N \int_0^t \frac{1}{N'} \sum_x \sum_y |x-y| p_N(y-x)\xi_s^N(y)\,ds$$

$$+ \bar{\theta} \int_0^t X_s^N(\log N h_n(y) f_0^N(\xi_s^N))\,ds$$

$$\leq \bar{\theta} C_h \log N \left(\frac{2\sigma^2}{N}\right)^{1/2} \int_0^t X_s^N(\mathbf{1})\,ds$$

$$+ \bar{\theta} \int_0^t X_s^N(h_n \log N f_0^N(\xi_s^N))\,ds,$$

where we have used the symmetry of $p_N(\cdot)$ and the bound $\sum_{x \in \mathbf{S}_N} |x| p_N(x) \leq (2\sigma^2/N)^{1/2}$.

The inequalities above, (4.2) and Burkholder's inequality show that, with $\eta_N'(T) = C_{4.2}(T)(\eta_N + \bar{\theta} C_h \log N (2\sigma^2/N)^{1/2})T$,

$$E\left(\sup_{t \leq T} X_t^N(h_n)\right) \leq X_0^N(h_n) + 2(E\langle M^N(h_n)\rangle_T)^{1/2} + \eta_N' X_0^N(\mathbf{1})$$

(4.17)
$$+ C_h \int_0^T E(X_s^N(h_{n-1}))\,ds$$

$$+ 2\bar{\theta} \int_0^T E(X_s^N(h_n \log N f_0^N(\xi_s^N)))\,ds.$$

[Note $\lim_{N \to \infty} \eta_N'(T) = 0$.]

To estimate the last integral above, we apply Proposition 4.5 with $p = 3$ and $\varepsilon = (\log N)^{-3}$. By the Markov property and (4.2),

$$E\left(\int_0^T X_s^N(h_n \log N f_0^N(\xi_s^N))\,ds\right)$$

$$= E\left(\int_0^\epsilon X_s^N(h_n \log N f_0^N(\xi_s^N))\,ds\right)$$

$$+ E\left(\int_\epsilon^T E(X_s^N(h_n \log N f_0^N(\xi_s^N)) \mid X_{s-\epsilon}^N)\,ds\right)$$

$$\leq \varepsilon(\log N)C_{4.2}(T)X_0^N(\mathbf{1}) + \frac{C_{4.8}\|h_n\|_{\text{Lip}}}{\log N}\int_\varepsilon^T E(X_{s-\varepsilon}^N(\mathbf{1}))\,ds$$

$$+ C_{4.8}\int_\varepsilon^T E(X_{s-\varepsilon}^N(h_n))\,ds.$$

Using (4.2) again, and letting $\eta_N''(T) = C_{4.2}(T)[(\log N)^{-2} + C_{4.8}C_hT/\log N]$, we have

(4.18)
$$E\Big(\int_0^T X_s^N(h_n\log N f_0^N(\xi_s^N))\,ds\Big)$$
$$\leq \eta_N''(T)X_0^N(\mathbf{1}) + C_{4.8}\int_0^T E(X_s^N(h_n))\,ds.$$

[Note $\lim_{N\to\infty}\eta_N''(T) = 0$.]

Next, by the above inequality, Lemma 4.8 and the bound (4.12),

$$E(\langle M^N(h_n)\rangle_T)$$
$$= E\Big(\int_0^T X_s^N(2\log N h_n^2 f_0^N(\xi_s^N))\,ds + \int_0^T (m_{1,s}^N(h_n) + m_{2,s}^N(h_n))\,ds\Big)$$
$$\leq 2\eta_N''(T)X_0^N(\mathbf{1}) + 2C_{4.8}\int_0^T E(X_s^N(h_n))\,ds$$
$$+ C_{4.12}(T)\|h_n\|_{\text{Lip}}^2\frac{\log N}{\sqrt{N}}TX_0^N(\mathbf{1}).$$

That is, letting $\eta_N'''(T) = 2\eta_N''(T) + C_{4.12}(T)TC_h^2\log N/\sqrt{N}$,

(4.19) $\quad E(\langle M^N(h_n)\rangle_T) \leq \eta_N'''(T)X_0^N(\mathbf{1}) + 2C_{4.8}\int_0^T E(X_s^N(h_n))\,ds.$

[Note $\lim_{N\to\infty}\eta_N''' = 0$.]

An application of Proposition 4.4 with $p = 3$ implies

(4.20)
$$\int_0^T E(X_s^N(h_{n-1}))\,ds \leq e^{(\log N)^{-2}}\int_0^T e^{C_{4.7}s}X_0^N(P_s^{N,*}h_{n-1})\,ds$$
$$+ e^{C_{4.7}T}(\log N)^{-1}C_hX_0^N(\mathbf{1}).$$

Finally, let $B_t^{N,*}$ be the continuous time random walk with semigroup $P_t^{N,*}$ defined before Proposition 4.4, $B_0^{N,*} = 0$, and note that

$$E(|B_s^{N,*}|^2) = 2\sigma^2 s \cdot (1 + \bar{\theta}\log N/N) \leq 4\sigma^2 T$$

for all $0 \leq s \leq T$ for large $N$. With this bound, we have

(4.21)
$$X_0^N(P_s^{N,*}h_{n-1}) \leq X_0^N\Big(B\Big(0,\frac{n-1}{2}\Big)^c\Big) + X_0^N(\mathbf{1})P\Big(|B_s^{N,*}| \geq \frac{n-1}{2}\Big)$$
$$\leq X_0^N\Big(B\Big(0,\frac{n-1}{2}\Big)^c\Big) + X_0^N(\mathbf{1})16\sigma^2T/(n-1)^2.$$

By assembling the inequalities (4.17)–(4.21), we have the following: for any $T, \epsilon' > 0$, there is an $N_0$ so that,

$$\text{for } N \geq N_0, n \geq N_0, \qquad E\left(\sup_{t \leq T} X_t^N(h_n)\right) < \epsilon'.$$

The required compact containment follows. □

PROOF OF PROPOSITION 4.1. The $C$-tightness of $\{P_N, N \in \mathbb{N}\}$ is now immediate from Propositions 4.11 and 4.12 above, and Theorem II.4.1 in [8]. □

PROOF OF PROPOSITION 4.2. It is now a simple matter to take a subsequential limit in the semimartingale decomposition of $X_t^N(\phi)$ in Proposition 3.1 to show that any limit point satisfies the martingale problem (3.10) characterizing the law $P_{X_0}^{4\pi\sigma^2, \theta, \sigma^2}$ of super-Brownian motion. Let $\phi \in C_b^\infty(\mathbb{R}^2)$. Proposition 4.7(b) implies $D_t^{N,2}(\phi) - D_t^{N,3}(\phi)$ approaches the drift term involving $\theta$ in (3.10). Proposition 4.7(a), Lemma 4.8(a), (b) and Proposition 4.3(a) show that the square function of the martingale part of $X_t^N(\phi)$ [given by (3.9)] approaches the square function of the martingale part of (3.10) with branching rate $b = 4\pi\sigma^2$. The other terms are handled just as for $d \geq 3$ in [2]. We refer the reader to the proof of Proposition 3.2 there. The only difference in this part of the proof is that we will take $1 < p < 2$ in Proposition 4.7, while in [2] $p = 2$. This leads to only trivial changes. □

**5. Voter, biased voter and Lotka–Volterra bounds.** As in [2], we will obtain bounds on the Lotka–Volterra model by obtaining bounds on the more tractable *biased voter model*. In turn, these bounds depend on good *voter model* bounds. In this section we will work with voter and biased voter models, which we now define.

For $b, v \geq 0$, the 1-biased voter model $\bar{\xi}_t$ is the Feller process taking values in $\{0, 1\}^{\mathbb{Z}^d}$, with rate function

$$(5.1) \qquad \bar{c}(x, \xi) = \begin{cases} (v+b)f_1(x, \xi), & \text{if } \xi(x) = 0, \\ vf_0(x, \xi), & \text{if } \xi(x) = 1, \end{cases}$$

where $f_i(x, \xi)$ is as in (1.1). Similarly, the 0-biased voter model is the Feller process $\underline{\xi}_t$ taking values in $\{0, 1\}^{\mathbb{Z}^d}$, with rate function

$$(5.2) \qquad \underline{c}(x, \xi) = \begin{cases} vf_1(x, \xi), & \text{if } \xi(x) = 0, \\ (v+b)f_0(x, \xi), & \text{if } \xi(x) = 1. \end{cases}$$

The voter model $\hat{\xi}_t$ is the 1-biased voter model with bias $b = 0$, that is, its rate function is $\hat{c}(x, \eta) = vf_i(x, \xi)$ if $\xi(x) = 1 - i$.

It is simple to check that

$$\xi(x) = 0 \quad \text{implies } \underline{c}(x,\xi) \leq \hat{c}(x,\xi) \leq \bar{c}(x,\xi)$$

and

$$\xi(x) = 1 \quad \text{implies } \bar{c}(x,\xi) \leq \hat{c}(x,\xi) \leq \underline{c}(x,\xi).$$

Therefore, as in Theorem III.1.5 of [6], assuming $\underline{\xi}_0 = \hat{\xi}_0 = \bar{\xi}_0$,

(5.3)
we may define $\underline{\xi}_t, \hat{\xi}_t$ and $\bar{\xi}_t$ on a common probability space so that $\underline{\xi}_t \leq \hat{\xi}_t \leq \bar{\xi}_t$ for all $t \geq 0$.

For $\xi, \zeta \in \{0,1\}^{\mathbb{Z}^2}$, $\xi \leq \zeta$ means $\xi(x) \leq \zeta(x)$ for all $x \in \mathbb{Z}^2$.

In Section 5.1 we will obtain the required voter model bounds. In Section 5.2 we will use these bounds to obtain good biased voter model bounds.

NOTE. We will assume throughout the rest of this section that (5.3) is in force.

5.1. *Voter model estimates.* We recall the voter model duality; see, for instance, [5] or [6]. Recall also the system of coalescing random walks $\{\hat{B}^x : x \in \mathbb{Z}^2\}$ from Section 1. The basic duality equation for the rate one ($v = 1$) voter model is as follows: for finite $A \subset \mathbb{Z}^2$,

(5.4) $$P(\hat{\xi}_t(x) = 1 \ \forall x \in A) = P(\hat{\xi}_0(\hat{B}_t^x) = 1 \ \forall x \in A).$$

Recall $\tilde{P} = \sum_e p(e) P^e$, $\tau_x$ and $H(t)$ from Section 2, and define the mean range of the random walk $B_t^0$ by

$$R(t) = E\left(\sum_x 1\{B_s^0 = x \text{ for some } s \leq t\}\right).$$

A *last time at 0* decomposition (see, e.g., Lemma A.2 of [1]) yields $R(t) = 1 + \int_0^t H(s)\,ds$, and [via (2.4)] the well-known asymptotic behavior

(5.5) $$\lim_{t \to \infty} \frac{R(t)}{t/\log t} = 2\pi\sigma^2.$$

Let $P_t, t \geq 0$, be the semigroup of a rate 1 random walk with step distribution $p(\cdot)$. We slightly abuse our earlier notation and for $\phi : \mathbb{Z}^2 \to \mathbb{R}$ and $\xi \in \{0,1\}^{\mathbb{Z}^2}$, let

$$\xi(\phi) = \sum_x \phi(x)\xi(x).$$

LEMMA 5.1. *Let $\hat{\xi}_t$ denote the rate-$v$ voter model. Then for all bounded $\phi: Z^2 \to \mathbb{R}^+$ and $t \geq 0$,*

$$E(\hat{\xi}_t(\phi)) = \hat{\xi}_0(P_{vt}\phi), \tag{5.6}$$

$$E(|\hat{\xi}_t|^2) \leq |\hat{\xi}_0|^2 + 2vt|\hat{\xi}_0|, \tag{5.7}$$

$$E(\hat{\xi}_t(\phi f_0(\hat{\xi}_t))) \leq (2\sigma^2 vt H(2vt))^{1/2} |\phi|_{\text{Lip}} |\bar{\xi}_0| + H(2vt)\hat{\xi}_0(\phi), \tag{5.8}$$

$$E(|\hat{\xi}_t|\hat{\xi}_t(f_0(\hat{\xi}_t))) \leq H(2vt)|\hat{\xi}_0|^2 + R(2vt)|\hat{\xi}_0|. \tag{5.9}$$

REMARK 5.2. For $\phi = \mathbf{1}$, the right-hand side of (5.8) is just $H(2vt)|\hat{\xi}_0|$.

PROOF OF LEMMA 5.1. By scaling, it suffices to consider the case $v = 1$. Also, the first two formulas are well known [the latter follows from (5.14) below], so we prove only the last two. By the duality equation (5.4), symmetry and translation invariance,

$$\begin{aligned} E(\hat{\xi}_t(\phi f_0(\hat{\xi}_t))) &= \sum_{x,e} \phi(x) p(e) P(\hat{B}_t^x \in \hat{\xi}_0, \hat{B}_t^{x+e} \notin \hat{\xi}_0) \\ &\leq \sum_{x,e,z} \hat{\xi}_0(z) \phi(x) p(e) P(\hat{B}_t^x = z, \hat{\tau}(x, x+e) > t) \\ &= \sum_{x,e,z} \hat{\xi}_0(z) \phi(x) p(e) P(B_t^0 = x - z, \tau(0, e) > t) \\ &= \sum_{e,z} \hat{\xi}_0(z) p(e) E(\phi(z + B_t^0) 1\{\tau(0, e) > t\}). \end{aligned}$$

For any $z$,

$$\begin{aligned} \sum_e p(e) E(\phi(z + B_t^0) 1\{\tau(0,e) > t\}) \\ &\leq \sum_e p(e) E((|\phi|_{\text{Lip}} |B_t^0| + \phi(z)) 1\{\tau(0,e) > t\}) \\ &\leq |\phi|_{\text{Lip}} \left( E(|B_t^0|^2) \sum_e p(e) P(\tau(0,e) > t) \right)^{1/2} \\ &\quad + \phi(z) \sum_e p(e) P(\tau(0,e) > t). \end{aligned}$$

Since $E(|B_t^0|^2) = 2\sigma^2 t$, this proves (5.8).

We expand the left-hand side of (5.9) and use duality to obtain

$$E(|\hat{\xi}_t|\hat{\xi}_t(f_0(\hat{\xi}_t))) = \Gamma_1 + \Gamma_2,$$

where

(5.10)    $\Gamma_1 = \sum\limits_{x,y,e} p(e) P(\hat{B}_t^y \in \hat{\xi}_0, \hat{B}_t^x \in \hat{\xi}_0, \hat{B}_t^{x+e} \notin \hat{\xi}_0, \hat{\tau}(x,y) > t),$

(5.11)    $\Gamma_2 = \sum\limits_{x,y,e} p(e) P(\hat{B}_t^x \in \hat{\xi}_0, \hat{B}_t^{x+e} \notin \hat{\xi}_0, \hat{\tau}(x,y) \le t).$

Consider $\Gamma_1$ first, which we expand in the form

$$\sum_{x,y,w,z,e} 1(y \ne x+e) p(e)$$
$$\times P(\hat{B}_t^y = w, \hat{B}_t^x = z, \hat{B}_t^{x+e} \notin \hat{\xi}_0, \hat{\tau}(x,y) > t) \hat{\xi}_0(w) \hat{\xi}_0(z).$$

By replacing the condition $\hat{B}_t^{x+e} \notin \hat{\xi}_0$ with $\hat{\tau}(x, x+e) > t$, switching to the independent random walk system [dropping the condition that $\hat{\tau}(x,y) > t$], it follows that $\Gamma_1$ is bounded above by

$$\sum_{x,y,w,z,e} 1(y \ne x+e) p(e) P(\hat{B}_t^y = w, \hat{B}_t^x = z, \hat{\tau}(x,y) > t, \hat{\tau}(x, x+e) > t)$$
$$\times \hat{\xi}_0(w) \hat{\xi}_0(z)$$
$$\le \sum_{x,y,w,z,e} 1(y \ne x+e, y \ne x) p(e) P(B_t^y = w, B_t^x = z, \tau(x, x+e) > t)$$
$$\times \hat{\xi}_0(w) \hat{\xi}_0(z).$$

Now changing variables, we have

$$\Gamma_1 \le \sum_{x,y,w,z,e} 1(y \ne 0, y \ne e) p(e) P(B_t^{y+x} = w+x, B_t^x = z+x, \tau(x, x+e) > t)$$
$$\times \hat{\xi}_0(w+x) \hat{\xi}_0(z+x)$$
$$= \sum_{x,y,w,z,e} 1(y \ne 0, y \ne e) p(e) P(B_t^y = w, B_t^0 = z, \tau(0,e) > t)$$
$$\times \hat{\xi}_0(w+x) \hat{\xi}_0(z+x)$$
$$= \sum_{x,y,w,z,e} 1(y \ne 0, y \ne e) p(e) P(B_t^y = w) P(B_t^0 = z, \tau(0,e) > t)$$
$$\times \hat{\xi}_0(w+x) \hat{\xi}_0(z+x).$$

Since $P(B_t^y = w) = P(B_t^0 = w - y)$, summing in order over $y, w, x, z$ shows that the last sum above is at most $|\hat{\xi}_0|^2 P^*(\tau > 2t)$. Thus, to prove (5.9), it suffices to show $\Gamma_2 \le R(2t) |\hat{\xi}_0|$.

In the definition of $\Gamma_2$ we drop the restriction $\hat{B}_t^{x+e} \notin \hat{\xi}_0$ and then sum over $e$ to obtain

$$\Gamma_2 \le \sum_{x,y,z} P(\hat{B}_t^x = z, \hat{\tau}(x,y) \le t) \hat{\xi}_0(z)$$

$$= \sum_{x,y,z} P(B_t^x = z+x, \tau(x, y+x) \leq t)\hat{\xi}_0(z+x)$$

$$= \sum_{x,y,z} P(B_t^0 = z, \tau(0, y) \leq t)\hat{\xi}_0(z+x),$$

where we have again changed variables. If we sum in order over $x, z, y$, we obtain $\Gamma_2 \leq |\hat{\xi}_0| R(2t)$, and we are done. $\square$

5.2. *Biased voter model bounds.* We first recall the following from Lemma 4.1 of [2]. If $\bar{\xi}_t$ is the 1-biased voter model with rate function (5.1), then

(5.12) $$E(|\bar{\xi}_t|) \leq e^{bt}|\bar{\xi}_0|,$$

(5.13) $$E(|\bar{\xi}_t|^2) \leq e^{2bt}\left(|\bar{\xi}_0|^2 + \frac{2v+b}{b}(1 - e^{-bt})|\bar{\xi}_0|\right).$$

Since $1 - e^{-bt} \leq bt$, the last inequality implies

(5.14) $$E(|\bar{\xi}_t|^2) \leq e^{2bt}(|\bar{\xi}_0|^2 + (2v+b)t|\bar{\xi}_0|).$$

These bounds must be improved. In (6.2) below we will compare the Lotka–Volterra model $\xi_t^N$ defined in the Introduction with the biased voter models $\underline{\xi}_t^N$, $\bar{\xi}_t^N$ on $\mathbf{S}_N$. In order to construct the coupling $\underline{\xi}_t^N \leq \xi_t^N \leq \bar{\xi}_t^N$, we must assume that the voting and bias rates $v_N$ and $b_N$ are

(5.15) $$v = v_N = N - \bar{\theta} \log N \quad \text{and} \quad b = b_N = 2\bar{\theta} \log N.$$

With this coupling, the bounds on $X_t^N(\mathbf{1})$ in Section 4 are then consequences of analogous bounds on $\bar{X}_t^N(\mathbf{1}) = (1/N)\sum_x \bar{\xi}_t^N(x)\delta_x$. However, for the above rates, the bound (5.12) implies only that $E(\bar{X}_t^N(\mathbf{1})) \leq e^{b_N t} X_0^N(\mathbf{1})$, not that $E(\bar{X}_t^N(\mathbf{1}))$ [and hence, $E(X_t^N(\mathbf{1}))$] is bounded in $N$, a fact we will need. Nevertheless, the estimates (5.12) and (5.13) are useful over short time periods, and will play an important role in deriving better bounds.

To state our improved versions of (5.12) and (5.13), we define several constants and functions depending on $b$ and $v$. For $p \geq 2$, define

(5.16)
$$\kappa_p = \kappa_p(b, v) = 3(bH(2v/b^p) + e^2) \quad \text{and} \quad \kappa = \kappa_3,$$
$$A = A(b, v) = bR(2v/b^3) + 3e^2(1 + 2v/b),$$
$$B_p = B_p(b, v) = (2\sigma^2 v b^{2-p} H(2v/b^p))^{1/2}$$
$$\qquad + bH(2v/b^p)(2\sigma^2(1 + v/b^p))^{1/2}$$

and

(5.17)
$$h_1(b,v)(t) = e^2 t^{-1/3} + \kappa 2 e^{2+2\kappa t},$$
$$h_2(b,v)(t) = e^2 t^{-1/3}(1 + 2v/b) + 5\kappa A e^{1+3\kappa t}.$$

Also, let $P\phi(x) = \sum_y p(y-x)\phi(y)$ and define the operators

(5.18) $$\bar{\mathcal{A}}\phi = v(P\phi - \phi) \quad \text{and} \quad \mathcal{A}^* = (1+b/v)\bar{\mathcal{A}}.$$

Let $\bar{\mathcal{A}}$ (resp., $\mathcal{A}^*$) have associated semigroup $\bar{P}_t$, $t \geq 0$ (resp., $P_t^*$, $t \geq 0$).

REMARK 5.3. The constants $\kappa_p, A, B_p$ and the functions $h_1, h_2$ are used in many bounds below. These bounds are not sharp, but they are adequate for our purposes. Note that for the parameters $v = v_N, b = b_N$ in (5.15), we have [by (2.4) and (5.5)] $\kappa_p = O(1)$, $A = O(N/\log N)$ and $B_p = O(N^{1/2}(\log N)^{(1-p)/2})$ as $N \to \infty$.

PROPOSITION 5.4. *Assume $b \geq 1$ and $p \geq 2$. For all $t \geq 0$,*

(5.19) $$E(|\bar{\xi}_t|) \leq e^{b^{1-p} + \kappa_p t}|\bar{\xi}_0|,$$

(5.20) $$E(|\bar{\xi}_t|^2) \leq e^{2+2\kappa t}|\bar{\xi}_0|^2 + 4Ae^{1+3\kappa t}|\bar{\xi}_0|,$$

(5.21) $$bE(\bar{\xi}_t(f_0(\bar{\xi}_t))) \leq h_1(t)|\bar{\xi}_0|,$$

(5.22) $$bE(|\bar{\xi}_t|\bar{\xi}_t(f_0(\bar{\xi}_t))) \leq h_1(t)|\bar{\xi}_0|^2 + h_2(t)|\bar{\xi}_0|.$$

*For all bounded $\phi : \mathbb{Z}^2 \to [0, \infty)$ and $p \geq 3$,*

(5.23) $$\begin{aligned}E(\bar{\xi}_t(\phi)) &\leq e^{b^{1-p} + (1+\kappa_p)t} \\ &\quad \times (\bar{\xi}_0(P_t^*(\phi)) + [\kappa_p b^{2-p}\|\phi\|_\infty + B_p|\phi|_{\mathrm{Lip}}]|\bar{\xi}_0|).\end{aligned}$$

To derive these bounds, we first state a special case of Proposition 2.3 of [2]. The biased voter models $\bar{\xi}_t$ and $\underline{\xi}_t$ are special cases of the general voter model perturbations introduced in [2]; for $\bar{\xi}$, take $N = v$, $\beta(\{e\}) = bp(e)$, $\beta(A) = 0$ for $\mathrm{card}(A) \neq 1$ and $\delta \equiv 0$, while for $\underline{\xi}$ take $N = v$, $\delta(\{e\}) = -bp(e)$, $\delta(\varnothing) = b$, $\delta(A) = 0$ for $\mathrm{card}(A) \neq 0$ or 1, and $\beta \equiv 0$. This notation is only important if you want to verify that Lemma 5.5 is indeed a special case of Proposition 2.3 of [2].

LEMMA 5.5. *Let $T > 0$ and $\phi : [0,T] \times Z^2 \to \mathbb{R}$, where $\phi, \dot{\phi}$ are both bounded and continuous. Then for $0 \leq t \leq T$,*

(a)

(5.24) $$\begin{aligned}\bar{\xi}_t(\phi_t) &= \bar{\xi}_0(\phi_0) + \int_0^t \bar{\xi}_s(\bar{\mathcal{A}}(\phi_s) + \dot{\phi}_s)\,ds \\ &\quad + b\int_0^t \sum_x \phi(s,x)[1 - \bar{\xi}_s(x)]f_1(x, \bar{\xi}_s)\,ds + \bar{M}_t(\phi),\end{aligned}$$

where $\bar{M}_t(\phi)$ is a square-integrable $(\mathcal{F}_t^{\bar{\xi}})$-martingale with predictable square function

(5.25) $$\langle \bar{M}(\phi)\rangle_t = \langle \bar{M}(\phi)\rangle_{1,t} + \langle \bar{M}(\phi)\rangle_{2,t},$$

with

$$\langle \bar{M}(\phi)\rangle_{1,t} = \int_0^t v \sum_x \phi_s(x)^2 [\bar{\xi}_s(x) f_0(x,\bar{\xi}_s) + (1-\bar{\xi}_s(x)) f_1(x,\bar{\xi}_s)]\, ds$$

and

$$\langle \bar{M}(\phi)\rangle_{2,t} = \int_0^t b \sum_x \phi_s(x)^2 (1-\bar{\xi}_s(x)) f_1(x,\bar{\xi}_s)\, ds.$$

(b)

(5.26) $$\underline{\xi}_t(\phi_t) = \underline{\xi}_0(\phi_0) + \int_0^t \underline{\xi}_s(\bar{\mathcal{A}}(\phi_s) + \dot{\phi}_s)\, ds$$
$$- b\int_0^t \underline{\xi}_s(\phi_s f_0(\underline{\xi}_s))\, ds + \underline{M}_t(\phi),$$

where $\underline{M}_t(\phi)$ is a square-integrable $(\mathcal{F}_t^{\underline{\xi}})$-martingale with predictable square function

(5.27) $$\langle \underline{M}(\phi)\rangle_t = \langle \underline{M}(\phi)\rangle_{1,t} + \langle \underline{M}(\phi)\rangle_{2,t},$$

with

$$\langle \underline{M}(\phi)\rangle_{1,t} = \int_0^t v \sum_x \phi_s(x)^2 [\underline{\xi}_s(x) f_0^N(x,\underline{\xi}_s) + (1-\underline{\xi}_s(x)) f_1(x,\underline{\xi}_s)]\, ds$$

and

$$\langle \underline{M}(\phi)\rangle_{2,t} = \int_0^t b\underline{\xi}_s(\phi_s^2 f_0(\underline{\xi}_s))\, ds.$$

If we set $\phi \equiv 1$ in Lemma 5.5 and do a bit of stochastic calculus, we see that

(5.28) $$|\bar{\xi}_t| - |\bar{\xi}_0| - b\int_0^t \bar{\xi}_s(f_0(\bar{\xi}_s))\, ds$$

and

(5.29) $$|\bar{\xi}_t|^2 - |\bar{\xi}_0|^2 - (2v+b)\int_0^t \bar{\xi}_s(f_0(\bar{\xi}_s))\, ds - 2b\int_0^t |\bar{\xi}_s|\bar{\xi}_s(f_0(\bar{\xi}_s))\, ds$$

are $(\mathcal{F}_t^{\bar{\xi}})$-martingales.

In the proof of Proposition 5.4 we will need a few properties of the function $\xi(f_0(\xi))$ which will be important to "transfer" the voter model bounds of Lemma 5.1 to the biased voter model.

LEMMA 5.6. *Assume $\xi, \eta \in \{0,1\}^{\mathbb{Z}^d}$ satisfy $\xi \leq \eta$.*
(a) *If $|\eta| < \infty$, then*

$$|\eta(f_0(\eta)) - \xi(f_0(\xi))| \leq |\eta| - |\xi| \tag{5.30}$$

*and*

$$||\eta|\eta(f_0(\eta)) - |\xi|\xi(f_0(\xi))| \leq |\eta|^2 - |\xi|^2. \tag{5.31}$$

(b) *If $\phi: \mathbb{Z}^2 \to \mathbb{R}^+$ is bounded, then*

$$\sum_x |\eta(x) f_0(x, \eta) - \xi(x) f_0(x, \xi)| \phi(x) \leq (\eta - \xi)(\phi + P\phi), \tag{5.32}$$

$$\sum_x |\eta(x) f_0(x, \eta)^2 - \xi(x) f_0(x, \xi)^2| \phi(x) \leq (\eta - \xi)(\phi + 2P\phi) \tag{5.33}$$

*and*

$$\sum_x |(1 - \eta(x)) f_1(x, \eta)^2 - (1 - \xi(x)) f_1(x, \xi)^2| \phi(x) \tag{5.34}$$
$$\leq (\eta - \xi)(\phi + 2P\phi).$$

PROOF. (a) The first step is

$$\eta(f_0(\eta)) - \xi(f_0(\xi))$$
$$= \sum_x (\eta(x) - \xi(x)) f_0(x, \eta) + \sum_x \xi(x) (f_0(x, \eta) - f_0(x, \xi)).$$

Since $f_0(x, \eta) - f_0(x, \xi) = f_1(x, \xi) - f_1(x, \eta)$, and $p(\cdot)$ is symmetric,

$$\sum_x \xi(x)(f_0(x,\eta) - f_0(x,\xi)) = \sum_x \xi(x)(f_1(x,\xi) - f_1(x,\eta))$$
$$= \sum_{x,y} \xi(x) p(y-x)(\xi(y) - \eta(y))$$
$$= \sum_y (\xi(y) - \eta(y)) f_1(y, \xi).$$

Thus,

$$\eta(f_0(\eta)) - \xi(f_0(\xi)) = \sum_x (\eta(x) - \xi(x))(f_0(x,\eta) - f_1(x,\xi)). \tag{5.35}$$

Since $|f_0(x,\eta) - f_1(x,\xi)| \leq 1$, (5.30) is an immediate consequence of (5.35). It follows from (5.30) that $|\xi| - \xi(f_0(\xi)) \leq |\eta| - \eta(f_0(\eta))$. Multiplying the left-hand side of this equality by $|\xi|$ and the right-hand side by $|\eta|$, we obtain [note that $|\xi| - \xi(f_0(\xi)) \geq 0$]

$$|\xi|^2 - |\xi|\xi(f_0(\xi)) \leq |\eta|^2 - |\eta|\eta(f_0(\eta)),$$

or $|\eta|\eta(f_0(\eta)) - |\xi|\xi(f_0(\xi)) \leq |\eta|^2 - |\xi|^2$. Since (5.30) also implies that $|\xi| + \xi(f_0(\xi)) \leq |\eta| + \eta(f_0(\eta))$, multiplying the left-hand side by $|\xi|$ and the right-hand side by $|\eta|$ yields the inequality $|\eta|\eta(f_0(\eta)) - |\xi|\xi(f_0(\xi)) \geq -(|\eta|^2 - |\xi|^2)$.

(b) These are similar so we only prove (5.34). The left-hand side is bounded above by

$$\sum_x [(1-\eta(x))(f_1(x,\eta)^2 - f_1(x,\xi)^2)\phi(x) + f_1(x,\xi)^2(\eta(x) - \xi(x))\phi(x)]$$

$$\leq \sum_x (1-\eta(x))2(f_1(x,\eta) - f_1(x,\xi))\phi(x) + (\eta(x) - \xi(x))\phi(x)$$

$$\leq 2\sum_x \sum_y p(y-x)\phi(x)(\eta(y) - \xi(y)) + (\eta - \xi)(\phi),$$

which equals the right-hand side of (5.34). □

PROOF OF (5.19) AND (5.23). Let $\phi$ be as in Lemma 5.5 and rewrite the second integral in (5.24) in the form

$$b\int_0^t \sum_{x,y} (\phi(s,x) - \phi(s,y))p(y-x)[1 - \bar{\xi}_s(x)]\bar{\xi}_s(y)\,ds$$

$$+ b\int_0^t \sum_y \phi(s,y)\bar{\xi}_s(y)f_0(y,\bar{\xi}_s)\,ds$$

$$= \frac{b}{v}\int_0^t \bar{\xi}_s(\bar{\mathcal{A}}\phi_s)\,ds + \int_0^t \bar{\xi}_s(b\phi_s f_0(\bar{\xi}_s))\,ds.$$

We have used the fact $\sum_{x,y}(\phi(s,x) - \phi(s,y))p(x-y)\bar{\xi}_s(x)\bar{\xi}_s(y) = 0$ by symmetry. Therefore, (5.24) becomes

(5.36)
$$\bar{\xi}_t(\phi_t) = \bar{\xi}_0(\phi_0) + \int_0^t \bar{\xi}_s(\mathcal{A}^*\phi_s + \dot{\phi}_s)\,ds$$
$$+ \int_0^t \bar{\xi}_s(\phi_s b f_0(\bar{\xi}_s))\,ds + \bar{M}_t(\phi).$$

Fix $c \geq 0$, $t > 0$ and $\phi: \mathbb{Z}^2 \to \mathbb{R}^+$, and set $\phi_s = e^{-cs}P_{t-s}^*\phi$. Then (5.36) implies

(5.37)
$$E(\bar{\xi}_t(\phi))e^{-ct} = \bar{\xi}_0(P_t^*\phi) + \int_0^t E(\bar{\xi}_s(bP_{t-s}^*\phi f_0(\bar{\xi}_s)))e^{-cs}\,ds$$
$$- \int_0^t E(\bar{\xi}_s(cP_{t-s}^*\phi))e^{-cs}\,ds.$$

This implies

(5.38) $$\bar{\xi}_0(P_t^*\phi) \leq E(\bar{\xi}_t(\phi)) \leq e^{bt}\bar{\xi}_0(P_t^*\phi),$$

where we have set $c = 0$ for the first inequality and $c = b$ for the second.

We now work at upgrading the second inequality in (5.38) to (5.19) by estimating the first integral in (5.37) via a comparison with the voter model. Put $\varepsilon = b^{-p}$ and assume $\phi \geq 0$. It follows from the coupling (5.3) and the inequalities (5.32), (5.6) and (5.12) that

$$
\begin{aligned}
E(|\bar{\xi}_\varepsilon(b\phi f_0(\bar{\xi}_\varepsilon)) - \hat{\xi}_\varepsilon(b\phi f_0(\hat{\xi}_\varepsilon))|) \\
\leq 2b\|\phi\|_\infty E(|\bar{\xi}_\varepsilon| - |\hat{\xi}_\varepsilon|) \leq 2b(e^{b\varepsilon} - 1)\|\phi\|_\infty |\bar{\xi}_0|.
\end{aligned}
\tag{5.39}
$$

In view of this bound, $b \geq 1$, and the voter model estimate (5.8), we have

$$
\begin{aligned}
E(\bar{\xi}_\varepsilon(b\phi f_0(\bar{\xi}_\varepsilon))) \leq 2eb^2\varepsilon\|\phi\|_\infty |\bar{\xi}_0| \\
+ b(2\sigma^2 v\varepsilon H(2v\varepsilon))^{1/2}|\phi|_{\mathrm{Lip}}|\bar{\xi}_0| + bH(2v\varepsilon)\bar{\xi}_0(\phi).
\end{aligned}
\tag{5.40}
$$

We have used the elementary inequality $e^u - 1 \leq eu$ for $0 \leq u \leq 1$ above (and will make use of it again without comment). The Markov property now implies, for $s \geq \varepsilon$, that

$$
\begin{aligned}
E(\bar{\xi}_s(b\phi f_0(\bar{\xi}_s)) \mid \mathcal{F}_{s-\varepsilon}) \\
\leq (2eb^2\varepsilon\|\phi\|_\infty + b(2\sigma^2 v\varepsilon H(2v\varepsilon))^{1/2}|\phi|_{\mathrm{Lip}})|\bar{\xi}_{s-\varepsilon}| \\
+ bH(2v\varepsilon)\bar{\xi}_{s-\varepsilon}(\phi).
\end{aligned}
\tag{5.41}
$$

We can now derive (5.19). Taking expectations in (5.41) for the function $\phi = \mathbf{1}$, we obtain, for $s \geq \varepsilon$,

$$
E(\bar{\xi}_s(bf_0(\bar{\xi}_s))) \leq \kappa_p E(|\bar{\xi}_{s-\varepsilon}|).
\tag{5.42}
$$

Using this inequality in (5.36), we get, for $t \geq \varepsilon$,

$$
E(|\bar{\xi}_t|) \leq E(|\bar{\xi}_\varepsilon|) + \kappa_p \int_\varepsilon^t E(|\bar{\xi}_{s-\varepsilon}|)\,ds \leq e^{b\varepsilon}|\bar{\xi}_0| + \kappa_p \int_0^t E(|\bar{\xi}_s|)\,ds,
$$

where (5.38) is used in the second inequality. This bound also holds for $t \leq \varepsilon$, and hence, (5.19) follows by Gronwall's inequality.

To prove (5.23), we will take expectations in (5.41), with $\phi$ replaced by $P^*_{t-s}\phi$, and substitute this in (5.37). However, we must first alter the last term of (5.41). For bounded $\psi : \mathbb{Z}^2 \to \mathbb{R}^+$, (5.38) implies that

$$
|E(\bar{\xi}_\varepsilon(\psi)) - \bar{\xi}_0(\psi)| \leq (e^{b\varepsilon} - 1)\bar{\xi}_0(P^*_\varepsilon\psi) + |\bar{\xi}_0(P^*_\varepsilon\psi) - \bar{\xi}_0(\psi)|.
$$

Since $|P^*_\varepsilon\psi(x) - \psi(x)| \leq |\psi|_{\mathrm{Lip}} E(|B^0_{v\varepsilon(1+b/v)}|) \leq |\psi|_{\mathrm{Lip}}(2\sigma^2\varepsilon(v+b))^{1/2}$,

$$
|E(\bar{\xi}_\varepsilon(\psi)) - \bar{\xi}_0(\psi)| \leq (eb\varepsilon\|\psi\|_\infty + (2\sigma^2\varepsilon(v+b))^{1/2}|\psi|_{\mathrm{Lip}})|\bar{\xi}_0|.
$$

Consequently, by the Markov property, for $s \geq \varepsilon$,

$$
E(\bar{\xi}_{s-\varepsilon}(\psi)) \leq E(\bar{\xi}_s(\psi)) + (eb\varepsilon\|\psi\|_\infty + (2\sigma^2\varepsilon(v+b))^{1/2}|\psi|_{\mathrm{Lip}})E(|\bar{\xi}_{s-\varepsilon}|).
$$

Using this inequality in (5.41), with $\psi = P^*_{t-s}\phi$ replacing $\phi$, we have for $s \geq \varepsilon$

(5.43)
$$E(\bar{\xi}_s(bP^*_{t-s}\phi f_0(\bar{\xi}_s)))$$
$$\leq (\kappa_p b^2 \varepsilon \|\phi\|_\infty + B_p |\phi|_{\mathrm{Lip}}) E(|\bar{\xi}_{s-\varepsilon}|) + \kappa_p E(\bar{\xi}_s(P^*_{t-s}\phi)).$$

Plugging this into (5.37) with $c = 1 + \kappa_p$, and using (5.38) for $s \leq \varepsilon$, and (5.19), we obtain

$$E(\bar{\xi}_t(\phi))e^{-ct} \leq \bar{\xi}_0(P^*_t \phi) + \int_0^\varepsilon b e^{bs} \bar{\xi}_0(P^*_s P^*_{t-s}\phi)e^{-cs}\,ds$$

$$+ \int_\varepsilon^t (\kappa_p b^2 \varepsilon \|\phi\|_\infty + B_p|\phi|_{\mathrm{Lip}})E(|\bar{\xi}_{s-\varepsilon}|)e^{-cs}\,ds$$

$$+ \kappa_p \int_\varepsilon^t E(\bar{\xi}_s(P^*_{t-s}\phi))e^{-cs}\,ds - c\int_\varepsilon^t E(\bar{\xi}_s(P^*_{t-s}\phi))e^{-cs}\,ds$$

$$\leq \bar{\xi}_0(P^*_t \phi)(1 + e^{b\varepsilon} - 1)$$

$$+ (\kappa_p b^2 \varepsilon \|\phi\|_\infty + B_p|\phi|_{\mathrm{Lip}})\int_0^t e^{b\varepsilon + \kappa_p s - (1+\kappa_p)s}|\bar{\xi}_0|\,ds,$$

that is,

$$E(\bar{\xi}_t(\phi)) \leq e^{b\varepsilon + (1+\kappa_p)t}(\bar{\xi}_0(P^*_t(\phi)) + (\kappa_p b^2 \varepsilon \|\phi\|_\infty + B_p|\phi|_{\mathrm{Lip}})|\bar{\xi}_0|).$$

This proves (5.23). □

PROOF OF (5.20). By (5.29),

(5.44)
$$E(|\bar{\xi}_t|^2) = |\bar{\xi}_0|^2 + (2v+b)\int_0^t E(\bar{\xi}_s(f_0(\bar{\xi}_s)))\,ds$$
$$+ 2b\int_0^t E(|\bar{\xi}_s|\bar{\xi}_s(f_0(\bar{\xi}_s)))\,ds.$$

The first integral in (5.44) is easy to handle. By (5.28) and (5.19),

(5.45)
$$(2v+b)\int_0^t E(\bar{\xi}_s(f_0(\bar{\xi}_s)))\,ds = (2v+b)\frac{E(|\bar{\xi}_t|) - |\bar{\xi}_0|}{b}$$
$$\leq \left(1 + \frac{2v}{b}\right)e^{1+\kappa t}|\bar{\xi}_0| \leq Ae^{1+\kappa t}|\bar{\xi}_0|.$$

We continue to write $\varepsilon = b^{-p}$, but now set $p = 3$. Since the integrand of the second integral in (5.44) is bounded by $E(|\bar{\xi}_s|^2)$, (5.14) implies

(5.46)
$$2b\int_0^\varepsilon E(|\bar{\xi}_s|\bar{\xi}_s(f_0(\bar{\xi}_s)))\,ds \leq 2b\int_0^\varepsilon e^{2bs}(|\bar{\xi}_0|^2 + (2v+b)s|\bar{\xi}_0|)\,ds$$
$$\leq (e^{2b\varepsilon} - 1)(|\bar{\xi}_0|^2 + (2v+b)\varepsilon|\bar{\xi}_0|)$$
$$\leq (e^2 - 1)|\bar{\xi}_0|^2 + A|\bar{\xi}_0|.$$

For the integral over $[\varepsilon, t]$, we use (5.31) and the voter model estimate (5.9),

$$E(|\bar{\xi}_\varepsilon|\bar{\xi}_\varepsilon(f_0(\bar{\xi}_\varepsilon))) \leq E(|\hat{\xi}_\varepsilon|\hat{\xi}_\varepsilon(f_0(\hat{\xi}_\varepsilon))) + E(|\bar{\xi}_\varepsilon|^2 - |\hat{\xi}_\varepsilon|^2)$$
$$\leq H(2v\varepsilon)|\bar{\xi}_0|^2 + R(2v\varepsilon)|\bar{\xi}_0| + E(|\bar{\xi}_\varepsilon|^2) - |\bar{\xi}_0|^2$$

since $E(|\hat{\xi}_\varepsilon|^2) \geq |\hat{\xi}_0|^2$ ($|\hat{\xi}_t|^2$ is a submartingale). The bound (5.14) now implies

$$bE(|\bar{\xi}_\varepsilon|\bar{\xi}_\varepsilon(f_0(\bar{\xi}_\varepsilon))) \leq \kappa |\bar{\xi}_0|^2 + A|\bar{\xi}_0|. \tag{5.47}$$

Consequently, by the Markov property and (5.19),

$$2b \int_\varepsilon^t E(|\bar{\xi}_s|\bar{\xi}_s(f_0(\bar{\xi}_s))) \, ds \leq 2\kappa \int_\varepsilon^t E(|\bar{\xi}_{s-\varepsilon}|^2) \, ds + 2A \int_\varepsilon^t E(|\bar{\xi}_{s-\varepsilon}|) \, ds$$
$$\leq 2\kappa \int_\varepsilon^t E(|\bar{\xi}_{s-\varepsilon}|^2) \, ds + 2Ae^{1+\kappa t}|\bar{\xi}_0|. \tag{5.48}$$

Combining this with (5.46) gives

$$2b \int_0^t E(|\bar{\xi}_s|\bar{\xi}_s(f_0(\bar{\xi}_s))) \, ds$$
$$\leq (e^2 - 1)|\bar{\xi}_0|^2 + (A + 2Ae^{1+\kappa t})|\bar{\xi}_0| + 2\kappa \int_0^t E(|\bar{\xi}_s|^2) \, ds. \tag{5.49}$$

In view of (5.44), (5.45) and (5.49),

$$E(|\bar{\xi}_t|^2) \leq e^2|\bar{\xi}_0|^2 + 4Ae^{1+\kappa t}|\bar{\xi}_0| + 2\kappa \int_0^t E(|\bar{\xi}_s|^2) \, ds \tag{5.50}$$

for all $t \geq 0$. It follows from Gronwall's inequality that

$$E(|\bar{\xi}_t|^2) \leq \{e^2|\bar{\xi}_0|^2 + 4Ae^{1+\kappa t}|\bar{\xi}_0|\} e^{2\kappa t}. \qquad \square$$

PROOF OF (5.21). First suppose $t \leq \varepsilon = b^{-3}$. Then $b \leq t^{-1/3}$ and $bt \leq 1$, and using (5.12),

$$bE(\bar{\xi}_t(f_0(\bar{\xi}_t))) \leq bE(|\bar{\xi}_t|) \leq et^{-1/3}|\bar{\xi}_0|. \tag{5.51}$$

For $t \geq \varepsilon$, it follows from (5.42) and (5.19) that

$$bE(\bar{\xi}_t(f_0(\bar{\xi}_t))) \leq \kappa E(|\bar{\xi}_{t-\varepsilon}|) \leq \kappa e^{1+\kappa t}|\bar{\xi}_0|. \tag{5.52}$$

The inequalities (5.51) and (5.52) imply (5.21). $\square$

PROOF OF (5.22). For $t \leq \varepsilon = b^{-3}$, we may apply the second moment estimate (5.14) to obtain

$$E(|\bar{\xi}_t|\bar{\xi}_t(f_0(\bar{\xi}_t))) \leq E(|\bar{\xi}_\varepsilon|^2) \leq e^{2b\varepsilon}(|\bar{\xi}_0|^2 + \varepsilon(2v+b)|\bar{\xi}_0|). \tag{5.53}$$

Arguing as before leads to

$$bE(|\bar{\xi}_t|\bar{\xi}_t(f_0(\bar{\xi}_t))) \leq e^2 t^{-1/3}(|\bar{\xi}_0|^2 + (1+2v/b)|\bar{\xi}_0|). \tag{5.54}$$

For $t \geq \varepsilon$, we apply the Markov property, (5.47) and the submartingale property of $|\bar{\xi}_t|^2$ and $|\bar{\xi}_t|$ to obtain

$$bE(|\bar{\xi}_t|\bar{\xi}_t(f_0(\bar{\xi}_t))) \leq \kappa E(|\bar{\xi}_t|^2) + AE(|\bar{\xi}_t|).$$

After using the bounds (5.19) and (5.20) for $E(|\bar{\xi}_t|)$ and $E(|\bar{\xi}_t|^2)$, respectively, and rearranging, we arrive at inequality (5.22). $\square$

We need only a few bounds for the 0-biased voter model $\underline{\xi}_t$. Using (5.26) directly with $\phi_s = e^{cs} \bar{P}_{t-s} \phi$ [recall (5.18)] for bounded $\phi : \mathbb{Z}^2 \to \mathbb{R}^+$ and $c \geq 0$, we get the analogue of (5.37),

$$E(\underline{\xi}_t(\phi))e^{ct} = \underline{\xi}_0(\bar{P}_t \phi) - \int_0^t E(\underline{\xi}_s(b\bar{P}_{t-s}\phi f_0(\underline{\xi}_s)))e^{cs}\,ds$$
$$+ \int_0^t E(\underline{\xi}_s(c\bar{P}_{t-s}\phi))e^{cs}\,ds. \tag{5.55}$$

Setting $\phi = \mathbf{1}$ gives

$$E(|\underline{\xi}_t|)e^{ct} = |\underline{\xi}_0| - b\int_0^t E(\underline{\xi}_s(f_0(\underline{\xi}_s)))e^{cs}\,ds + c\int_0^t E(|\underline{\xi}_s|)e^{cs}\,ds. \tag{5.56}$$

By setting $c = b$ above, we get the simple bound

$$E(|\underline{\xi}_t|) \geq e^{-bt}|\underline{\xi}_0|, \tag{5.57}$$

which implies

$$E(|\underline{\xi}_t|^2) \geq e^{-2bt}|\underline{\xi}_0|^2. \tag{5.58}$$

Since $|\underline{\xi}_t| \leq |\bar{\xi}_t|$, the above and (5.14) therefore imply

$$|E(|\underline{\xi}_t|^2) - |\underline{\xi}_0|^2| \leq (e^{2bt}-1)|\underline{\xi}_0|^2 + e^{2bt}(2v+b)t|\underline{\xi}_0|. \tag{5.59}$$

The next step is to improve the bound (5.57).

LEMMA 5.7. *Assume $b \geq 1$. Then for all $t \geq 0$,*

$$E(|\underline{\xi}_t|) \geq \exp(-b^{-2} - \kappa e^{b^{-2}}t)|\underline{\xi}_0|. \tag{5.60}$$

PROOF. With $\varepsilon = b^{-3}$, by the coupling $\underline{\xi}_\varepsilon \leq \hat{\xi}_\varepsilon$ in (5.3) and the bounds (5.30), (5.6), Remark 5.2 and (5.57), we get

$$bE(\underline{\xi}_\varepsilon(f_0(\underline{\xi}_\varepsilon))) \leq (bH(2v\varepsilon) + b^2\varepsilon)|\underline{\xi}_0| \leq \kappa|\underline{\xi}_0|. \tag{5.61}$$

For $t \geq \varepsilon$, (5.56) implies

$$E(|\underline{\xi}_t|)e^{ct} = e^{c\varepsilon}E(|\underline{\xi}_\varepsilon|) + \int_\varepsilon^t E(ce^{cs}|\underline{\xi}_s| - be^{cs}\underline{\xi}_s(f_0(\underline{\xi}_s))) \, ds.$$

By (5.57), (5.61) and the Markov property,

$$E(|\underline{\xi}_t|)e^{ct} \geq e^{(c-b)\varepsilon}|\underline{\xi}_0| + \int_\varepsilon^t E(ce^{-b\varepsilon}|\underline{\xi}_{s-\varepsilon}| - \kappa|\underline{\xi}_{s-\varepsilon}|)e^{cs} \, ds,$$

and hence,

$$E(|\underline{\xi}_t|) \geq e^{-b\varepsilon}e^{-ct}|\bar{\xi}_0| + e^{-c(t-\varepsilon)} \int_0^{t-\varepsilon} (ce^{-b\varepsilon} - \kappa)E|\underline{\xi}_s|e^{cs} \, ds.$$

Put $c = \kappa e^{b\varepsilon}$ to obtain (5.60). □

We conclude this section with two final inequalities which are useful for small $t \geq 0$. The first one follows easily from (5.12) and (5.57), the second from (5.14) and (5.58). For $\underline{\xi}_0 = \bar{\xi}_0 = \xi_0$,

(5.62) $\quad 0 \leq E(|\bar{\xi}_t|) - E(|\underline{\xi}_t|) \leq 2(e^{bt} - 1)|\xi_0|,$

(5.63) $\quad 0 \leq E(|\bar{\xi}_t|^2) - E(|\underline{\xi}_t|^2) \leq 2(e^{2bt} - 1)|\xi_0|^2 + e^{2bt}(2v + b)t|\xi_0|.$

**6. Proofs of Propositions 4.3–4.5.** Return now to the rescaled regime of the Section 1 and let $\xi_t^N \in \{0,1\}^{\mathbf{S}_N}$ be the rescaled Lotka–Volterra model with rate function $c_N(x, \xi)$ given by (1.6), where (1.9) continues to hold. Let $\bar{\xi}_t \in \{0,1\}^{\mathbb{Z}^2}$ be the 1-biased voter model and $\underline{\xi}_t \in \{0,1\}^{\mathbb{Z}^2}$ the 0-biased voter model with rates $v = v_N = N - \bar{\theta} \log N$ and $b = b_N = 2\bar{\theta} \log N$ given in (5.15), and set $\bar{\xi}_t^N(x) = \bar{\xi}_t(x\sqrt{N})$ and $\underline{\xi}_t^N(x) = \underline{\xi}_t(x\sqrt{N})$ for $x \in \mathbf{S}_N$. Thus, $\bar{\xi}_t^N \in \{0,1\}^{\mathbf{S}_N}$ has rate function

$$\bar{c}_N(x, \xi) = \begin{cases} (v_N + b_N)f_1^N(x, \xi), & \text{if } \xi(x) = 0, \\ v_N f_0^N(x, \xi), & \text{if } \xi(x) = 1, \end{cases}$$

and $\underline{\xi}_t^N \in \{0,1\}^{\mathbf{S}_N}$ has rate function

$$\underline{c}_N(x, \xi) = \begin{cases} v_N f_1^N(x, \xi), & \text{if } \xi(x) = 0, \\ (v_N + b_N)f_0^N(x, \xi), & \text{if } \xi(x) = 1. \end{cases}$$

We assume that $N$ is large enough ($N \geq N_0$) so that $v_N > 0$ and $b_N > 1$.

It is easy to check that $\underline{c}_N(x, \xi) \leq c_N(x, \xi) \leq \bar{c}_N(x, \xi)$ if $\xi(x) = 0$ and $\underline{c}_N(x, \xi) \geq c_N(x, \xi) \geq \bar{c}_N(x, \xi)$ if $\xi(x) = 1$. Thus, as in (5.3), assuming $\underline{\xi}_0^N = \xi_0^N = \bar{\xi}_0^N$, we may construct the three processes on one probability space so that

(6.1) $\qquad \underline{\xi}_t^N \leq \xi_t^N \leq \bar{\xi}_t^N \qquad$ for all $t \geq 0$.

Letting $\bar{X}_t^N = \frac{1}{N'}\sum_{x\in \mathbf{S}_N} \bar{\xi}_t^N(x)\delta_x$ and $\underline{X}_t^N = \frac{1}{N'}\sum_{x\in \mathbf{S}_N} \underline{\xi}_t^N(x)\delta_x$, it follows that $\underline{X}_0^N = X_0^N = \bar{X}_0^N$ and

(6.2) $$\underline{X}_t^N \leq X_t^N \leq \bar{X}_t^N \qquad \text{for all } t \geq 0.$$

We begin with bounds for the 1-biased voter model. Proposition 5.4 and Remark 5.3 imply that there are constants $C_{6.3}$ and $C_{4.1}$ such that if $g$ is as in (4.1), then for all $N \geq N_0, t \geq 0$,

(6.3) $$E(\bar{X}_t^N(\mathbf{1})) \leq C_{6.3} e^{C_{6.3} t} \bar{X}_0^N(\mathbf{1}),$$

(6.4) $$E(\bar{X}_t^N(\mathbf{1})^2) \leq C_{6.3} e^{C_{6.3} t}(\bar{X}_0^N(\mathbf{1})^2 + \bar{X}_0^N(\mathbf{1}))$$

and

(6.5) $$(\log N) E(\bar{X}_t^N(f_0^N(\cdot, \bar{\xi}_t^N))) \leq g(t)\bar{X}_0^N(\mathbf{1}),$$

(6.6) $$(\log N) E(\bar{X}_t^N(\mathbf{1})\bar{X}_t^N(f_0^N(\cdot, \bar{\xi}_t^N))) \leq g(t)(\bar{X}_0^N(\mathbf{1})^2 + \bar{X}_0^N(\mathbf{1})).$$

We also need the following comparison result, which follows from the coupling (6.1), inequality (5.19) and Lemma 5.7, our choice of $v_N$ and $b_N$, and Remark 5.3. There is a constant $C_{6.7}$ such that

(6.7) $$\begin{aligned} E(\bar{X}_t^N(\mathbf{1})) - E(\underline{X}_t^N(\mathbf{1})) \\ \leq C_{6.7}[(\log N)^{-2} + t]X_0^N(\mathbf{1}), \qquad 0 \leq t \leq 1. \end{aligned}$$

Since the coupling (6.2) does not allow us to compare $X_t^N(f_0(\xi_t^N))$ and $\bar{X}_t^N(f_0(\bar{\xi}_t^N))$, we need the following. Recall that $h_1, h_2$ are defined in (5.17).

PROPOSITION 6.1. *For $t \geq 0$,*

(6.8) $$bE(\xi_t^N(f_0(\xi_t^N))) \leq h_1(t)|\xi_0^N|$$

*and*

(6.9) $$bE(|\xi_t^N|\xi_t^N(f_0(\xi_t^N))) \leq 2(h_1(t)|\xi_0^N|^2 + h_2(t)|\xi_0^N|).$$

PROOF. For $t \leq \varepsilon = b^{-3}$, $bE(\xi_t^N(f_0(\xi_t))) \leq bE(|\xi_t^N|) \leq bE(|\bar{\xi}_t^N|)$, so just as in (5.51) we obtain

(6.10) $$bE(\xi_t^N(f_0(\xi_t))) \leq et^{-1/3}|\xi_0^N|.$$

To handle $t \geq \varepsilon$, we compare with the biased voter model using (5.30):

$$\begin{aligned} bE(\xi_\varepsilon^N(f_0(\xi_\varepsilon^N))) &\leq bE(\bar{\xi}_\varepsilon^N(f_0(\bar{\xi}_\varepsilon^N))) + bE(|\bar{\xi}_\varepsilon^N| - |\xi_\varepsilon^N|) \\ &\leq b(H(2v\varepsilon) + 2eb\varepsilon + 2(e^{b\varepsilon} - 1))|\xi_0^N| \\ &\leq \kappa |\xi_0^N|, \end{aligned}$$

by (5.40), (5.62) and the coupling (6.1), and the definition of $\kappa$. Our standard argument with the Markov property now gives, for $t \geq \varepsilon$,

$$bE(\xi_t^N(f_0(\xi_t^N))) \leq \kappa E(|\xi_{t-\varepsilon}^N|) \leq \kappa E(|\bar{\xi}_t^N|) \leq \kappa e^{1+\kappa t}|\xi_0^N|,$$

the last by (5.19). This inequality and (6.10) imply (6.8).

Now consider (6.9). For $t \leq \varepsilon$, (5.14) implies

$$\begin{aligned}
bE(|\xi_t^N|\xi_t^N(f_0(\xi_t^N))) &\leq bE(|\xi_t^N|^2) \leq bE(|\bar{\xi}_t^N|^2) \\
(6.11) \quad &\leq be^{2b\varepsilon}(|\xi_0^N|^2 + (2v+b)\varepsilon|\xi_0^N|) \\
&\leq e^2 t^{-1/3}(|\xi_0^N|^2 + (1+2v/b)|\xi_0^N|),
\end{aligned}$$

where the last inequality follows as in (5.51). By comparing with the biased voter model using (5.31),

$$bE(|\xi_\varepsilon^N|\xi_\varepsilon^N(f_0(\xi_\varepsilon^N))) \leq bE(|\bar{\xi}_\varepsilon^N|\bar{\xi}_\varepsilon^N(f_0(\bar{\xi}_\varepsilon^N))) + b(E(|\bar{\xi}_\varepsilon^N|^2) - E(|\xi_\varepsilon^N|^2)).$$

Consequently, by the above bound, (5.47) and (5.63),

$$bE(|\xi_\varepsilon^N|\xi_\varepsilon^N(f_0(\xi_\varepsilon^N))) \leq 2(\kappa|\xi_0^N|^2 + A|\xi_0^N|).$$

It follows from the Markov property that, for all $t \geq \varepsilon$,

$$\begin{aligned}
bE(|\xi_t^N|\xi_t^N(f_0(\xi_t^N))) &\leq 2(\kappa E(|\xi_{t-\varepsilon}^N|^2) + AE(|\xi_{t-\varepsilon}^N|)) \\
&\leq 2(\kappa E(|\bar{\xi}_t^N|^2) + AE(|\bar{\xi}_t^N|)).
\end{aligned}$$

We have used (5.28) and (5.29) in the last inequality. Using (5.19) and (5.20) in the above, and also recalling (6.11), we get (6.9), as required. $\square$

PROOF OF PROPOSITION 4.3. Part (a) follows from the coupling (6.2), the fact that $\bar{X}_t^N(\mathbf{1})^2$ is a submartingale [see (5.29)], the strong $L^2$ inequality for non-negative submartingales, and the bounds (6.3) and (6.4). Part (b) follows directly from the previous Proposition and Remark 5.3. $\square$

Proposition 4.4 is a direct consequence of the coupling (6.2) and the following biased voter model bound. Recall the notation introduced in (4.6).

PROPOSITION 6.2. *For $p \geq 3$, there is a constant $C_{6.12}(p)$ such that, for any $t \geq 0$ and $\psi : \mathbb{R}^2 \to \mathbb{R}^+$,*

$$\begin{aligned}
(6.12) \quad E(\bar{X}_t^N(\psi)) &\leq e^{(\log N)^{1-p}} e^{C_{6.12} t} \bar{X}_0^N(P_t^{N,*}\psi) \\
&\quad + C_{6.12} e^{C_{6.12} t} \|\psi\|_{\mathrm{Lip}} (\log N)^{(1-p)/2} \bar{X}_0^N(\mathbf{1}).
\end{aligned}$$

PROOF. Fix $N$, and recall that in this section $\bar{\xi}_t \in \{0,1\}^{\mathbb{Z}^2}$ is the biased voter model with rates $v = N - \bar{\theta}\log N$ and $b = 2\bar{\theta}\log N$, and that $\bar{\xi}_t^N(x) = \xi_t(x\sqrt{N})$, $x \in \mathbf{S}_N$. Define $\phi:\mathbb{Z}^2 \to \mathbb{R}^+$ by $\phi(x) = \psi(x/\sqrt{N})$. Then $\|\phi\|_\infty = \|\psi\|_\infty$, $|\phi|_{\text{Lip}} = N^{-1/2}|\psi|_{\text{Lip}}$, $P_t^*\phi(x) = P_t^{N,*}\psi(x/\sqrt{N})$ [$P_t^*$ is defined after (5.18)] and $\bar{\xi}_t^N(\psi) = \bar{\xi}_t(\phi)$. Applying (5.23),

$$E(\bar{\xi}_t^N(\psi)) \leq e^{b^{1-p}+(1+\kappa_p)t}(\bar{\xi}_0^N(P_t^{N,*}\psi)$$
$$+ [\kappa_p b^{2-p}\|\psi\|_\infty + B_p N^{-1/2}|\psi|_{\text{Lip}}]|\bar{\xi}_0^N|).$$

Since $p \geq 3$, it follows from Remark 5.3 that $\kappa_p b^{2-p} + B_p N^{-1/2} = O((\log N)^{(1-p)/2})$ as $N \to \infty$, and thus, (6.12) holds because $\bar{\theta} \geq 1$ implies $b \geq \log N$. □

PROOF OF PROPOSITION 4.5. Let $\varepsilon = b^{-p}$. By (5.32),

(6.13)
$$E(X_\varepsilon^N(b\phi f_0(\xi_\varepsilon^N)))$$
$$\leq E(\bar{X}_\varepsilon^N(b\phi f_0(\bar{\xi}_\varepsilon^N))) + 2b\|\phi\|_\infty(E(\bar{X}_\varepsilon^N(\mathbf{1}) - X_\varepsilon^N(\mathbf{1}))).$$

Now applying (5.40) and (5.62), we find that

$$E(X_\varepsilon^N(b\phi f_0(\xi_\varepsilon^N))) \leq (6eb^{2-p}\|\phi\|_\infty + B_p N^{-1/2}|\phi|_{\text{Lip}})X_0^N(\mathbf{1}) + \kappa_p X_0^N(\phi).$$

Using our standard asymptotics for $B_p$ and $\kappa_p$, we obtain (4.8). □

**7. Proof of Proposition 4.7—Part I.** For $\phi:\mathbb{R}^2 \to \mathbb{R}$, $\zeta \in \{0,1\}^{\mathbf{S}_N}$ and $X(\phi) = (1/N')\sum_x \phi(x)\zeta(x)$, define

$$\Delta_1^{N,+}(\phi,\zeta) = X(\log N \phi^2 f_0^N(\cdot,\zeta)),$$
$$\Delta_2^{N,+}(\phi,\zeta) = \frac{1}{N'}\sum_x (1-\zeta(x))\phi(x)\log N f_1^N(x,\zeta)^2,$$
$$\Delta_3^{N,+}(\phi,\zeta) = X(\log N \phi f_0^N(\cdot,\zeta)^2)$$

and

$$\Delta_j^N(\phi,\zeta) = \Delta_j^{N,+}(\phi,\zeta) - \gamma_j X(\phi), \qquad j = 1,2,3,$$

where $\gamma_1 = 2\pi\sigma^2$ and $\gamma_2 = \gamma_3 = \gamma^*$. An easy calculation using (4.3) and (1.9) shows that to prove Proposition 4.7 it suffices to prove the following: for $0 < p < 2$, there exists $\eta_N = \eta_N(p) \downarrow 0$ as $N \to \infty$ such that, for all $K, T > 0$, there exists $C(p, K, T)$ such that, for all $t \leq T$, $\phi:\mathbb{R}^2 \to [0,\infty)$ and $\xi_0^N$ such that

(7.1) $$\|\phi\|_{\text{Lip}} \vee X_0^N(\mathbf{1}) \leq K,$$

we have

$$E\left(\left|\int_0^t \Delta_j^N(\phi, \xi_s^N)\,ds\right|^p\right) \leq \eta_N C(p, K, T), \qquad j = 1, 2, 3. \tag{7.2}$$

We will do this by following the general strategy we have already used, comparing the Lotka–Volterra with the voter model over short time periods and estimating the difference via comparisons with the biased voter model. We must be more careful than before because we need more precise estimates than, say, the ones in Proposition 4.3.

Define the sequences

$$\varepsilon_N = (\log \log N)^{-1}, \qquad t_N = \frac{\varepsilon_N}{\log N},$$
$$K_N = (\log N)^{1/2}, \qquad \delta_N = K_N t_N. \tag{7.3}$$

We assume that $N$ is large enough so that $\varepsilon_N \vee t_N \vee \delta_N \leq 1$ and $\delta_N \leq (\log N)^{-1/2}$. Note that $\delta_N/\varepsilon_N \to 0$ as $N \to \infty$. We assume the Lotka–Volterra model on $\mathbf{S}_N$, $\xi_\cdot^N$, and the 1-biased voter model $\bar{\xi}_\cdot^N$ on $\mathbf{S}_N$ are as in the previous section.

This section is devoted to proving the following two results.

PROPOSITION 7.1. *There is a constant $C_{7.4}(K)$ and sequence $\eta_{7.4}(N) \downarrow 0$ so that, for all $\phi$ satisfying (7.1) and $j = 1, 2, 3$,*

$$|E(\Delta_j^N(\phi, \xi_{t_N}^N))| \leq C_{7.4}(K)$$
$$\times \left( \eta_{7.4}(N)(X_0^N(\mathbf{1}) + X_0^N(\mathbf{1})^2) \right.$$
$$\left. + \varepsilon_N^{-1} \int\!\!\int 1(|w - z| \leq \sqrt{\delta_N})\,dX_0^N(w)\,dX_0^N(z) \right). \tag{7.4}$$

PROPOSITION 7.2. *There is a constant $C_{7.5}$ such that, for all $0 \leq t \leq T$,*

$$E\left(\int\!\!\int 1(|x - y| \leq \sqrt{\delta_N})\,dX_t^N(x)\,dX_t^N(y)\right)$$
$$\leq C_{7.5} e^{C_{7.5} T}(X_0^N(\mathbf{1}) + X_0^N(\mathbf{1})^2)$$
$$\times \left[ \frac{\delta_N}{\delta_N + t}(1 + t^{2/3}) + \delta_N t^{-1/3} \log\left(1 + \frac{t}{\delta_N}\right) \right]. \tag{7.5}$$

We divide our work into four parts. In the first we collect together various random walk estimates that we will need. In the second we will establish Proposition 7.5, the voter model version of Proposition 7.1. In the third part we prove Proposition 7.1 by comparison with the voter model, and in the fourth we prove a biased voter model analogue of Proposition 7.2, which implies Proposition 7.2. The key inequality (7.2) is then proved in Section 8, completing the proof of Proposition 4.7.

7.1. *Random walk estimates.* Recall from Section 1 that $B_t^0$ is a rate 1 continuous time random walk with step distribution $p(\cdot)$, starting at 0.

LEMMA 7.3. (a) *There is a constant $C_{7.6}$, depending on $p(\cdot)$, such that*

$$\sup_{x \in \mathbb{Z}^2} P(B_t^0 = x) \leq C_{7.6}(t+1)^{-1} \qquad \text{for all } t \geq 0 \tag{7.6}$$

*and*

$$\sup_{t \geq 0} P(B_t^0 = x) \leq C_{7.6}(|x|^2 + 1)^{-1} \qquad \text{for all } x \in \mathbb{Z}^2. \tag{7.7}$$

(b) *If $z_T \in \mathbb{Z}^2$ and $t_T > 0$ satisfy*

$$\lim_{T \to \infty} \frac{z_T}{\sqrt{T}} = z \quad \text{and} \quad \lim_{T \to \infty} \frac{t_T}{T} = s > 0, \tag{7.8}$$

*then*

$$\lim_{T \to \infty} TP(B_{t_T}^0 = z_T) = \frac{e^{-|z|^2/2\sigma^2 s}}{2\pi\sigma^2 s}. \tag{7.9}$$

(c) *For each $K > 0$, there is a constant $C_{7.10}(K) > 0$ so that*

$$\liminf_{T \to \infty} \inf_{|x| \leq K\sqrt{T}} TP(B_T^0 = x) \geq C_{7.10}(K). \tag{7.10}$$

REMARK 7.4. If $p(\cdot)$ is a kernel on $\mathbb{Z}^d$ (any $d$) satisfying the conditions of Section 1 and also $\sum_x |x|^d p(x) < \infty$, then part (a) and it's proof go through, where $|x|^2$ may now be replaced with $|x|^d$ in (7.7).

PROOF OF LEMMA 7.3. The first inequality is standard for discrete time walks; see (A.7) in [1] for the simple adaptation to continuous time. For the second, let $Y_n$ be the discrete time random walk with step distribution $p(\cdot)$ introduced in Section 2. Then P7.10 in [10] implies there is a constant $C_{7.6}$ such that

$$\sup_{n \geq 0} P^0(Y_n = x) \leq C_{7.6}(|x|^2 + 1)^{-1} \qquad \text{for all } x \in \mathbb{Z}^2.$$

Since $P(B_t^0 = x) = \sum_{n=0}^\infty e^{-t} t^n P(Y_n = x)/n!$, (7.7) follows immediately.

For $\varepsilon > 0$, let $Y_n^\varepsilon, n = 0, 1, 2, \ldots,$ be the discrete time random walk with step distribution $p^\varepsilon(x) = P(B_\varepsilon = x)$. Applying the discrete time local central limit theorem (P7.10 in [10]) to this random walk, we can conclude that if (7.8) holds, assuming that $t_T \in \varepsilon \mathbb{Z}^+$ for all $T$, then (7.9) must hold (note the step variance of $Y^\varepsilon$ is $\sigma^2 \varepsilon$). The fact that $p(\cdot)$ is symmetric implies that $P(B_t^0 = 0)$ is decreasing in $t$, and therefore, the Markov property implies that

$$P(B_{n\varepsilon}^0 = x)p^\varepsilon(0) \leq P(B_u^0 = x) \leq \frac{P(B_{(n+1)\varepsilon}^0 = x)}{p^\varepsilon(0)}$$

for $u \in [n\varepsilon, (n+1)\varepsilon]$.

This inequality and an argument by contradiction shows that (7.8) implies (7.9) without the restriction $t_T \in \varepsilon \mathbb{Z}^+$.

For the walk $Y_n$, let $c_n(K) = \inf_{|x| \leq K\sqrt{n}} nP^0(Y_n = x)$. Then P7.10 of [10] implies that $\liminf_{n \to \infty} c_n(K) > 0$ for every $K > 0$. Let $S(t), t \geq 0$, be a rate one Poisson process. For all $|x| \leq K\sqrt{t}$,

$$tP(B_t^0 = x) \geq \sum_{|t-n|<t/2} tP(S(t) = n) \inf_{|y| \leq K\sqrt{2n}} P(Y_n = y)$$
$$\geq P(|S(t) - t| < t/2)(2/3) \inf_{n>t/2} c_n(K\sqrt{2}).$$

This is enough to prove (7.10). $\square$

7.2. *Voter model estimates.* Let $\varepsilon_N, t_N, K_N, \delta_N$ be as in (7.3). For $N$ fixed, let $\hat{\xi}_t$ be the rate $v_N = N - \bar{\theta} \log N$ voter model on $\mathbb{Z}^2$ with rates as in (5.1) for $b = 0$ and $v = v_N$. Define $\hat{\xi}_t^N(x) = \hat{\xi}_t(x\sqrt{N}), x \in \mathbf{S}_N$, the rate $v_N$ voter model on $\mathbf{S}_N$. We introduce some rather trivial notation which will be used frequently:

(7.11) $\qquad m(1) = 2 \quad \text{and} \quad m(2) = m(3) = 1.$

Our goal here is to prove the following analogue of Proposition 7.1 for $\hat{\xi}_t^N$.

PROPOSITION 7.5. *There is a constant $C_{7.12}$ and a sequence $\eta_{7.12}(N) \downarrow 0$ so that, for $j = 1, 2, 3$, if $\phi : \mathbb{R}^2 \to \mathbb{R}$, then*

$$|E(\Delta_j^N(\phi, \hat{\xi}_{t_N}^N))|$$
(7.12) $\qquad \leq \eta_{7.12}(N)(\hat{X}_0^N(\mathbf{1}) + \hat{X}_0^N(\mathbf{1})^2) \|\phi\|_{\mathrm{Lip}}^{m(j)}$
$$+ \frac{C_{7.12} \|\phi\|_\infty^{m(j)}}{\varepsilon_N} \iint 1(|w-z| \leq \sqrt{\delta_N}) \, d\hat{X}_0^N(w) \, d\hat{X}_0^N(z).$$

To prepare for the proof of this result, we introduce rescaled versions of the independent and coalescing random walks systems $\{B_t^x\}$ and $\{\hat{B}_t^x\}$ introduced in Section 1 as follows: for $x, y \in \mathbf{S}_N$,

(7.13) $\qquad B_t^{N,x} = B_{v_N t}^{x\sqrt{N}}/\sqrt{N}, \qquad \hat{B}_t^{N,x} = \hat{B}_{v_N t}^{x\sqrt{N}}/\sqrt{N},$

and
$$\tau^N(x,y) = \tau(\sqrt{N}x, \sqrt{N}y)/v_N, \qquad \hat{\tau}^N(x,y) = \hat{\tau}(\sqrt{N}x, \sqrt{N}y)/v_N.$$

We will need the following estimate.

LEMMA 7.6. *There is a constant $C_{7.14}$ such that*

(7.14)
$$\frac{\log N}{N'} \sum_{x,e} p_N(e) P(\hat{\xi}_0^N(B_{t_N}^{N,x}) = \hat{\xi}_0^N(B_{t_N}^{N,x+e}) = 1, \tau^N(x, x+e) > t_N)$$
$$\leq C_{7.14} \varepsilon_N^{-1} \iint 1(|w-z| \leq \sqrt{\delta_N}) \, d\hat{X}_0^N(w) \, d\hat{X}_0^N(z)$$
$$+ C_{7.14} \frac{\hat{X}_0^N(\mathbf{1})^2}{K_N \varepsilon_N}.$$

REMARK 7.7. In Section 5.1, working with the (unscaled) voter model, we estimated quantities like the left-hand side above by, in effect, dropping the condition $\hat{\xi}_0^N(B_{t_N}^{N,x+e}) = 1$. The bound this produces is too crude for our current needs.

PROOF. By translation invariance and symmetry, the left-hand side of (7.14) is

(7.15)
$$(N')^{-2} \sum_{w,z} \hat{\xi}_0^N(w) \hat{\xi}_0^N(z) \sum_e p_N(e)$$
$$\times \left[ \sum_x N P(B_{t_N}^{N,0} = w - x, B_{t_N}^{N,e} = z - x, \tau^N(0,e) > t_N) \right]$$
$$= (N')^{-2} \sum_{w,z} \hat{\xi}_0^N(w) \hat{\xi}_0^N(z) \sum_e p_N(e) N P(B_{2t_N}^{N,e} = z - w, \tau_0^{N,e} > 2t_N)$$
$$\equiv \Sigma_d^N + \Sigma_c^N,$$

where $\tau_0^{N,e} = \inf\{s : B_s^{N,e} = 0\}$, and $\Sigma_d^N$, respectively, $\Sigma_c^N$, denotes the contribution to (7.15) from $w, z$ satisfying $|w - z| \leq \sqrt{K_N t_N}$, respectively, $|w - z| > \sqrt{K_N t_N}$. Define $\tilde{P}^N$ by

$$\tilde{P}^N((B_\cdot^N, \tau_0^N) \in \cdot) = \sum_e p_N(e) P((B_\cdot^{N,e}, \tau_0^{N,e}) \in \cdot),$$

and note that
$$\tilde{E}^N(|B_s^N|^2) = \sum_{e \in \mathbf{S}_N} p_N(e) E(|e + B_s^{N,0}|^2)$$
$$= \sum_{e \in \mathbf{S}_N} p_N(e)(|e|^2 + E(2\langle e, B_s^{N,0}\rangle + |B_s^{N,0}|^2)) = 2\sigma^2(N^{-1} + s).$$

For $\Sigma_d^N$, use (7.6) and the Markov property at time $t_N$ to see that

$$N\tilde{P}^N(B_{2t_N}^N = z - w, \tau_0^N > 2t_N)$$
$$\leq N\tilde{E}(1(\tau_0^N > t_N)(\omega)P(B_{t_N}^{N,0} = z - w - B_{t_N}^N(\omega)))$$
$$\leq N\tilde{P}(\tau_0^N > t_N)C_{7.6}(v_N t_N)^{-1}$$
$$\leq C_{7.6}\frac{N}{v_N}H(v_N t_N)/t_N.$$

In view of (2.4), there is a constant $C_{7.16}$ such that

$$(7.16) \quad \Sigma_d^N \leq C_{7.16}\varepsilon_N^{-1}\iint 1(|w - z| \leq \sqrt{K_N t_N})\,d\hat{X}_0^N(w)\,d\hat{X}_0^N(z).$$

To bound $\Sigma_c^N$, let $\hat{\eta}_N = e^{-\sqrt{\log N}}$ and use the Markov property at time $\hat{\eta}_N t_N$ and the bounds in Lemma 7.3(a) [recall (7.13)] to see that, for $|w - z| > \sqrt{K_N t_N}$,

$$\tilde{P}^N(B_{2t_N}^N = w - z, \tau_0^N > 2t_N)$$
$$\leq \tilde{E}^N\left(1\left(\tau_0^N > \hat{\eta}_N t_N, |B_{\hat{\eta}_N t_N}^N| > \frac{\sqrt{K_N t_N}}{2}\right)\right)\sup_{x'}P(B_{(2-\hat{\eta}_N)t_N}^{N,0} = x')$$
$$+ \tilde{E}^N\left(1\left(\tau_0^N > \hat{\eta}_N t_N, |B_{\hat{\eta}_N t_N}^N| \leq \frac{\sqrt{K_N t_N}}{2}\right)\right.$$
$$\left.\times P(B_{(2-\hat{\eta}_N)t_N}^{N,0} = w - z - B_{\hat{\eta}_N t_N}^N)\right)$$
$$\leq 4\frac{\tilde{E}^N(|B_{\hat{\eta}_N t_N}^N|^2)}{K_N t_N}\frac{C_{7.6}}{v_N(2-\hat{\eta}_N)t_N} + 4H(v_N\hat{\eta}_N t_N)\frac{C_{7.6}}{NK_N t_N}.$$

Using $\tilde{E}^N(|B_{\hat{\eta}_N t_N}^N|^2) = 2\sigma^2(\hat{\eta}_N t_N + N^{-1})$ and $H(v) \sim c/\log N$ as $v \to \infty$ [i.e., (2.4)], and then plugging in the value of these constants, one shows that each of the last two terms above is $O((N(\log N)K_N t_N)^{-1})$ as $N \to \infty$. Hence, there is a constant $C_{7.17}$ so that

$$(7.17) \quad \Sigma_c^N \leq \frac{C_{7.17}\hat{X}_0^N(\mathbf{1})^2}{\varepsilon_N K_N}.$$

Put (7.17) and (7.16) into (7.15) to get the required bound. $\square$

PROOF OF PROPOSITION 7.5. Fix $\phi:\mathbb{R}^2 \to \mathbb{R}$ such that $\|\phi\|_{\mathrm{Lip}} < \infty$, and let $e, f$ denote independent random variables with law $p_N(\cdot)$. The structure of the proof is as follows. We will use duality to decompose each $E(\Delta_j^{N,+})$ into a sum of simpler terms [defined below in (7.20)–(7.22)],

$$
\begin{aligned}
E(\Delta_1^{N,+}(\phi,\hat{\xi}_{t_N}^N)) &= \Sigma_1^{1,N} - \Sigma_1^{2,N}, \\
E(\Delta_2^{N,+}(\phi,\hat{\xi}_{t_N}^N)) &= \Sigma_2^{1,N} - \Sigma_2^{2,N} + \Sigma_2^{3,N}, \\
E(\Delta_3^{N,+}(\phi,\hat{\xi}_{t_N}^N)) &= \Sigma_3^{1,N} - \Sigma_3^{2,N} + \Sigma_3^{3,N}.
\end{aligned} \tag{7.18}
$$

We will show that there is a sequence $\eta(N) \to 0$ such that, for $j = 1, 2, 3$,

$$|\Sigma_j^{1,N} - \gamma_j \hat{X}_0^N(\phi)| \leq \eta(N) \|\phi\|_{\text{Lip}}^{m(j)} \hat{X}_0^N(\mathbf{1}), \tag{7.19}$$

and that the remaining $\Sigma_j^{i,N}$ are bounded above by terms of the form given in the right-hand side of (7.12). This will prove (7.12). Note that it is (7.19) which identifies the parameters of the limiting super-Brownian motion of Theorem 1.2.

The $\Sigma_j^{i,N}$ are given by

$$\Sigma_1^{1,N} = \frac{1}{N'} \sum_x \phi(x)^2 \log NP(\hat{\xi}_0^N(\hat{B}_{t_N}^{N,x}) = 1, \hat{\tau}^N(x, x+e) > t_N),$$

$$\Sigma_1^{2,N} = \frac{1}{N'} \sum_x \phi(x)^2 \log NP(\hat{\xi}_0^N(\hat{B}_{t_N}^{N,x}) = \hat{\xi}_0^N(\hat{B}_{t_N}^{N,x+e}) = 1, \tag{7.20}$$

$$\hat{\tau}^N(x, x+e) > t_N),$$

$$\Sigma_2^{1,N} = \frac{1}{N'} \sum_x \phi(x) \log NP(\hat{\xi}_0^N(\hat{B}_{t_N}^{N,x+e}) = 1,$$
$$\hat{\tau}^N(x, x+e) \wedge \hat{\tau}^N(x, x+f) > t_N,$$
$$\hat{\tau}^N(x+e, x+f) \leq t_N),$$

$$\Sigma_2^{2,N} = \frac{1}{N'} \sum_x \phi(x) \log NP(\hat{\xi}_0^N(\hat{B}_{t_N}^{N,x}) = \hat{\xi}_0^N(\hat{B}_{t_N}^{N,x+e}) = 1,$$

$$\hat{\tau}^N(x, x+e) \wedge \hat{\tau}^N(x, x+f) > t_N, \tag{7.21}$$
$$\hat{\tau}^N(x+e, x+f) \leq t_N),$$

$$\Sigma_2^{3,N} = \frac{1}{N'} \sum_x \phi(x) \log NP(\xi_0^N(\hat{B}_{t_N}^{N,x}) = 0,$$
$$\hat{\xi}_0^N(\hat{B}_{t_N}^{N,x+e}) = \hat{\xi}_0^N(\hat{B}_{t_N}^{N,x+f}) = 1,$$
$$\hat{\tau}^N(x, x+e) \wedge \hat{\tau}^N(x, x+f)$$
$$\wedge \hat{\tau}^N(x+e, x+f) > t_N)$$

and

$$\Sigma_3^{1,N} = \frac{1}{N'} \sum_x \phi(x) \log NP(\hat{\xi}_0^N(\hat{B}_{t_N}^{N,x}) = 1,$$

$$\hat{\tau}^N(x, x+e) \wedge \hat{\tau}^N(x, x+f) > t_N,$$
$$\hat{\tau}^N(x+e, x+f) \leq t_N),$$

(7.22) $\Sigma_3^{2,N} = \Sigma_2^{2,N},$

$$\Sigma_3^{3,N} = \frac{1}{N'} \sum_x \phi(x) \log NP(\hat{\xi}_0^N(\hat{B}_{t_N}^{N,x}) = 1,$$

$$\hat{\xi}_0^N(\hat{B}_{t_N}^{N,x+e}) = \hat{\xi}_0^N(\hat{B}_{t_N}^{N,x+f}) = 0,$$
$$\hat{\tau}^N(x, x+e) \wedge \hat{\tau}^N(x, x+f)$$
$$\wedge \hat{\tau}^N(x+e, x+f) > t_N).$$

Although the expressions for the $\Sigma_j^{i,N}$ are lengthy, verification of (7.18) is a straightforward application of duality. We will prove the decomposition for $E(\Delta_3^{N,+}(\phi, \hat{\xi}_{t_N}^N))$, the others are proved similarly.

Using duality,

$$E(\Delta_3^{N,+}(\phi, \hat{\xi}_{t_N}^N))$$

$$= \frac{1}{N'} \sum_x \phi(x) \log NP(\hat{\xi}_{t_N}^N(x) = 1, \hat{\xi}_{t_N}^N(x+e) = \hat{\xi}_{t_N}^N(x+f) = 0)$$

$$= \frac{1}{N'} \sum_x \phi(x) \log NP(\hat{\xi}_0^N(\hat{B}_{t_N}^{N,x}) = 1, \hat{\xi}_0^N(\hat{B}_{t_N}^{N,x+e}) = \hat{\xi}_0^N(\hat{B}_{t_N}^{N,x+f}) = 0,$$

$$\hat{\tau}^N(x, x+e) \wedge \hat{\tau}^N(x, x+f) > t_N).$$

The possibility $\hat{\tau}^N(x+e, x+f) > t_N$ gives rise to $\Sigma_3^{3,N}$. For the complement, letting $E_N = \{\hat{\tau}^N(x, x+e) \wedge \hat{\tau}^N(x, x+f) > t_N, \hat{\tau}^N(x+e, x+f) \leq t_N\}$,

$$\{\hat{\xi}_0^N(\hat{B}_{t_N}^{N,x}) = 1, \hat{\xi}_0^N(\hat{B}_{t_N}^{N,x+e}) = \hat{\xi}_0^N(\hat{B}_{t_N}^{N,x+f}) = 0, E_N\}$$
$$= \{\hat{\xi}_0^N(\hat{B}_{t_N}^{N,x}) = 1, E_N\} \setminus \{\hat{\xi}_0^N(\hat{B}_{t_N}^{N,x}) = \hat{\xi}_0^N(\hat{B}_{t_N}^{N,x+e}) = 1, E_N\}.$$

Taking expectations, this gives the terms $\Sigma_3^{1,N}$ and $\Sigma_3^{2,N}$, proving the decomposition for $E(\Delta_3^{N,+})$ in (7.18).

Before tackling (7.19), let us dispense with the error terms $\Sigma_j^{i,N}$, $i \neq 1$. By Lemma 7.6, we see, that for $j = 1, 2, 3$,

(7.23) $$\Sigma_j^{2,N} \leq C_{7.14} \|\phi\|_\infty^{m(j)} \Big[ \varepsilon_N^{-1} \iint 1(|w-z| \leq \sqrt{\delta}_N) \, d\hat{X}_0^N(w) \, d\hat{X}_0^N(z)$$
$$+ (K_N \varepsilon_N)^{-1} \hat{X}_0^N(\mathbf{1})^2 \Big].$$

Next, consider $\Sigma_2^{3,N}$, and recall that $\tau^N(x,y) = \tau(x\sqrt{N}, y\sqrt{N})/v_N$. This term is no larger than

$$\frac{\|\phi\|_\infty}{N'} \sum_x \log NP(\hat{\xi}_0^N(\hat{B}_{t_N}^{N,x+e}) = 1,$$

$$\hat{\tau}^N(x, x+e) \wedge \hat{\tau}^N(x, x+f) \wedge \hat{\tau}^N(x+e, x+f) > t_N)$$

$$= \frac{\|\phi\|_\infty}{N'} \sum_w \hat{\xi}_0^N(w) \log N \sum_x P(B_{t_N}^{N,e} = w - x,$$

$$\tau^N(0,e) \wedge \tau^N(0,f) \wedge \tau^N(e,f) > t_N)$$

$$= \|\phi\|_\infty \hat{X}_0^N(\mathbf{1}) \log NP(\tau^N(0,e) \wedge \tau^N(0,f) \wedge \tau^N(e,f) > t_N).$$

Note that $\tau^N(e,f) > t_N$ implies $e \neq f$, and also that $\tau^N(x,y) > t_N$ is equivalent to $\tau(\sqrt{N}, x\sqrt{N}y) > v_N t_N$. Therefore, by Lemma 2.5, (2.7) and (2.4), we may conclude that, for a constant $C_{7.24}$ depending on $p(\cdot)$,

(7.24) $$\Sigma_2^{3,N} \leq C_{7.24} \|\phi\|_\infty \hat{X}_0^N(\mathbf{1})(\log N)^{-1/2}.$$

Virtually the same reasoning gives

(7.25) $$\Sigma_3^{3,N} \leq C_{7.24} \|\phi\|_\infty \hat{X}_0^N(\mathbf{1})(\log N)^{-1/2}.$$

On account of (7.23), (7.24) and (7.25), the proof of Proposition 7.5 will be complete once we establish (7.19).

Consider first the $j = 2$ case of (7.19). Then

$$|\Sigma_2^{1,N} - \gamma^* \hat{X}_0^N(\phi)|$$

$$= \left| \frac{1}{N'} \sum_w \hat{\xi}_0^N(w) \sum_x \phi(x) \log NP(\hat{B}_{t_N}^{N,e} = w - x, \right.$$

$$\hat{\tau}^N(0,e) \wedge \hat{\tau}^N(0,f) > t_N,$$

$$\left. \hat{\tau}^N(e,f) \leq t_N) - \gamma^* \hat{X}_0^N(\phi) \right|$$

$$= \left| \frac{1}{N'} \sum_w \hat{\xi}_0^N(w)(\log NE(\phi(w - \hat{B}_{t_N}^{N,e}) \right.$$

$$\times 1(\hat{\tau}^N(0,e) \wedge \hat{\tau}^N(0,f) > t_N,$$

$$\left. \hat{\tau}^N(e,f) \leq t_N) - \gamma^* \phi(w))) \right|.$$

Recalling the notation $q_T$ from (2.1), and using Cauchy–Schwarz in the second inequality below, we see that the above implies

$$|\Sigma_2^{1,N} - \gamma^* \hat{X}_0^N(\phi)|$$

$$\leq \frac{1}{N'} \sum_w \hat{\xi}_0^N(w) \log N$$
$$\times E(|\phi(w - \hat{B}_{t_N}^{N,e}) - \phi(w)|1(\hat{\tau}^N(0,e) \wedge \hat{\tau}^N(0,f) > t_N,$$
$$\hat{\tau}^N(e,f) \leq t_N))$$
$$+ \left| \frac{1}{N'} \sum_w \hat{\xi}_0^N(w)\phi(w)((\log N)q_{v_N t_N} - \gamma^*) \right|$$
$$\leq \hat{X}_0^N(\mathbf{1}) \log N E(|\phi|_{\mathrm{Lip}}^2 |B_{t_N}^{N,e}|^2)^{1/2} q_{v_N t_N}^{1/2}$$
$$+ \|\phi\|_\infty \hat{X}_0^N(\mathbf{1})|(\log N)q_{v_N t_N} - \gamma^*|$$
$$\leq |\phi|_{\mathrm{Lip}} \hat{X}_0^N(\mathbf{1}) \log N (2\sigma^2(N^{-1} + t_N)H(2v_N t_N))^{1/2}$$
$$+ \|\phi\|_\infty \hat{X}_0^N(\mathbf{1})|(\log N)q_{v_N t_N} - \gamma^*|.$$

Thus, by (2.4) and Proposition 2.1,

(7.26) $\quad |\Sigma_2^{1,N} - \gamma^* \hat{X}_0^N(\phi)| \leq \eta_{7.26}(N)\|\phi\|_{\mathrm{Lip}}\hat{X}_0^N(\mathbf{1}),$

where $\eta_{7.26}(N) \to 0$ as $N \to \infty$. Virtually the same argument gives the same bound for $|\Sigma_3^{1,N} - \gamma^* \hat{X}_0^N(\phi)|$.

Finally, arguing as we did for $\Sigma_2^{1,N}$, we have

$$|\Sigma_1^{1,N} - 2\pi\sigma^2 \hat{X}_0^N(\phi^2)|$$
$$= \left| \frac{1}{N'} \sum_w \hat{\xi}_0^N(w)[\log N E(\phi^2(w - B_{t_N}^{N,0})1(\tau^N(0,e) > t_N)) - 2\pi\sigma^2 \phi^2(w)] \right|$$
$$\leq \frac{1}{N'} \sum_w \hat{\xi}_0^N(w)[\log N E(|\phi^2(w - B_{t_N}^{N,0}) - \phi^2(w)|1(\tau^N(0,e) > t_N))]$$
$$+ \frac{1}{N'} \sum_w \hat{\xi}_0^N(w)\phi^2(w)|\log N P(\tau^N(0,e) > t_N) - 2\pi\sigma^2|$$
$$\leq \hat{X}_0^N(\mathbf{1})(\log N E((2\|\phi\|_\infty |\phi|_{\mathrm{Lip}}|B_{t_N}^{N,0}|)^2)^{1/2} H(2v_N t_N)^{1/2}$$
$$+ \|\phi\|_\infty^2 |(\log N)H(2v_N t_N) - 2\pi\sigma^2|).$$

Again, since $E(|B_{t_N}^{N,0}|^2) \leq 2\sigma^2 t_N$, (2.4) implies that

(7.27) $\quad |\Sigma_1^{1,N} - 2\pi\sigma^2 \hat{X}_0^N(\phi^2)| \leq \eta_{7.27}(N)\hat{X}_0^N(\mathbf{1})\|\phi\|_{\mathrm{Lip}}^2,$

where $\eta_{7.27}(N) \to 0$ as $N \to \infty$. The required result, (7.19), has been proved thanks to (7.26), its analogue for $j = 3$, and (7.27). This completes the proof of Proposition 7.5. $\square$

The next result will allow us to deduce Proposition 7.1 from Proposition 7.5.

LEMMA 7.8. *There is a constant $C_{7.28}$ so that, for $j = 1, 2, 3$, $\phi: \mathbb{R}^2 \to \mathbb{R}$, and all $0 \leq t \leq 1$,*

$$|E(\Delta_j^{N,+}(\phi, \xi_t^N)) - E(\Delta_j^{N,+}(\phi, \hat{\xi}_t^N))| \tag{7.28}$$
$$\leq C_{7.28} \|\phi\|_\infty^{m(j)} [(\log N)^{-1} + t \log N] X_0^N(\mathbf{1}),$$

*where $m(j)$ is as in (7.11).*

PROOF. The proofs are quite similar, so we only consider $j = 2$. If $\xi, \eta \in \{0,1\}^{\mathbf{S}_N}$ satisfy $\xi \leq \eta$, then (5.34) with $\phi = \delta_x$ implies

$$|(1 - \eta(x))f_1^N(x, \eta)^2 - (1 - \xi(x))f_1^N(x, \xi)^2|$$
$$\leq \eta(x) - \xi(x) + 2[f_1^N(x, \eta) - f_1^N(x, \xi)].$$

If we apply the above with $\eta = \bar{\xi}_t^N$ and $\xi = \xi_t^N$, we get

$$|E(\Delta_2^{N,+}(\phi, \bar{\xi}_t^N)) - E(\Delta_2^{N,+}(\phi, \xi_t^N))|$$
$$\leq \frac{1}{N'} \left| E\left[ \sum_x [(1 - \bar{\xi}_t^N(x))f_1^N(x, \bar{\xi}_t^N)^2 \right.\right.$$
$$\left.\left. - (1 - \xi_t^N(x))f_1^N(x, \xi_t^N)^2] \log N \phi(x) \right] \right|$$
$$\leq \|\phi\|_\infty \log N E\left[ \frac{1}{N'} \sum_x \bar{\xi}_t^N(x) - \xi_t^N(x) + 2(f_1^N(x, \bar{\xi}_t^N) - f_1^N(x, \xi_t^N)) \right]$$
$$\leq 3\|\phi\|_\infty \log N [E(\bar{X}_t^N(\mathbf{1})) - E(X_t^N(\mathbf{1}))]$$
$$\leq 3\|\phi\|_\infty C_{6.7}[(\log N)^{-1} + t \log N] X_0^N(\mathbf{1}) \qquad \text{for } 0 \leq t \leq 1,$$

applying the bound (6.7) in the last inequality above. This argument gives the same bound for $|E(\Delta_2^{N,+}(\phi, \bar{\xi}_t^N)) - E(\Delta_2^{N,+}(\phi, \hat{\xi}_t^N))|$, and (7.28) follows by the triangle inequality. □

7.3. *Proof of Proposition 7.1.* Let $m(j)$ be as in (7.11). By writing $\phi = \phi^+ - \phi^-$ if $m(j) = 1$, we may assume $\phi^{m(j)} \geq 0$. Combine Proposition 7.5 and Lemma 7.8 to conclude that

$$|E(\Delta_j^N(\phi, \xi_{t_N}^N))| \leq \eta_{7.12}(N)(X_0^N(\mathbf{1}) + X_0^N(\mathbf{1})^2)\|\phi\|_{\text{Lip}}^{m(j)}$$
$$+ \frac{C_{7.12}\|\phi\|_\infty^{m(j)}}{\varepsilon_N} \iint \mathbf{1}(|w - z| \leq \sqrt{\delta_N})\, dX_0^N(w)\, dX_0^N(z)$$

$$+ C_{7.28}\|\phi\|_\infty^{m(j)}[(\log N)^{-1} + t_N \log N]X_0^N(\mathbf{1}),$$
$$+ \gamma_j |E(X_{t_N}^N(\phi^{m(j)}) - \hat{X}_{t_N}^N(\phi^{m(j)}))|.$$

By our coupling and (6.7),

$$|E(X_{t_N}^N(\phi^{m(j)}) - \hat{X}_{t_N}^N(\phi^{m(j)}))| \le E(\bar{X}_{t_N}^N(\phi^{m(j)}) - \underline{X}_{t_N}^N(\phi^{m(j)}))$$
$$\le \|\phi\|_\infty^{m(j)} E(\bar{X}_{t_N}^N(\mathbf{1}) - \underline{X}_{t_N}^N(\mathbf{1}))$$
$$\le C_{6.7}\|\phi\|_\infty^{m(j)}[(\log N)^{-2} + t_N]X_0^N(\mathbf{1}).$$

If we insert the above into the previous inequality we obtain the required upper bound (7.4).

7.4. *Proof of Proposition 7.2.* We will work with the biased voter model and prove an analogous result for it. Namely, there is constant $C_{7.5}$ such that

$$E\left(\iint 1(|x-y| \le \sqrt{\delta_N}) d\bar{X}_t^N(x) d\bar{X}_t^N(y)\right)$$
(7.29)
$$\le C_{7.5} e^{C_{7.5} t}(\bar{X}_0^N(\mathbf{1})^2 + \bar{X}_0^N(\mathbf{1}))$$
$$\times \left[\frac{\delta_N}{\delta_N + t}(1 + t^{2/3}) + \delta_N t^{-1/3} \log\left(1 + \frac{t}{\delta_N}\right)\right].$$

Since $\xi_t^N \le \bar{\xi}_t^N$, (7.29) would imply (7.5).

To prove (7.29), we need two estimates which are simple consequences of Proposition 7.3. We express these estimates in terms of the random walk $B_\cdot^{N,*}$ which takes steps according to $p_N(\cdot)$ at rate $v_N + b_N = N + \bar{\theta} \log N$ and has semigroup $P_t^{N,*}$, introduced just before Proposition 4.4.

COROLLARY 7.9. (a) *For all $x \in \mathbf{S}_N$ and $t \ge 0$,*

(7.30) $$P(B_t^{N,*} = x) \le \frac{C_{7.6}}{1 + Nt}.$$

(b) *Assume $\delta_N' \downarrow 0$ and $N\delta_N' \to \infty$. For each $K > 0$, there is a constant $C_{7.31}(K) > 0$ so that*

(7.31) $$\inf_{N\ge 1, w \in \mathbf{S}_N, |w| \le K\sqrt{\delta_N'}} N\delta_N' P(B_{2\delta_N'}^{N,*} = w) \ge C_{7.31}(K).$$

PROOF. Inequality (7.30) is a direct consequence of (7.6). The inequality (7.31) is a direct consequence of (7.9) and the fact that $P(B_t^0 = x) > 0$ for all $x \in \mathbb{Z}^2, t > 0$. □

For $x, w \in \mathbf{S}_N$, let $p_t^{N,w}(x) = NP(B_t^{N,*} = x - w)$ and $p_t^N(x) = p_t^{N,0}(x)$. For fixed $\varepsilon, t > 0$, let $\phi_s^z(x) = p_{\varepsilon+t-s}^{N,z}(x)$. If $\bar{M}$ is as in Lemma 5.5 and $\phi: \mathbb{R}_+ \times \mathbf{S}_N \to \mathbb{R}$, let $\phi_N(s, x) = \phi(s, x/\sqrt{N}), x \in \mathbb{Z}^2$, and $\bar{M}_t^N(\phi) = (N')^{-1} \bar{M}_t(\phi_N)$. By (5.36) and integration by parts,

$$\bar{X}_t^N(p_\varepsilon^{N,z})^2 = \bar{X}_0^N(p_{\varepsilon+t}^{N,z})^2 + 2 \int_0^t \bar{X}_{s-}^N(\phi_s^z) \, d\bar{M}_s^N(\phi^z)$$

$$+ 2 \int_0^t \bar{X}_s^N(p_{\varepsilon+t-s}^{N,z}) \bar{X}_s^N(p_{\varepsilon+t-s}^{N,z} b_N f_0^N(\cdot, \bar{\xi}_s^N)) \, ds$$

$$+ [\bar{M}^N(\phi^z)]_t.$$

Then $E([\bar{M}^n(\phi^z)]_t) = N'^{-2} E(\langle \bar{M}(\phi^z)_N\rangle_t)$ and so we may take means and use Lemma 5.5(a) to conclude

$$E(\bar{X}_t^N(p_\varepsilon^{N,z})^2)$$

(7.32)
$$\leq \bar{X}_0^N(p_{\varepsilon+t}^{N,z})^2 + 2E\left( \int_0^t \bar{X}_s^N(p_{\varepsilon+t-s}^{N,z}) \bar{X}_s^N(p_{\varepsilon+t-s}^{N,z} b_N f_0^N(\bar{\xi}_s^N)) \, ds \right)$$

$$+ E\left( \int_0^t C_{7.32} \frac{\log N}{N'} \sum_x p_{\varepsilon+t-s}^{N,z}(x)^2 [\bar{\xi}_s^N(x) f_0^N(x, \bar{\xi}_s^N) \right.$$

$$\left. + (1 - \bar{\xi}_s^N(x)) f_1^N(x, \bar{\xi}_s^N) ] \, ds \right).$$

Sum over $z$, multiply by $\varepsilon/N$, and use $\frac{1}{N} \sum_z p_s^{N,z}(x) p_s^{N,z}(y) = p_{2s}^N(y - x)$ to see that

$$E\left( \int \int \varepsilon p_{2\varepsilon}^N(y - x) \, d\bar{X}_t^N(x) \, d\bar{X}_t^N(y) \right)$$

$$\leq \int \int \varepsilon p_{2(\varepsilon+t)}^N(y - x) \, dX_0^N(x) \, dX_0^N(y)$$

$$+ 2E\left( \int_0^t \int \int \varepsilon p_{2(\varepsilon+t-s)}^N(y - x) b_N f_0^N(y, \bar{\xi}_s^N) \, d\bar{X}_s^N(x) \, d\bar{X}_s^N(y) \, ds \right)$$

$$+ C_{7.32} \int_0^t \varepsilon p_{2(\varepsilon+t-s)}^N(0) E\left( \frac{\log N}{N'} \sum_x [\bar{\xi}_s^N(x) f_0^N(x, \bar{\xi}_s^N) \right.$$

$$\left. + (1 - \bar{\xi}_s^N(x)) f_1^N(x, \bar{\xi}_s^N)] \right) ds$$

$$\leq \frac{C_{7.6} \varepsilon N}{1 + 2N(\varepsilon + t)} X_0^N(\mathbf{1})^2$$

$$+ \int_0^t \frac{2C_{7.6} \varepsilon N}{1 + 2N(\varepsilon + t - s)} E(\bar{X}_s^N(\mathbf{1}) b_N \bar{X}_s^N(f_0^N(\bar{\xi}_s^N))) \, ds$$

$$+ 2C_{7.32}C_{7.6} \int_0^t \frac{\varepsilon N}{1+2N(\varepsilon+t-s)} E(\log N \bar{X}_s^N(f_0^N(\bar{\xi}_s^N))) \, ds,$$

where we have used (7.30) in the last line. Next use the mass bounds (6.5) and (6.6) and the definition of $g$ in (4.1) to bound the sum of the last two integrals by

$$Ce^{C_{4.1}t} \int_0^t \frac{\varepsilon}{\varepsilon+t-s} s^{-1/3} \, ds \, [X_0^N(\mathbf{1})^2 + X_0^N(\mathbf{1})].$$

Consequently, there is a constant $C_{7.33}$ such that, for $0 \le t \le T$,

(7.33)
$$E\left(\int\int \varepsilon p_{2\varepsilon}^N(y-x) \, d\bar{X}_t^N(x) \, d\bar{X}_t^N(y)\right)$$
$$\le C_{7.33} e^{C_{7.33}T}[X_0^N(\mathbf{1})^2 + X_0^N(\mathbf{1})]\left(\frac{\varepsilon}{\varepsilon+t} + \int_0^t \frac{\varepsilon}{\varepsilon+t-s} s^{-1/3} \, ds\right)$$
$$\le C_{7.33} e^{C_{7.33}T}[X_0^N(\mathbf{1})^2 + X_0^N(\mathbf{1})]\left(\frac{\varepsilon}{\varepsilon+t}(1+3t^{2/3})\right.$$
$$\left. + 2\varepsilon t^{-1/3}\log(1+t/\varepsilon)\right)$$

(split the integral at $t/2$). Now set $\varepsilon = \delta_N$ and note that the above upper bound is the right-hand side of (7.29). The left-hand side of (7.33) is bounded below by

$$E\left(\int\int \delta_N p_{2\delta_N}(y-x) 1(|y-x| \le \sqrt{\delta_N}) \, d\bar{X}_t^N(x) \, d\bar{X}_t^N(y)\right)$$
$$\ge C_{7.31}(1) E(1(|y-x| \le \sqrt{\delta_N}) \, d\bar{X}_t^N(x) \, d\bar{X}_t^N(y)),$$

where $C_{7.31}(1) > 0$ by (7.31). Combine this with (7.33) to complete the proof of (7.29), and hence, Proposition 7.2.

**8. Proof of Proposition 4.7—Part II.** We will make frequent use of the following elementary estimate. If $\phi: \mathbb{R}^d \to \mathbb{R}$ such that $\|\phi\|_\infty \le K$, then there is a constant $C_{8.1} = C_{8.1}(K)$ such that, for $j = 1, 2, 3$ and $s > 0$,

(8.1) $\quad |\Delta_j^N(\phi, \xi_s^N)| \le C_{8.1}(K)[X_s^N(\log N f_0^N(\cdot, \xi_s^N)) + X_s^N(\mathbf{1})].$

Now let $J \ge 1$, $T_J^N = \inf\{s : X_s^N(\mathbf{1}) \ge J\}$, $t \le T$, $1 < p < 2$, and recall the sequences defined in (7.3). Assume also that $\|\phi\|_{\text{Lip}} \vee X_0^N(\mathbf{1}) \le K$. We first show there is a constant $C_{8.2}(T, K, p)$ such that, for $j = 1, 2, 3$,

(8.2) $\quad E\left(\left|\int_0^t 1(s \ge T_J^N + t_N) \Delta_j^N(\phi, \xi_s^N) \, ds\right|^p\right) \le \frac{C_{8.2}(T, K, p)}{J^{2-p}}.$

By Hölder's inequality, the left-hand side is at most

$$P(T_J^N \leq t)^{(2-p)/2} E\bigg(\bigg(\int_0^t \Delta_j^N(\phi, \xi_s^N)\, ds\bigg)^2\bigg)^{p/2}$$

$$\leq 2^{p/2} P\bigg(\sup_{s \leq t} X_s^N(\mathbf{1}) > J\bigg)^{(2-p)/2}$$

$$\times \bigg[E\bigg(\int_0^t \Delta_j^N(\phi, \xi_{s_1}^N)\, ds_1 \int_{s_1}^t \Delta_j^N(\phi, \xi_{s_2}^N)\, ds_2\bigg)\bigg]^{p/2}$$

$$\leq 2^{p/2} C_{8.1}^p \frac{E(\sup_{s \leq t} X_s^N(\mathbf{1})^2)^{(2-p)/2}}{J^{2-p}}$$

$$\times \bigg[E\bigg(\int_0^t (X_{s_1}^N(\log N f_0^N(\xi_{s_1}^N)) + X_{s_1}^N(\mathbf{1}))\, ds_1$$

$$\times \int_{s_1}^t E_{X_{s_1}^N}(X_{s_2-s_1}^N(\log N f_0^N(\xi_{s_2-s_1}^N)) + X_{s_2-s_1}^N(\mathbf{1}))\, ds_2\bigg)\bigg]^{p/2}$$

$$\leq \frac{C(T, K, p)}{J^{2-p}} \bigg[E\bigg(\int_0^T (X_{s_1}^N(\log N f_0^N(\xi_{s_1}^N)) + X_{s_1}^N(\mathbf{1}))$$

$$\times \int_{s_1}^T [(g(s_2 - s_1) + 1) X_{s_1}^N(\mathbf{1})]\, ds_2\, ds_1\bigg)\bigg]^{p/2}$$

for a constant $C(T, K, p)$. In the next to last line we used (8.1) and in the last line we used (4.2), (4.3) and (4.4). After simplification and applying (4.5) and (4.3), we have that the last line above is at most

$$\frac{C'(T, K, p)}{J^{2-p}} E\bigg(\int_0^T g(s_1)(X_0^N(\mathbf{1})^2 + X_0^N(\mathbf{1}))\, ds_1\bigg)^{p/2}$$

for a constant $C'(T, K, p)$. This proves (8.2).

We can now use $L^2$ estimates. Let

$$(8.3) \quad E\bigg(\bigg[\int_0^t 1(s \leq T_J^N + t_N) \Delta_j^N(\phi, \xi_s^N)\, ds\bigg]^2\bigg) = I_1(N, J, t) + I_2(N, J, t),$$

where

$$I_1(N, J, t) = 2E\bigg(\int_0^t 1(s_1 \leq T_J^N + t_n) \Delta_j^N(\phi, \xi_{s_1}^N)$$

$$\times \int_{s_1}^{t \wedge (s_1 + t_N)} 1(s_2 \leq T_J^N + t_N) \Delta_j^N(\phi, \xi_{s_2}^N)\, ds_1\, ds_2\bigg)$$

$$I_2(N, J, t) = 2E\bigg(\int_0^t 1(s_1 \leq T_J^N + t_n) \Delta_j^N(\phi, \xi_{s_1}^N)$$

$$\times \int_{s_1}^t 1(s_1 + t_N \leq s_2 \leq T_J^N + t_N) \Delta_j^N(\phi, \xi_{s_2}^N)\, ds_2\, ds_1\bigg).$$

By (8.1), the Markov property and (4.2) and (4.4), it follows that $I_1(N, J, t)$ is at most

$$2C_{8.1}^2(K)E\bigg(\int_0^T [X_{s_1}^N(\log N f_0^N(\xi_{s_1}^N)) + X_{s_1}^N(\mathbf{1})]$$
$$\times \int_{s_1}^{s_1+t_N} E_{X_{s_1}^N}(X_{s_2-s_1}^N(\log N f_0^N(\xi_{s_2-s_1}^N)) + X_{s_2-s_1}^N(\mathbf{1}))\, ds_2 \, ds_1\bigg)$$
$$\leq 2C_{8.1}^2 E\bigg(\int_0^T [X_{s_1}^N(\log N f_0^N(\xi_{s_1}^N)) + X_{s_1}^N(\mathbf{1})]$$
$$\times \int_{s_1}^{s_1+t_N} (g(s_2 - s_1) + C_{4.2}(T))X_{s_1}^N(\mathbf{1})\, ds_2 \, ds_1\bigg).$$

In view of the definition of $g(s)$, (4.3) and (4.5), it follows that there is a constant $C_{8.4}(K,T)$ such that

$$(8.4) \qquad I_1(N, J, T) \leq C_{8.4}(K,T) t_N^{2/3}.$$

Turning to $I_2(N, J, T)$, we may use Proposition 7.1 and (8.1) to see that, for $s_1 + t_N \leq s_2$,

$$|E(1(s_1 \leq T_J^N)1(s_2 \leq T_J^N + t_N)\Delta_j^N(\phi, \xi_{s_1}^N)\Delta_j^N(\phi, \xi_{s_2}^N))|$$
$$\leq E(1(s_1 < T_J^N)1(s_2 \leq T_J^N + t_N)|\Delta_j^N(\phi, \xi_{s_1}^N)||E_{X_{s_2-t_N}^N}(\Delta_j^N(\phi, \xi_{t_N}^N))|)$$
$$\leq C_{8.1}(K)$$
$$\times E\bigg(1(s_1 < T_J^N)[X_{s_1}^N(\log N f_0^N(\xi_{s_1}^N)) + X_{s_1}^N(\mathbf{1})]$$
$$\times C_{7.4}\Big[\eta_{7.4}(N)(X_{s_2-t_N}^N(\mathbf{1})^2 + X_{s_2-t_N}^N(\mathbf{1}))$$
$$+ (\varepsilon_N)^{-1}$$
$$\times \int\int 1(|w-z| \leq \sqrt{\delta_N})\, dX_{s_2-t_N}^N(z)\, dX_{s_2-t_N}^N(w)\Big]\bigg).$$

Since $1\{T_J^N < s_1 < T_J^N + t_N\}1\{s_1 + t_N < s_2 < T_J^N + t_N\} = 0$, this implies that

$$I_2(N, J, t)$$
$$\leq 2C_{8.1}(K)$$
$$\times \int_0^T E\bigg[1(s_1 < T_J^N)[X_{s_1}^N(\log N f_0^N(\xi_{s_1}^N)) + X_{s_1}^N(\mathbf{1})]$$
$$\times C_{7.4}\bigg(\eta_{7.4}(N)E_{X_{s_1}^N}\bigg(\int_0^{T-s_1}(X_{s_2}^N(\mathbf{1})^2 + X_{s_2}^N(\mathbf{1}))\, ds_2\bigg)$$

$$+ (\varepsilon_N)^{-1}$$
$$\times E_{X_{s_1}^N}\left(\int_0^{T-s_1}\int\int 1(|w-z| \leq \sqrt{\delta_N})\, dX_{s_2}^N(z)\, dX_{s_2}^N(w)\, ds_2\right)\bigg)\bigg]\, ds_1.$$

Apply (4.3) and Proposition 7.2 to see that there is a constant $C(T, K)$ such that $I_2(N, J, T)$ is bounded above by

$$C(T, K) \int_0^T E\bigg(1(s_1 < T_J^N)[X_{s_1}^N(\log N f_0^N(\xi_{s_1}^N)) + X_{s_1}^N(\mathbf{1})]$$
$$\times \bigg\{\eta_{7.4}(N)(X_{s_1}^N(\mathbf{1})^2 + X_{s_1}^N(\mathbf{1}))$$
$$+ \varepsilon_N^{-1} C_{7.5}(X_{s_1}^N(\mathbf{1})^2 + X_{s_1}^N(\mathbf{1}))$$
$$\times \bigg[\int_0^{T-s_1}\bigg(\frac{\delta_N}{\delta_N + s_2}(1 + s_2^{2/3})$$
$$+ \delta_N s_2^{-1/3}\log\bigg(1 + \frac{s_2}{\delta_N}\bigg)\bigg)\, ds_2\bigg]\bigg\}\bigg)\, ds_1$$
$$\leq C(T,K)(J^2 + J)$$
$$\times \int_0^T E(X_{s_1}^N(\log N f_0^N(\xi_{s_1}^N)) + X_{s_1}^N(\mathbf{1}))$$
$$\times \bigg\{\eta_{7.4}(N)$$
$$+ (\delta_N/\varepsilon_N)\int_0^{T-s_1}\bigg(\frac{1 + T^{2/3}}{\delta_N + s_2} + s_2^{-1/3}\log\bigg(1 + \frac{T}{\delta_N}\bigg)\bigg)\, ds_2\bigg\}\, ds_1.$$

By using (4.2) and (4.4) and evaluating the remaining deterministic integrals, we see that there is a constant $C_{8.5}(K, T)$ such that

(8.5) $\qquad I_2(N, J, T) \leq C_{8.5}(K, T)(J^2 + J)\eta_{8.5}(N),$

where $\eta_{8.5}(N) = \frac{\delta_N}{\varepsilon_N}\log(1 + \frac{T}{\delta_N}) + \eta_{7.4}(N) \to 0$ as $N \to \infty$ [recall (7.3)].

By standard inequalities,

$$E\bigg(\bigg|\int_0^t \Delta_j^N(\phi, s, \xi^N)\, ds\bigg|^p\bigg) \leq \bigg(E\bigg(\bigg|\int_0^t \Delta_j^N(\phi, s, \xi^N)\, ds\bigg|^2\bigg)\bigg)^{p/2}$$
$$\leq (I_1(N,T,K)^{1/2} + I_2(N,T,K)^{1/2})^p.$$

We now choose $J = J_N \to \infty$ such that $J_N^2 \eta_{8.5}(N) \to 0$ as $N \to \infty$. Then (7.2) follows from (8.2), (8.3), (8.4), (8.5) and the last inequality. As was noted at the beginning of Section 7, Proposition 4.7 is then immediate.

**9. Proof of Theorem 1.3.** The proof of survival in [3] was given for general voter model perturbations assuming $d \geq 3$ and $N' \equiv N$. Here we are concerned only with Lotka–Volterra models, but are working in dimension $d = 2$ with mass normalization $N' = N/\log N$. In this section we will state and prove analogues of Lemma 3.2 and Proposition 4.2 of [3]. Given these results, the argument in Section 5 of [3] applies without further change to complete the proof of Theorem 1.3. We will content ourselves with proving survival only, and not derive a lower bound on the probability of survival as given in Corollary 3 in [3].

For any $K > 2$, $L > 1$ and $N \in \mathbb{N}$, define $I = [-L, L]^2$, $I_1 = (2L, 0) + I$, $I_{-1} = (-2L, 0) + I$, $I' = (-KL, KL)^2$, and $I'_N = (-KL\sqrt{N}, KL\sqrt{N})^2$. For $\tilde{\xi}_0 \in \{0,1\}^{\mathbb{Z}^2}$, supported on $I'_N$, let $\{\tilde{\xi}_t(x) : x \in \mathbb{Z}^2, t \geq 0\}$ be the Lotka–Volterra model where all sites $x \notin I'_N$ are set to 0 for all time. We may construct $\xi$. and $\tilde{\xi}$. as the solutions of a stochastic differential equation as in [3] [see (SDE)(I′) in Proposition 2.1 of that paper] so that if $|\xi_0| < \infty$ and $\tilde{\xi}_0 \leq \xi_0$, which we will assume throughout, then $\tilde{\xi}_t \leq \xi_t$ for all $t \geq 0$. For $x \in \mathbf{S}_N$, let $\tilde{\xi}_t^N(x) = \tilde{\xi}_{Nt}(x\sqrt{N})$, $\tilde{X}_t^N = \frac{1}{N'}\sum_{x \in \mathbf{S}_N} \tilde{\xi}_t^N(x)\delta_x$, and $\xi_\cdot^N$, $X_\cdot^N$ be as usual. Note that $\tilde{\xi}_\cdot^N$ and $\tilde{X}_\cdot^N$ are supported on $I'$.

The main technical step in the proof of Theorem 1.3 is the following version of Lemma 3.2 of [3]. Among our standing assumptions (1.8)–(1.10), it only requires (1.8)(a) and (1.10). Let $\|\cdot\|$ be the sup norm on $\mathbb{R}^2$. Recalling the independent random walk family $B_t^x$ introduced just before Theorem 1.1, we define here

$$\beta_t^{N,x} = B_{tN}^{x\sqrt{N}}, \qquad x \in \mathbf{S}_N. \tag{9.1}$$

Also define

$$\delta'_N = \frac{\bar{\theta}\log N}{N}, \tag{9.2}$$

from the following:

LEMMA 9.1. *There exists a nondecreasing $C_{9.3}: \mathbb{R}_+ \to \mathbb{R}_+$, depending only on $\bar{\theta}$ and $p(\cdot)$, such that, for any $N \in \mathbb{N}$, $t \geq 0$, $K > 2$, $L > 1$, if $X_0^N = \tilde{X}_0^N$ is supported by $I$, then*

$$\begin{aligned}E(X_t^N(\mathbf{1}) - \tilde{X}_t^N(\mathbf{1})) \\ \leq C_{9.3}(t)X_0^N(\mathbf{1})\bigg[\frac{KL}{\log N} + P\bigg(\sup_{u \leq t(1+\delta'_N)} \|\beta_u^{N,0}\| > (K-2)L\bigg)\bigg].\end{aligned} \tag{9.3}$$

Given this bound, the next step is the following analogue of Proposition 4.2 in [3]. Recall the definition of $S^\eta$ before Theorem 1.3. For $\alpha = (\alpha_0, \alpha_1)$, let $\|\alpha\|_1 = |\alpha_0 - 1| + |\alpha_1 - 1|$, and for $K \geq 1$, let $\gamma_K = 6^{-4(2K+1)^2}$.

PROPOSITION 9.2. *Assume $0 < \eta < 1$. There are $L, K, J \in \mathbb{N}$, $T \geq 1$ and $r \in (0, e^{-4})$ depending on $\eta$ such that if*

(9.4) $\quad \alpha \in S^\eta, \qquad \|\alpha\|_1 < r \quad \text{and} \quad N = N(\alpha) = \left\lfloor \left( \frac{\log(1/\|\alpha\|_1)}{\|\alpha\|_1} \right)^{1/2} \right\rfloor^2,$

*then $X_0^N(I) = X_0^N(\mathbf{1}) \geq J$ implies*

(9.5) $\qquad\qquad P^\alpha(\tilde{X}_T^N(I_1) \wedge \tilde{X}_T^N(I_{-1}) \geq J) \geq 1 - \gamma_K.$

Inequality (9.5) is just what is needed to show that the Lotka–Volterra process "dominates" a super-critical oriented percolation, and hence, survives. The details of this argument are spelled out in Section 5 of [3], and apply without change to the current setting. Therefore, to prove Theorem 1.3, it suffices now to prove Lemma 9.1 and Proposition 9.2. We will start with the proof of the second result, assuming the validity of the first one. Our argument closely follows the proofs of Proposition 4.2 and Theorem 8.3 in [3].

PROOF OF PROPOSITION 9.2. We now choose certain constants which depend only on $\eta > 0$ and $p(\cdot)$. Let $c = c(\sigma) \geq 1$ be large enough so that

(9.6) $\qquad \exp(-c^2 K^2 / 17\sigma^2) < \gamma_K / 4 \qquad \text{for all } K \geq 1.$

As in Lemma 4.3 of [3], we may choose $T = T(\eta) > 2$ and $L = c\sqrt{T} \in \mathbb{N}$ such that if $X_t$ is super-Brownian motion with branching rate $2\gamma_e$, diffusion rate $\sigma^2$ and drift $d_0 \in [\eta\gamma^*/24, \gamma^*]$ [recall $\gamma^*$ from (1.12)], then there is a constant $C_{9.7} = C_{9.7}(\eta)$ such that

(9.7) $\qquad P(X_T(I_1) \wedge X_T(I_{-1}) \leq 3X_0(I)) \leq C_{9.7}/X_0(I).$

Next, let $K > 4 + \frac{4\sigma}{c}$ be large enough so that

(9.8) $\qquad 8 C_{9.3}(T) e^{-c^2 K^2 / 16\sigma^2} < e^{-c^2 K^2 / 17\sigma^2},$

and let $J \in \mathbb{N}$ be large enough so that

(9.9) $\qquad\qquad C_{9.7}/J < e^{-c^2 K^2 / 17\sigma^2}.$

By monotonicity of $\tilde{X}^N$ (Proposition 2.1 (b)(ii) of [3]), we may assume $X_0^N(I) = X_0^N(\mathbf{1}) = J$. We claim that, with $c, T, L$ and $K$ defined above, there exists $r \in (0, e^{-4})$ such that if $\alpha \in S^\eta$, $\|\alpha\|_1 < r$, and $N = N(\alpha)$ is defined as in (9.4), then

(9.10) $\qquad P^\alpha(X_T^N(I_1) \wedge X_T^N(I_{-1}) \leq 3J) \leq 2e^{-c^2 K^2 / 17\sigma^2}$

and

(9.11) $\qquad P^\alpha(X_T^N(I_j) - \tilde{X}_T^N(I_j) > J) \leq e^{-c^2 K^2 / 17\sigma^2}, \qquad j = \pm 1.$

Given these estimates, it is straightforward using (9.6) to complete the proof of (9.5).

Before beginning the proofs of (9.10) and (9.11), we note that

$$\alpha_0 - \alpha_1 \geq \frac{\eta}{3}\|\alpha\|_1 \quad \text{for all } \alpha \in S^\eta \tag{9.12}$$

and also that for $\alpha \in S_\eta$ and $N' = N'(\alpha) = N(\alpha)/\log N(\alpha)$,

$$N'\|\alpha\|_1 \in [\tfrac{1}{8}, 1]. \tag{9.13}$$

These estimates follow from (crude) elementary calculations which we omit. (9.12) uses $\eta < 1$ and (9.13) uses $\|\alpha\|_1 < r < e^{-4}$.

If (9.10) fails, then we may suppose there exists a sequence $\alpha^m \in S^\eta$ and initial states $X_0^{N_m}$ supported on $I$ with $X_0^{N_m}(I) = J$ such that $\|\alpha^m\|_1 \to 0$ as $m \to \infty$, and

$$P^{\alpha^m}(X_T^{N_m}(I_1) \wedge X_T^{N_m}(I_{-1}) \leq 3J) > 2e^{-c^2 K^2/17\sigma^2} \quad \text{for all } m. \tag{9.14}$$

We may assume, by taking an appropriate subsequence, that $X_0^{N_m} \to X_0$ for some $X_0 \in \mathcal{M}_f$ supported on $I$ satisfying $X_0(I) = J$. We may also assume, in view of (9.13), that $N'_m(\alpha_0^m - 1, \alpha_1^m - 1) \to (\theta_0, \theta_1)$ for some $(\theta_0, \theta_1) \in \mathbb{R}^2$. The inequalities (9.12) and (9.13), and the fact that $|\alpha_0 - \alpha_1| \leq \|\alpha\|_1$ imply that

$$\gamma^*(\theta_0 - \theta_1) = \lim_{m \to \infty} \gamma^* N'_m(\alpha_0^m - \alpha_1^m) \in [\eta\gamma^*/24, \gamma^*].$$

Now let $X_t$ denote super-Brownian motion with branching rate $2\gamma_e$, diffusion rate $\sigma^2$ and drift $\gamma^*(\theta_0 - \theta_1)$. By Theorem 1.2, the fact that $X_T(\partial I_1) = 0$ a.s., and the inequalities (9.7) and (9.9), it follows that

$$\limsup_{m \to \infty} P^{\alpha^{N_m}}(X_T^{N_m}(I_1) \leq 3J) \leq P_{X_0}(X_T(I_1) \leq 3J)$$

$$< C_{9.7}/J \leq e^{-c^2 K^2/17\sigma^2}.$$

Since the same estimate is valid for $I_{-1}$, (9.14) cannot hold. This proves (9.10).

Now consider (9.11) for $I_1$. Following the proof of (4.11) in [3] (the odd lower bound on $K$ is used here), the invariance principle for Brownian motion and a standard Gaussian estimate imply there exists a sequence $\varepsilon_N \to 0$, depending ultimately on $\eta$, $\sigma^2$ and $\bar{\theta}$, such that

$$P\left(\sup_{u \leq T(1+\delta'_N)} \|\beta_u^{N,0}\| > (K-2)L\right) \leq 8\exp(-c^2 K^2/16\sigma^2) + \varepsilon_N.$$

Combining this estimate with (9.3), the coupling $\tilde{\xi}_t^N \leq \xi_t^N$ and Chebyshev's inequality imply

$$P(X_T^N(I_1) - \tilde{X}_T^N(I_1) > J) \tag{9.15}$$
$$\leq C_{9.3}(T)X_0^N(I)\left[\frac{KL}{\log N} + 8e^{-c^2 K^2/16\sigma^2} + \varepsilon_N\right]\Big/J.$$

If (9.11) fails for $I_1$, then we may assume the existence of sequences $\alpha^{N_m}$ and $X_0^{N_m}$ as before, such that, for all $m$,

(9.16) $\quad P^{\alpha^{N_m}}(X_T^{N_m}(I_1) - \tilde{X}_T^{N_m}(I_1) > J) > \exp(-c^2 K^2/17\sigma^2).$

However, the estimate (9.15) implies that [recall $X_0^N(I) = J$]

$$\limsup_{m \to \infty} P^{\alpha^{N_m}}(X_T^{N_m}(I_1) - \tilde{X}_T^{N_m}(I_1) > J) \leq 8C_{9.3}(T)e^{-c^2K^2/16\sigma^2}$$

$$< \exp(-c^2 K^2/17\sigma^2),$$

by (9.8), and this contradicts (9.16). □

To prove Lemma 9.1, we must also work with the rescaled voter and biased voter models $\hat{\xi}_\cdot^N$, $\bar{\xi}_\cdot^N$, $\hat{X}_\cdot^N$, $\bar{X}_\cdot^N$ from Sections 5 and 6 with rates and bias $v_N$ and $b_N$, respectively, as in (5.15), as well as their counterparts with 0 boundary values off of $I'$, $\tilde{\hat{\xi}}_\cdot^N$, $\tilde{\bar{\xi}}_\cdot^N$, $\tilde{\hat{X}}_\cdot^N$ and $\tilde{\bar{X}}_\cdot^N$. We will assume that $\xi_0^N = \bar{\xi}_0^N = \hat{\xi}_0^N$ and $\tilde{\xi}_0^N = \tilde{\bar{\xi}}_0^N = \tilde{\hat{\xi}}_0^N$ so that the construction of these particle systems via (SDE)($I'$) as in [3] ensures that

(9.17) $\quad\quad\quad \hat{\xi}_t^N \leq \bar{\xi}_t^N, \quad\quad \xi_t^N \leq \bar{\xi}_t^N \quad\quad \text{for all } t \geq 0,$

as well as

(9.18) $\quad\quad\quad \tilde{\hat{\xi}}_t^N \leq \tilde{\bar{\xi}}_t^N, \quad\quad \tilde{\xi}_t^N \leq \tilde{\bar{\xi}}_t^N \quad\quad \text{for all } t \geq 0.$

We use $E_{\xi_0^N, \tilde{\xi}_0^N}$ to denote expectations for our initial conditions $\tilde{\xi}_0^N \leq \xi_0^N$ as above.

PROOF OF LEMMA 9.1. In this proof constants $C$ and functions $C(\cdot)$ will depend only on $\bar{\theta}$ and $p(\cdot)$, and may change from line to line. The $C(T)$ will always be assumed to be an increasing function from $\mathbb{R}_+$ to $\mathbb{R}_+$. Implicit use of Remark 4.6 will be made to ensure this from time to time.

Let $\tilde{P}_t^N$ denote the semigroup of $\beta_t^{N,0}$, killed when it exits $I'$. Arguing just as in the proof of Lemma 3.2 of [3], we get

(9.19)
$$\tilde{X}_t^N(\mathbf{1}) = \tilde{X}_0^N(\tilde{P}_t^N \mathbf{1}) + \tilde{M}_t^N$$
$$+ \log N \theta_0^N \int_0^t \frac{1}{N'} \sum_x \tilde{P}_{t-s}^N \mathbf{1}(x)(1 - \tilde{\xi}_s^N(x)) f_1^N(x, \tilde{\xi}_s^N)^2 \, ds$$
$$- \log N \theta_1^N \int_0^t \frac{1}{N'} \sum_x \tilde{P}_{t-s}^N \mathbf{1}(x) \tilde{\xi}_s^N(x) f_0^N(x, \tilde{\xi}_s^N)^2 \, ds,$$

where $\tilde{M}_t^N$ is a square integrable martingale with mean 0. Since $\tilde{X}_0^N = X_0^N$,

we may take differences with (3.7) (with $\phi_s = 1$) to conclude that

$$E(X_t^N(\mathbf{1}) - \tilde{X}_t^N(\mathbf{1}))$$
$$= X_0^N(\mathbf{1} - \tilde{P}_t^N \mathbf{1})$$
$$+ \int_0^t E\bigg(\theta_0^N \log N \frac{1}{N'} \sum_x (1 - \xi_s^N(x)) f_1^N(x, \xi_s^N)^2 (1 - \tilde{P}_{t-s}^N \mathbf{1}(x))$$
$$- \theta_1^N \log N \frac{1}{N'} \sum_x \xi_s^N(x) f_0^N(x, \xi_s^N)^2 (1 - \tilde{P}_{t-s}^N \mathbf{1}(x)) \bigg) ds$$

(9.20) $\quad + \theta_0^N \log N \int_0^t E\bigg(\frac{1}{N'} \sum_x \tilde{P}_{t-s}^N \mathbf{1}(x)$
$$\times [(1 - \xi_s^N(x)) f_1^N(x, \xi_s^N)^2$$
$$- (1 - \tilde{\xi}_s^N(x)) f_1^N(x, \tilde{\xi}_s^N)^2] \bigg) ds$$
$$- \theta_1^N \log N \int_0^t E\bigg(\frac{1}{N'} \sum_x \tilde{P}_{t-s}^N \mathbf{1}(x) [\xi_s^N(x) f_0^N(x, \xi_s^N)^2$$
$$- \tilde{\xi}_s^N(x) f_0^N(x, \tilde{\xi}_s^N)^2] \bigg) ds.$$

Choose $\phi_{K,L} : \mathbb{R}_+ \to \mathbb{R}_+$ so that

$$1(x \geq KL) \leq \phi_{K,L}(x) \leq 1(x \geq (K-1)L) \quad \text{and} \quad |\phi_{K,L}|_{\text{Lip}} \leq 1.$$

Then

(9.21) $\quad \mathbf{1} - \tilde{P}_{t-s}^N \mathbf{1}(x) \leq E\bigg(\phi_{K,L}\bigg(\sup_{u \leq t-s} \|x + B_u^{N,0}\|\bigg)\bigg) \equiv \bar{h}_{t-s}(x)$

and

(9.22) $\quad\quad\quad\quad |\bar{h}_{t-s}|_{\text{Lip}} \leq 1, \quad \|\bar{h}_{t-s}\|_{\text{Lip}} \leq 2.$

For $s \leq t$ and $\tilde{\xi} \leq \xi$ are in $\{0,1\}^{\mathbf{S}_N}$, define

$$H_0(t-s, \xi) = \frac{\log N}{N'} \sum_x \bar{h}_{t-s}(x)[(1 - \xi(x)) f_1^N(x, \xi)^2 + \xi(x) f_0^N(x, \xi)^2]$$

$$H_1(\xi, \tilde{\xi}) = \frac{\log N}{N'} \sum_x |\xi(x) f_0^N(x, \xi)^2 - \tilde{\xi}(x) f_0^N(x, \tilde{\xi})^2|$$

$$H_2(\xi, \tilde{\xi}) = \frac{\log N}{N'} \sum_x |(1 - \xi(x)) f_1^N(x, \xi)^2 - (1 - \tilde{\xi}(x)) f_1^N(x, \tilde{\xi})^2|.$$

From (9.20) and (9.21) we get

$$E(X_t^N(\mathbf{1}) - \tilde{X}_t^N(\mathbf{1}))$$
$$(9.23) \qquad \leq X_0^N(\mathbf{1} - \tilde{P}_t^N \mathbf{1})$$
$$+ \bar{\theta}\left[\int_0^t E(H_0(t-s, \xi_s^N))\, ds + \int_0^t \sum_{i=1}^2 H_i(\xi_s^N, \tilde{\xi}_s^N)\, ds\right].$$

Our goal is to show, by estimating each of the above terms, that there are constants $c_0(t)$ and $c_1$ such that

$$E(X_t^N(\mathbf{1}) - \tilde{X}_t^N(\mathbf{1}))$$
$$(9.24) \qquad \leq c_0(t) X_0^N(\mathbf{1})\left[\frac{KL}{\log N} + P\left(\sup_{u \leq t(1+\delta_N')} \|\beta_u^{N,0}\| > (K-2)L\right)\right]$$
$$+ c_1 \int_0^t E(X_s^N(\mathbf{1}) - \tilde{X}_s^N(\mathbf{1}))\, ds.$$

Applying Gronwall's lemma, we obtain (9.3).

*Step* 1. The first term in (9.23) is simple. Since $X_0^N$ is supported on $I = [-L, L]^2$,

$$(9.25) \qquad X_0^N(\mathbf{1} - \tilde{P}_t^N \mathbf{1}) \leq X_0^N(\mathbf{1}) P\left(\sup_{u \leq t} \|\beta_u^{N,0}\| \geq (K-1)L\right).$$

*Step* 2. Let $\varepsilon = (\log N)^{-p}$, where $p = 18$, and consider the $H_0$ term in (9.23). (The choice $p = 18$ is used only in the last line of this proof.) We first note that

$$\frac{N'}{\log N} H_0(t-s, \xi)$$
$$\leq \sum_x \sum_y [(1-\xi(x))\xi(y) + \xi(x)(1-\xi(y))] p_N(y-x) \bar{h}_{t-s}(x)$$
$$= 2\xi(f_0^N(\xi)\bar{h}_{t-s})$$
$$+ \sum_x \sum_y [(1-\xi(x))\xi(y) p_N(x-y)(\bar{h}_{t-s}(x) - \bar{h}_{t-s}(y))]$$
$$\leq 2\xi(f_0^N(\xi)\bar{h}_{t-s}) + \sqrt{2}\sigma|\xi|N^{-1/2},$$

using (9.22) and the covariance assumption on $p$ to bound $\sum_z |z| p_N(z)$ by $\sqrt{2}\sigma N^{-1/2}$ in the last line. Now we may use the Markov property and the

above to see that

$$E\left(\int_0^t H_0(t-s,\xi_s^N)\,ds\right)$$
$$\leq E\left(\int_0^\varepsilon H_0(t-s,\xi_s^N)\,ds\right)$$
$$+ E\left(\int_\varepsilon^{t\vee\varepsilon} 2E_{X_{s-\varepsilon}^N}(X_\varepsilon^N(\log N f_0^N(\xi_\varepsilon^N))\bar{h}_{t-s})\,ds\right)$$
$$+ \sqrt{2}\sigma\frac{\log N}{\sqrt{N}}\int_\varepsilon^{t\vee\varepsilon} E(X_s^N(\mathbf{1}))\,ds$$
$$\leq 2\log N E\left(\int_0^\varepsilon X_s^N(\mathbf{1})\,ds\right) + \frac{2C_{4.8}}{\log N}E\left(\int_\varepsilon^{t\vee\varepsilon} X_{s-\varepsilon}^N(\mathbf{1})\,ds\right)$$
$$+ 2C_{4.8}E\left(\int_\varepsilon^{t\vee\varepsilon} X_{s-\varepsilon}^N(\bar{h}_{t-s})\,ds\right) + \sqrt{2}\sigma\frac{\log N}{\sqrt{N}}\int_\varepsilon^{t\vee\varepsilon} E(X_s^N(\mathbf{1}))\,ds,$$

where we have used the trivial bound $H_0(t-s,\xi_s^N) \leq 2\log N X_s^N(\mathbf{1})$, Proposition 4.5 and (9.22) in the last line. Recall the definition of $P_t^{N,*}$ before Proposition 4.4. Next use (4.2), Proposition 4.4, (9.22) and a bit of arithmetic to bound the above by

$$C(t)\left[\frac{X_0^N(\mathbf{1})}{\log N} + \int_\varepsilon^{t\vee\varepsilon} X_0^N(P_{s-\varepsilon}^{N,*}\bar{h}_{t-s})\,ds\right]$$

for some $C(t)$. In view of (9.1), $P_t^{N,*} = P_{t(1+\delta_N')}^N$, where $P_t^N$ is the semigroup of $\beta_t^{N,0}$, and we readily see that

$$P_{s-\varepsilon}^{N,*}\bar{h}_{t-s}(x) \leq P\left(\sup_{u\leq t(1+\delta_N')} \|\beta_u^{N,x}\| > (K-1)L\right).$$

Now use the fact that $X_0^N$ is supported on $I$ to conclude that

(9.26)
$$E\left(\int_0^t H_0(t-s,\xi_s^N)\,ds\right)$$
$$\leq C_{9.26}(t)X_0^N(\mathbf{1})\left[(\log N)^{-1} + P\left(\sup_{u\leq t(1+\delta_N')} \|\beta_u^{N,0}\| > (K-2)L\right)\right]$$

for some $C_{9.26}(t)$.

*Step* 3. Turn next to the $H_i$ terms in (9.23), $i=1,2$. With an application of the Markov property in mind, let us for the moment consider more general initial conditions than that in the theorem and assume only $\tilde{\xi}_0^N \leq \xi_0^N$ are such that $|\xi_0^N| < \infty$ and $\tilde{\xi}_0^N$ is supported on $I'$. We couple $\tilde{\xi}_t^N \leq \xi_t^N$ and $\tilde{\bar{\xi}}_t^N \leq \bar{\xi}_t^N$ as described above so that (9.17) and (9.18) hold. Let $q = 1/6$ and for $\delta > 0$, let

$$I'(\delta) = \{w \in I' : d(w, \partial I') \leq \delta\},$$

where $d(w, \partial I')$ is the distance from $w$ to the boundary of $I'$ in the supremum norm on $\mathbb{R}^2$. The goal of this step is to prove that there is a constant $C_{9.27}$ such that

$$\begin{aligned}
E_{\xi_0^N, \tilde{\xi}_0^N}\left(\sum_{i=1}^2 H_i(\xi_\varepsilon^N, \tilde{\xi}_\varepsilon^N)\right) \\
\leq C_{9.27}\left[X_0^N(\mathbf{1}) - \tilde{X}_0^N(\mathbf{1}) + \frac{X_0^N(\mathbf{1})}{\log N} + \log N \tilde{X}_0^N(I'(2\varepsilon^q))\right].
\end{aligned}$$
(9.27)

As before, we proceed via comparisons with the biased voter model and voter model, and hence, need the decomposition

(9.28) $\quad E_{\xi_0^N, \tilde{\xi}_0^N}(H_i(\xi_\varepsilon^N, \tilde{\xi}_\varepsilon^N)) \leq \Gamma_1 + \Gamma_2 + E_{\xi_0^N, \tilde{\xi}_0^N}(H_i(\hat{\xi}_\varepsilon^N, \tilde{\hat{\xi}}_\varepsilon^N)),$

where

$$\Gamma_1 = \Gamma_1^i = E_{\xi_0^N, \tilde{\xi}_0^N}(|H_i(\xi_\varepsilon^N, \tilde{\xi}_\varepsilon^N) - H_i(\bar{\xi}_\varepsilon^N, \tilde{\bar{\xi}}_\varepsilon^N)|),$$

$$\Gamma_2 = \Gamma_2^i = E_{\xi_0^N, \tilde{\xi}_0^N}(|H_i(\bar{\xi}_\varepsilon^N, \tilde{\bar{\xi}}_\varepsilon^N) - H_i(\hat{\xi}_\varepsilon^N, \tilde{\hat{\xi}}_\varepsilon^N)|).$$

After expanding $\Gamma_1$ using the definitions of $H_1$ and $H_2$, and rearranging using the inequality $||a - b| - |c - d|| \leq |a - c| + |b - d|$, we have

$$\Gamma_1 \leq E_{\xi_0^N, \tilde{\xi}_0^N}(H_i(\bar{\xi}_\varepsilon^N, \xi_\varepsilon^N)) + E_{\xi_0^N, \tilde{\xi}_0^N}(H_i(\tilde{\bar{\xi}}_\varepsilon^N, \tilde{\xi}_\varepsilon^N)).$$

We estimate these two terms as follows.

By (5.33) and (5.34), we have for $i = 1, 2$,

$$\begin{aligned}
E_{\xi_0^N}(H_i(\bar{\xi}_\varepsilon^N, \xi_\varepsilon^N)) &\leq 3 \log N E(\bar{X}_\varepsilon^N(\mathbf{1}) - X_\varepsilon^N(\mathbf{1})) \\
&\leq 3C_{6.7} \log N[(\log N)^{-2} + \varepsilon]X_0^N(\mathbf{1}) \\
&\leq C_{9.29} X_0^N(\mathbf{1})/\log N
\end{aligned}$$
(9.29)

for a constant $C_{9.29}$, where in the next to last line we used the basic coupling and (6.7).

Again, as in the proof of Lemma 3.2 in [3], and arguing as in (5.36), we may show that, for some square integrable martingale $\tilde{M}_t^N$,

$$\tilde{X}_t^N(\mathbf{1}) = \tilde{X}_0^N(\tilde{P}_t^{N,*}\mathbf{1}) + \tilde{M}_t^N + \int_0^t 2\bar{\theta} \log N \tilde{X}_s^N(\tilde{P}_{t-s}^{N,*}\mathbf{1} f_0^N(\tilde{\xi}_s^N))\, ds,$$

where $\tilde{P}_t^{N,*} = \tilde{P}_{t(1+\delta_N')}^N$ and $\delta_N'$ is as in (9.2). Take differences with (9.19) and use the fact that $\tilde{P}_t^N \mathbf{1}$ is decreasing in $t$ to see that

$$E(\tilde{\bar{X}}_\varepsilon^N(\mathbf{1}) - \tilde{X}_\varepsilon^N(\mathbf{1})) \leq \tilde{X}_0^N(\tilde{P}_{\varepsilon(1+\delta_N')}^N \mathbf{1} - \tilde{P}_\varepsilon^N \mathbf{1})$$

$$+ 2\bar{\theta} \log N \int_0^\varepsilon E(\tilde{\bar{X}}_s^N(\mathbf{1}) + \tilde{X}_s^N(\mathbf{1}))\,ds$$

(9.30)

$$\leq 2\bar{\theta} \log N \int_0^\varepsilon E(\tilde{\bar{X}}_s^N(\mathbf{1}) + \tilde{X}_s^N(\mathbf{1}))\,ds$$

$$\leq C(\log N)^{-2} X_0^N(\mathbf{1}),$$

where (6.3) and the coupling (9.18) is used in the last line. Now repeat the derivation of (9.29) to see that, for $i = 1, 2$,

(9.31) $$E_{\tilde{\xi}_0^N}(H_i(\tilde{\bar{\xi}}_\varepsilon^N, \tilde{\xi}_\varepsilon^N)) \leq 3C(\log N)^{-1} X_0^N(\mathbf{1}).$$

We have therefore proved [(9.29) and (9.31)] there is a constant $C_{9.32}$ such that

(9.32) $$\Gamma_1 \leq C_{9.32}(\log N)^{-1} X_0^N(\mathbf{1}), \qquad i = 1, 2.$$

Similar reasoning leads to

(9.33) $$\Gamma_2 \leq C_{9.33}(\log N)^{-1} X_0^N(\mathbf{1}), \qquad i = 1, 2.$$

To finish the proof of (9.27), it remains only to prove the voter model inequality

(9.34)
$$E_{\xi_0^N, \tilde{\xi}_0^N}\left(\sum_{i=1}^2 H_i(\hat{\xi}_\varepsilon^N, \tilde{\hat{\xi}}_\varepsilon^N)\right) \leq C_{9.34}\left(X_0^N(\mathbf{1}) - \tilde{X}_0^N(\mathbf{1}) + \frac{\tilde{X}_0^N(\mathbf{1})}{\log N}\right)$$
$$+ 6 \log N \tilde{X}_0^N(I'(2\varepsilon^q))$$

for a constant $C_{9.34}$. We need an elementary lemma.

LEMMA 9.3. *If* $\tilde{\xi} \leq \xi \in \{0,1\}^{\mathbf{S}_N}$, $\bar{H}_1(\xi, \tilde{\xi}) = \frac{\log N}{N'} \sum_x (\xi(x) - \tilde{\xi}(x)) f_0^N(x, \xi)$ *and* $\bar{H}_2(\xi, \tilde{\xi}) = \frac{\log N}{N'} \sum_x \tilde{\xi}(x)[f_0^N(x, \tilde{\xi}) - f_0^N(x, \xi)]$, *then*

$$H_1(\xi, \tilde{\xi}) + H_2(\xi, \tilde{\xi}) \leq 3(\bar{H}_1(\xi, \tilde{\xi}) + \bar{H}_2(\xi, \tilde{\xi})).$$

PROOF.

$$H_2(\xi, \tilde{\xi}) \leq \frac{\log N}{N'} \sum_x (1 - \xi(x))[f_1(x, \xi)^2 - f_1(x, \tilde{\xi})^2] + (\xi(x) - \tilde{\xi}(x)) f_1(x, \tilde{\xi})^2$$

$$\leq \frac{2 \log N}{N'} \sum_x \sum_y (1 - \xi(x)) p_N(y - x)(\xi(y) - \tilde{\xi}(y))$$

$$+ \frac{\log N}{N'} \sum_x \sum_y (\xi(x) - \tilde{\xi}(x)) p_N(y - x) \tilde{\xi}(y)$$

TWO-DIMENSIONAL LOTKA–VOLTERRA MODEL 63$$= 2\bar{H}_1(\xi,\tilde{\xi}) + \frac{\log N}{N'} \sum_y \tilde{\xi}(y) \sum_x p_N(y-x)(1-\tilde{\xi}(x)) - (1-\xi(x)))$$

$$= (2\bar{H}_1 + \bar{H}_2)(\xi,\tilde{\xi}).$$

Simpler reasoning shows that $H_1(\xi,\tilde{\xi}) \le (\bar{H}_1 + 2\bar{H}_2)(\xi,\tilde{\xi})$. □

To prove (9.34), we will also need an extension of the voter model duality described in Sections 5 and 7. Define the killed walks $\tilde{B}_t^{N,x}$ and $\tilde{\hat{B}}_t^{N,x}$ by setting them equal to the rescaled walks $B_t^{N,x}$ and $\hat{B}_t^{N,x}$ of (7.13), but killed (set equal to a cemetery state $\Delta$) when they first leave $I'$ at time

$$\hat{\tau}_N(x,I') = \inf\{t \ge 0 : \hat{B}_t^{N,x} \notin I'\}.$$

We also set

(9.35) $$\tilde{\xi}_0^N(\Delta) = 0.$$

With this convention, the joint duality we need is as follows: for all $e_i, f_i \in \{0,1\}$ and $x_i, y_i \in \mathbf{S}_N$,

(9.36)
$$P_{\xi_0^N, \tilde{\xi}_0^N}(\hat{\xi}_\varepsilon^N(x_i) = e_i, \tilde{\hat{\xi}}_\varepsilon^N(y_i) = f_i, i=1,\ldots,M)$$
$$= P(\xi_0^N(\hat{B}_\varepsilon^{N,x_i}) = e_i, \tilde{\xi}_0^N(\tilde{\hat{B}}_\varepsilon^{N,y_i}) = f_i, i=1,\ldots,M).$$

This is readily obtained from the coupling of $\hat{\xi}^N$ and $\tilde{\hat{\xi}}^N$ through the stochastic differential equation in Proposition 2.1 of [3]. One can refine the dynamics there by using appropriately defined uniformly distributed random variables to identify the "parent" of a 0 or 1, enabling us to define the usual Poisson arrows, which in turn allow us to define the coalescing dual. If the dual random walk from a site $x$ leads to a site outside $I'$, we know that $\tilde{\xi}_\varepsilon^N(x) = 0$ by our 0 boundary conditions and so by freezing the dual random walk at $\Delta$, we ensure the validity of (9.36) thanks to (9.35). The details are standard and left for the reader.

As before, $e$ denote a random variable with distribution $p_N$, independent of our coalescing random walks. Let $I' - x$ be the obvious translation of $I'$. By (9.36), we have

$$E_{\xi_0^N,\tilde{\xi}_0^N}(\bar{H}_2(\hat{\xi}_\varepsilon^N, \tilde{\hat{\xi}}_\varepsilon^N))$$
$$= \frac{\log N}{N'} \sum_x P(\tilde{\hat{\xi}}_\varepsilon(x) = 1, \tilde{\hat{\xi}}_\varepsilon(x+e) = 0, \hat{\xi}_\varepsilon(x+e) = 1)$$
$$= \frac{\log N}{N'} \sum_x P(\tilde{\xi}_0^N(\tilde{\hat{B}}_\varepsilon^{N,x}) = 1, \tilde{\xi}_0^N(\tilde{\hat{B}}_\varepsilon^{N,x+e}) = 0, \xi_0^N(\hat{B}_\varepsilon^{N,x+e}) = 1)$$

$$= \frac{\log N}{N'} \sum_x [P(\tilde{\xi}_0^N(\hat{\tilde{B}}_\varepsilon^{N,x}) = 1, \tilde{\xi}_0^N(\hat{B}_\varepsilon^{N,x+e}) = 0, \xi_0^N(\hat{B}_\varepsilon^{N,x+e}) = 1)$$

$$+ P(\tilde{\xi}_0^N(\hat{\tilde{B}}_\varepsilon^{N,x}) = 1, \tilde{\xi}_0^N(\hat{B}_\varepsilon^{N,x+e}) = 1, \tilde{\xi}_0^N(\hat{\tilde{B}}_\varepsilon^{N,x+e}) = 0)],$$

since $\tilde{\xi}_0^N(\hat{B}_\varepsilon^{N,x+e}) = 0$ implies $\tilde{\xi}_0^N(\hat{\tilde{B}}_\varepsilon^{N,x+e}) = 0$, and $\tilde{\xi}_0^N(\hat{B}_\varepsilon^{N,x+e}) = 1$ implies $\xi_0^N(\hat{B}_\varepsilon^{N,x+e}) = 1$. Using translation invariance, and the fact that $\tilde{\xi}_0^N(\hat{\tilde{B}}_\varepsilon^{N,x}) = 1$ and $\tilde{\xi}_0^N(\hat{B}_\varepsilon^{N,x+e}) = 0$ imply $\hat{B}_\cdot^{N,x}$ and $\hat{B}_\cdot^{N,x+e}$ have not coalesced before $\varepsilon$, the above is not larger than

$$\frac{\log N}{N'} \sum_w (\xi_0^N(w) - \tilde{\xi}_0^N(w)) \sum_x P(\hat{B}_\varepsilon^{N,e} = w - x, \hat{\tau}^N(0,e) > \varepsilon)$$

$$+ \frac{\log N}{N'} \sum_w \tilde{\xi}_0^N(w) \sum_x P(\hat{B}_\varepsilon^{N,e} = w - x, \hat{\tau}^N(e, I' - x) \leq \varepsilon)$$

$$\equiv S_1 + S_2.$$

Summing over $x$ first, and using (2.4), we conclude that

(9.37)
$$S_1 = \log N(X_0^N(\mathbf{1}) - \tilde{X}_0^N(\mathbf{1}))H(N\varepsilon)$$
$$\leq C_{9.37}(X_0^N(\mathbf{1}) - \tilde{X}_0^N(\mathbf{1}))$$

for a constant $C_{9.37}$. In the last, recall from Step 2 that $\varepsilon = (\log N)^{-18}$.

For $S_2$, we first separate out the sums (1) $w \in I'(2\varepsilon^q)$ and (2) $w \in I' \setminus I'(2\varepsilon^q)$ and $\|w - x\| > \varepsilon^q$. Note that $d(x, \partial I') > \varepsilon^q$ for the remaining $x, w$. Therefore,

$$S_2 \leq \log N \tilde{X}_0^N(I'(2\varepsilon^q)) + \log N \tilde{X}_0^N(\mathbf{1}) P(\|\hat{B}_\varepsilon^{N,e}\| > \varepsilon^q)$$

$$+ \frac{\log N}{N'} \sum_w \sum_x 1(w \in I' \setminus I'(2\varepsilon^q), \|x - w\| \leq \varepsilon^q) \tilde{\xi}_0^N(w)$$

$$\times P\left(\hat{B}_\varepsilon^{N,e} = w - x, \sup_{u \leq \varepsilon} \|\hat{B}_u^{N,e}\| > \varepsilon^q\right)$$

(9.38)
$$\leq \log N \tilde{X}_0^N(I'(2\varepsilon^q)) + \log N \tilde{X}_0^N(\mathbf{1}) E(\|\hat{B}_\varepsilon^{N,e}\|^2) \varepsilon^{-2q}$$

$$+ \log N \tilde{X}_0^N(\mathbf{1}) P\left(\sup_{u \leq \varepsilon} \|\hat{B}_u^{N,e}\| > \varepsilon^q\right)$$

$$\leq \log N \tilde{X}_0^N(I'(2\varepsilon^q)) + C_{9.38} \log N \tilde{X}_0^N(\mathbf{1}) \varepsilon^{1-2q}$$

for a constant $C_{9.38}$, where we have used Doob's weak maximal inequality in the last. Combine (9.37) and (9.38), and use $p(1 - 2q) > 2$ to derive

(9.39)
$$E_{\xi_0^N, \tilde{\xi}_0^N}(\bar{H}_2(\hat{\xi}_\varepsilon^N, \hat{\tilde{\xi}}_\varepsilon^N))$$
$$\leq \log N \tilde{X}_0^N(I'(2\varepsilon^q)) + C_{9.39}[(\log N)^{-1} \tilde{X}_0^N(\mathbf{1}) + X_0^N(\mathbf{1}) - \tilde{X}_0^N(\mathbf{1})]$$

for a constant $C_{9.39}$.

For $\bar{H}_1$, we have, by similar reasoning,

$$E_{\xi_0^N, \tilde{\xi}_0^N}(\bar{H}_1(\hat{\xi}_\varepsilon^N, \tilde{\hat{\xi}}_\varepsilon^N))$$

$$= \frac{\log N}{N'} \sum_x P(\xi_0^N(\hat{B}_\varepsilon^{N,x}) = 1, \tilde{\xi}_0^N(\hat{\tilde{B}}_\varepsilon^{N,x}) = 0, \xi_0^N(\hat{B}_\varepsilon^{N,x+e}) = 0)$$

$$= \frac{\log N}{N'} \sum_x [P(\xi_0^N(\hat{B}_\varepsilon^{N,x}) = 1, \tilde{\xi}_0^N(\hat{B}_\varepsilon^{N,x}) = 0, \xi_0^N(\hat{B}_\varepsilon^{N,x+e}) = 0)$$

$$+ P(\tilde{\xi}_0^N(\hat{B}_\varepsilon^{N,x}) = 1, \tilde{\xi}_0^N(\hat{\tilde{B}}_\varepsilon^{N,x}) = 0, \xi_0^N(\hat{B}_\varepsilon^{N,x+e}) = 0)]$$

$$\leq \frac{\log N}{N'} \sum_w (\xi_0^N(w) - \tilde{\xi}_0^N(w)) \sum_x P(\hat{B}_\varepsilon^{N,0} = w - x, \hat{\tau}^N(0,e) > \varepsilon)$$

$$+ \frac{\log N}{N'} \sum_w \tilde{\xi}_0^N(w) \sum_x P(\hat{B}_\varepsilon^{N,0} = w - x, \hat{\tau}^N(0, I' - x) \leq \varepsilon)$$

$$\equiv S_1' + S_2'.$$

The $S_i'$ are actually slightly simpler than $S_i$ to handle (we have 0 in place of $e$) and so, as before, we may bound $E_{\xi_0^N, \tilde{\xi}_0^N}(\bar{H}_1(\hat{\xi}_\varepsilon^N, \tilde{\hat{\xi}}_\varepsilon^N))$ by the right-hand side of (9.39). Combine these two bounds and use Lemma 9.3 to complete the derivation of (9.34), and hence of (9.27).

After inserting bounds (9.25), (9.26) and (9.27) into (9.23), the Markov property implies that

$$E(X_t^N(\mathbf{1}) - \tilde{X}_t^N(\mathbf{1}))$$

$$\leq X_0^N(\mathbf{1}) P\left(\sup_{u \leq t} \|\beta_u^{N,0}\| > (K-1)L\right)$$

$$+ \bar{\theta} C_{9.26}(t) X_0^N(\mathbf{1})$$

$$\times \left[(\log N)^{-1} + P\left(\sup_{u \leq t(1+\delta_N')} \|\beta_u^{N,0}\| > (K-2)L\right)\right]$$

$$+ \bar{\theta} E\left(\int_0^\varepsilon (H_1 + H_2)(\xi_s^N, \tilde{\xi}_s^N) \, ds\right)$$

$$+ \bar{\theta} C_{9.27} E\left(\int_\varepsilon^{t \vee \varepsilon} E\left(X_{s-\varepsilon}^N(\mathbf{1}) - \tilde{X}_{s-\varepsilon}^N(\mathbf{1}) + \frac{X_{s-\varepsilon}^N(\mathbf{1})}{\log N} \right.\right.$$

$$\left.\left. + \log N \tilde{X}_{s-\varepsilon}^N(I'(2\varepsilon^q)) \right) ds\right).$$

Thus, in order to prove (9.24), it suffices now to prove the two bounds

$$(9.40) \quad \int_0^\varepsilon E((H_1 + H_2)(\xi_s^N, \tilde{\xi}_s^N))\, ds + \frac{1}{\log N} \int_0^t E(X_s^N(\mathbf{1}))\, ds$$
$$\leq \frac{C_{9.40}(t)}{\log N} X_0^N(\mathbf{1})$$

and

$$(9.41) \quad \log N \int_0^t E(\tilde{X}_s^N(I'(2\varepsilon^q)))\, ds \leq C_{9.41}(t)(KL/\log N) X_0^N(\mathbf{1}).$$

The bound (9.40) follows easily from (4.2), the definition of $\varepsilon$ and the fact that $H_i(\xi_s^N, \tilde{\xi}_s^N) \leq 2 \log N X_s^N(\mathbf{1})$. For (9.41), let $\bar{I}'(\varepsilon) = \{w : d(w, \partial I') \leq 4\varepsilon^q\}$ and choose $\psi_\varepsilon : \mathbb{R}^2 \to [0,1]$ such that

$$1_{I'(2\varepsilon^q)} \leq \psi_\varepsilon \leq 1_{\bar{I}'(\varepsilon)} \quad \text{and} \quad \|\psi_\varepsilon\|_{\text{Lip}} \leq \varepsilon^{-q}.$$

Then, using Proposition 4.4, the left-hand side of (9.41) is bounded above by

$$\log N E\left(\int_0^t X_s^N(\psi_\varepsilon)\, ds\right)$$
$$\leq c(t)\left[\varepsilon^{-q}(\log N)^{(3-p)/2} X_0^N(\mathbf{1}) + \log N \int_0^t X_0^N(P_s^{N,*}\psi_\varepsilon)\, ds\right]$$
$$\leq c(t)\left[(\log N)^{-9/2} X_0^N(\mathbf{1}) + \log N \int_0^t \int P(B_s^{N,*,x} \in \bar{I}'(\varepsilon)) X_0^N(dx)\, ds\right],$$

for some $c(t)$, where $B_s^{N,*,x}$ is a random walk starting at $x$ with semigroup $P_t^{N,*}$ as in Proposition 4.4. Use the bound on $P(B_t^{N,*,x} = w)$ from (7.30) to see that

$$\log N \int_0^t \int P(B_s^{N,*,x} \in \bar{I}'(\varepsilon)) X_0^N(dx)\, ds$$
$$\leq c X_0^N(\mathbf{1}) \log N \int_0^t \frac{\varepsilon^q N K L}{1 + Ns}\, ds$$
$$= c X_0^N(\mathbf{1}) K L (\log N)^{1 - p/6} \log(1 + Nt)$$
$$\leq c(t) K L X_0^N(\mathbf{1}) (\log N)^{-1}.$$

We have finally used our choice of $p = 18$. This proves (9.41), and hence, completes the proof of (9.24). □

**Acknowledgment.** We thank an anonymous referee for a thorough reading and many useful suggestions.

# REFERENCES

[1] Cox, J. T., Durrett, R. and Perkins, E. A. (2000). Rescaled voter models converge to super-Brownian motion. *Ann. Probab.* **28** 185–234. MR1756003
[2] Cox, J. T. and Perkins, E. A. (2005). Rescaled Lotka–Volterra models converge to super-Brownian motion. *Ann. Probab.* **33** 904–947. MR2135308
[3] Cox, J. T. and Perkins, E. A. (2007). Survival and coexistence in stochastic spatial Lotka–Volterra models. *Probab. Theory Related Fields* **139** 89–142. MR2322693
[4] Ethier, S. N. and Kurtz, T. G. (1986). *Markov Processes. Characterization and Convergence.* Wiley, New York. MR0838085
[5] Griffeath, D. (1978). *Additive and Cancellative Interacting Particle Systems.* Springer, Berlin. MR0538077
[6] Liggett, T. M. (1985). *Interacting Particle Systems.* Springer, New York. MR0776231
[7] Neuhauser, C. and Pacala, S. W. (1999). An explicitly spatial version of the Lotka–Volterra model with interspecific competition. *Ann. Appl. Probab.* **9** 1226–1259. MR1728561
[8] Perkins, E. (2002). Dawson–Watanabe superprocesses and measure-valued difffusions. *École d'Été de Probabilités de Saint Flour XXIX–1999. Lecture Notes Math.* **1781** 125–324. Springer, Berlin. MR1915445
[9] Presutti, E. and Spohn, H. (1983). Hydrodynamics of the voter model. *Ann. Probab.* **11** 867–875. MR0714951
[10] Spitzer, F. L. (1976). *Principles of Random Walk*, 2nd ed. Springer, New York. MR0388547

Department of Mathematics
Syracuse University
Syracuse, New York 13244
USA
E-mail: jtcox@syr.edu

Department of Mathematics
University of British Columbia
1984 Mathematics Road
Vancouver, British Columbia
Canada V6T 1Z2
E-mail: perkins@math.ubc.ca